\documentclass[12pt]{amsart}
\date{22 May 2005}
\usepackage{amsmath,latexsym,amsfonts,amssymb,amscd,amsxtra,euscript}

\theoremstyle{plain}
\begingroup
\newtheorem{thm}{Theorem}[section]
\newtheorem{cor}[thm]{Corollary}
\newtheorem{lem}[thm]{Lemma}
\newtheorem{prop}[thm]{Proposition}

\endgroup

\theoremstyle{definition}
\newtheorem{defn}[thm]{Definition}
\newtheorem{assumption}[thm]{Assumption}

\theoremstyle{remark}
\newtheorem{rem}[thm]{Remark}

\numberwithin{equation}{section}
\numberwithin{table}{section}

\DeclareMathOperator*{\Coeff}{Coeff}

\newcommand{\bm}{\bar{m}}
\newcommand{\bdd}{\bar{d}_0}
\newcommand{\bdm}{\bar{d}_M}
\newcommand{\bd}{\bar{d}_1}

\newcommand{\ima}{\mathbf{i}}

\newcommand{\abs}[1]{\lvert #1\rvert}
\renewcommand{\setminus}{\smallsetminus}

\newcommand{\dbar}{\bar{\partial}}
\newcommand{\vol}{\mathrm{vol}}
\newcommand{\Tr}{\mathrm{Tr}\,}
\newcommand{\norm}[1]{\lVert#1\rVert}
\newcommand{\rk}{\operatorname{rk}}
\newcommand{\pdeg}{\operatorname{pardeg}}
\newcommand{\pmu}{\operatorname{par\mu}}
\newcommand{\Hom}{\operatorname{Hom}}
\newcommand{\PH}{\operatorname{ParHom}}
\newcommand{\SPH}{\operatorname{SParHom}}
\newcommand{\PE}{\operatorname{ParEnd}}
\newcommand{\SPE}{\operatorname{SParEnd}}

\newcommand{\surj}{\twoheadrightarrow}
\newcommand{\inc}{\hookrightarrow}
\newcommand{\ar}{\rightarrow}
\newcommand{\too}{\longrightarrow}
\newcommand{\x}{\times}
\newcommand{\ox}{\otimes}
\newcommand{\iso}{\cong}

\newcommand{\End}{\mathrm{End}}

\newcommand{\Jac}{\mathrm{Jac}}

\newcommand{\Ad}{\mathrm{Ad}}

\newcommand{\Sym}{\mathrm{Sym}}

\newcommand{\Ext}{\mathrm{Ext}}

\newcommand{\Tor}{\mathrm{Tor}}

\newcommand{\Id}{\mathrm{Id}}

\newcommand{\Sl}{\mathrm{SL}}

\newcommand{\cE}{\mathcal{E}}

\newcommand{\cM}{\mathcal{M}}

\newcommand{\cN}{\mathcal{N}}

\newcommand{\cO}{\mathcal{O}}

\newcommand{\cS}{\mathcal{S}}
\newcommand{\cT}{\mathcal{T}}

\newcommand{\CC}{\mathbb{C}}

\newcommand{\HH}{\mathbb{H}}
\newcommand{\PP}{\mathbb{P}}
\newcommand{\QQ}{\mathbb{Q}}
\newcommand{\RR}{\mathbb{R}}

\newcommand{\ZZ}{\mathbb{Z}}
\newcommand{\WW}{\mathbb{W}}
\newcommand{\AAA}{{\curly A}}
\newcommand{\CCC}{{\curly C}}

\newcommand{\GGG}{{\curly G}}
\newcommand{\HHH}{{\curly H}}

\DeclareFontFamily{OT1}{rsfs}{}
\DeclareFontShape{OT1}{rsfs}{n}{it}{<->rsfs10}{}
\DeclareMathAlphabet{\curly}{OT1}{rsfs}{n}{it}

\renewcommand{\a}{\alpha}
\renewcommand{\b}{\beta}

\newcommand{\e}{\varepsilon}

\newcommand{\s}{\sigma}

\newcommand{\p}{\phi}

\title[Betti numbers of the moduli space of rank 3 parabolic Higgs
bundles]{Betti numbers of the moduli space of rank 3 parabolic
Higgs bundles}

\author{O. Garc\'{\i}a-Prada}
   \address{Departamento de Matem\'aticas \\
  CSIC \\ Serrano 113 bis
  \\ 28006 Madrid \\ Spain}
   \email{oscar.garcia-prada@uam.es}
\author{P.B. Gothen}
\address{Departamento de Matematica Pura \\
  Faculdade de Ci\^encias da Universidade do Porto \\
  Rua do Campo Alegre 687 \\ 4169-007 Porto \\ Portugal }
   \email{pbgothen@fc.up.pt}
\author{V. Mu\~noz}
  \address{Departamento de Matem\'aticas \\
  CSIC \\ Serrano 113 bis
  \\ 28006 Madrid \\ Spain}
  \email{vicente.munoz@imaff.cfmac.csic.es}

  \thanks{Partially supported by Ministerio de Educaci\'{o}n y
    Tecnolog\'{\i}a and Conselho de Reitores das Universidades
    Portuguesas through Acci\'{o}n Integrada Hispano-Lusa HP-2000-0015
    (Spain) / E--30/03 (Portugal) and by the European Research
    Training Networks EDGE (Contract no.\ HPRN-CT-2000-00101)
    and EAGER (Contract no.\ HPRN-CT-2000-00099).  Members of
    VBAC (Vector Bundles on Algebraic Curves). Second author partially supported
    by Centro de Matem\'atica da Universidade do Porto, financed by
    FCT (Portugal) through the programmes POCTI and POSI of the QCA
    III (2000--2006) with European Community (FEDER) and national
    funds. First and third authors partially supported by
    Ministerio de Educaci{\'o}n y Ciencia (Spain) through Project
    BFM2000-0024.
}

\begin{document}



\begin{abstract}
  Parabolic Higgs bundles on a Riemann surface are of interest for
  many reasons, one of them being their importance in the study of
  representations of the fundamental group of the punctured surface in
  the complex general linear group. In this paper we calculate the
  Betti numbers of the moduli space of rank 3 parabolic Higgs bundles
with fixed and non-fixed determinant,
  using Morse theory.  A key point is that certain critical
  submanifolds of the Morse function can be identified with moduli
  spaces of parabolic triples.  These moduli spaces come in families
  depending on a real parameter and we carry out a careful analysis of
  them by studying their variation with this parameter.  Thus we
  obtain in particular information about the topology of the moduli
  spaces of parabolic triples for the value of the parameter relevant
  to the study of parabolic Higgs bundles.  The remaining critical
  submanifolds are also described: one of them is the moduli space of
  parabolic bundles, while the remaining ones have a description in
  terms of symmetric products of the Riemann surface.  As another
  consequence of our Morse theoretic analysis, we obtain a proof of
  the parabolic version of a theorem of Laumon, which states that the
  nilpotent cone (the preimage of zero under the Hitchin map) is a
  Lagrangian subvariety of the moduli space of parabolic Higgs
  bundles.
\end{abstract}

\maketitle

\section{Introduction}
\label{sec:introduction}

Let $X$ be a connected, smooth projective complex algebraic curve
of genus $g$ and let $D=p_1+ p_2 +\cdots +p_n$ be a divisor, with
distinct points $p_1,\ldots, p_n$. Let $K$ be the canonical bundle
of $X$.  A {\em parabolic Higgs bundle} is a pair $(E,\Phi)$,
where $E$ is a parabolic bundle, that is a holomorphic bundle over
$X$ together with a weighted flag in the fibre of $E$ over each
$p\in D$, and $\Phi: E\to E\otimes K(D)$ is a strongly parabolic
homomorphism.  This means that $\Phi$ is a meromorphic
endomorphism valued one-form with simple poles along $D$ whose
residue at $p$ is nilpotent with respect to the flag.

Like in the non-parabolic case, there is a stability criterion
allowing the construction of moduli spaces of semistable parabolic
Higgs bundles \cite{Y1}. For generic weights, semistability and
stability coincide and the moduli space is a smooth quasiprojective
algebraic manifold.  The goal of this paper is to compute the Betti
numbers of this moduli space in the case in which the rank of the
bundle is 3.  The computation for rank $2$ was carried out by Boden
and Yokogawa \cite{BY}, and Nasatyr and Steer \cite{NS} in the case of
rational weights. In the non-parabolic case the Betti numbers had been
previously computed by Hitchin \cite{H} in rank $2$ and Gothen
\cite{G1} in rank $3$.  We have been informed by T. Hausel of a
conjecture for the Betti numbers of the moduli space of parabolic
Higgs bundles of any rank and degree, analogous to his conjecture for
the case of non-parabolic Higgs bundles \cite{Ha2}. His formula gives
the same result as ours in the cases that we have checked (cf.\
Remark~\ref{rem:hausel}), thus providing  support for
his conjecture.

Similarly to the non-parabolic case, the moduli space of parabolic
Higgs bundles has an extremely rich geometric structure. It can be
identified with the moduli space of solutions to the parabolic
version of Hitchin's equations:
 $$
 F(A)^{\perp} + [\Phi,\Phi^*] = 0\;\;\; \mbox{and}\;\;\; \dbar_A\Phi=0,
 $$
where $A$ is  a singular connection unitary with respect to a
singular hermitian metric on $E$  adapted to the parabolic
structure (see Section \ref{sec:moduli} for details). The moduli
space of parabolic Higgs bundles  contains the total space of the
cotangent bundle of the moduli space of parabolic  bundles, whose
natural holomorphic symplectic form can be extended  to the whole
moduli space. This form can be combined with the  real symplectic
form coming from the gauge-theoretic interpretation to endow the
moduli space with a hyperk\"ahler structure \cite{K,NS}.

An important motivation to study ordinary Higgs bundles comes from
their relation with  complex representations of the fundamental
group of the curve. This is established by identifying the moduli
space of solutions to Hitchin's equations with  the moduli space of
Higgs bundles \cite{H,S2} as well as with the moduli space of
complex connections with constant central curvature \cite{C,D}. In
the parabolic case, there is a similar correspondence proved by
Simpson \cite{S1}. This involves meromorphic complex connections
with simple poles at the points  and parabolic weights. At the
topological side one has to consider filtered local systems. The
natural context for the correspondence is a class larger than
parabolic  Higgs bundles in which the Higgs field $\Phi$ is allowed
to be parabolic and not necessarily strongly parabolic. In other
words, at a parabolic point the residue of $\Phi$ is parabolic with respect
to the flag.  Under this correspondence, parabolic Higgs bundles
(those for which the Higgs field is strongly parabolic) are in
bijection with  meromorphic flat connections whose holonomy
around each parabolic point defines a conjugacy class of an  element
 in the unitary group. These, in turn, correspond to representations of the
fundamental group of the punctured surface in the general linear
group, which send a small loop around each  parabolic point to an
element conjugate to a unitary element.

The main tool for our computation of the Betti numbers, as in the
previously studied cases, is the use of the Morse-theoretic
techniques introduced by Hitchin \cite{H}:  the $L^2$-norm of the
Higgs field defines a perfect Bott--Morse function on the moduli
space. We have to compute the Poincar\'e series and the indices of
the various critical subvarieties.  In fact the calculation of the
indices can be carried out for any rank, whereas the calculation of
the Poincar\'e series of the critical subvarieties depends crucially
on the rank $3$ assumption. Here is a description of the paper.

In Section \ref{sec:higgs} we review the basic definitions and basic
facts of parabolic Higgs bundles. In Section \ref{sec:morse} we
consider the Bott--Morse function on the moduli space and identify the
critical subvarieties.  These coincide with the fixed subvarieties
under the action of $S^1$ on the moduli space given by multiplying the
Higgs field. These in turn correspond, as shown by Simpson \cite{S1},
to variations of Hodge structures, in particular the bundle has to be
a direct sum of subbundles.  We then compute the indices --- this can
be done for any rank and leads to a parabolic version of the theorem
of Laumon, that the nilpotent cone in the moduli space of parabolic
Higgs bundles is a Lagrangian subvariety. In the rank $3$ case, the
possible decompositions of the vector bundle in a sum of subbundles
are of two types: a sum of three line bundles or a sum of a line
bundle and a rank $2$ vector bundle. The latter case gives rise to
so-called parabolic triples. These have been introduced in \cite{BiG}
and generalise the triples studied in \cite{BG} and \cite{BGG}.
Through Sections \ref{sec:stable-triples},
\ref{sec:critical-values-flips} and \ref{sec:thaddeus-program} we
study the moduli spaces of parabolic triples.  They depend on a real
parameter, relating to parabolic Higgs bundles when the value of this
parameter is $2g-2$. To compute the Betti numbers for a given value of
the parameter (in particular for $2g-2$) we follow the strategy
introduced by Thaddeus in \cite{T1}. After characterising the moduli
space for the largest value of the parameter, we need to analyse the
changes when we cross a finite number of values until we get to the
one we want.  In Section \ref{sec:(1,1,1)} we compute the Poincar\'e
polynomial and indices for the critical subvarieties for which the
vector bundle is a sum of three line bundles. In Sections
\ref{sec:(1,2)} and \ref{sec:(2,1)} we do the other cases using the
previous computations for the moduli space of parabolic triples.  An
important technical point is that the Betti numbers of the moduli
space of parabolic Higgs bundles do not depend on the degree and the
weights (certainly if the weights are generic).  So we can choose the
degree coprime with the rank and the weights as small as convenient,
to facilitate our computations. In Section \ref{par-rank3}, based on
the computations by Nitsure \cite{Ni} and Holla \cite{Ho} of the Betti
numbers of the moduli space of parabolic bundles, we work out the
formula for the rank 3 case. In Section \ref{sec:rank3} we collect all
the computations, to give the Poincar\'e polynomial of the rank $3$
moduli space of parabolic Higgs bundles. Finally, in Section
\ref{sec:fixed-det}, we compute the Poincar\'e polynomial of the rank
$3$ moduli space of parabolic Higgs bundles with fixed determinant. It
is interesting to observe that, like in the non-parabolic case, and in
contrast to the case of stable parabolic bundles, the Poincar\'e
polynomial of the non-fixed determinant moduli space does not split as
the product of those of the Jacobian and the fixed determinant moduli
space. In particular, it follows that tensoring by a line bundle gives
a non-trivial action of the group of elements of order three in the
Jacobian on the cohomology of the fixed determinant moduli space with
rational coefficients; in fact our methods allow us to determine
precisely the non-invariant part of the rational cohomology.


\textbf{Acknowledgements.} We wish to thank Tam\'as Hausel and Nigel
Hitchin for very useful comments.


\section{Parabolic Higgs bundles}
\label{sec:higgs}

\subsection{Definitions and basic facts}
\label{sec:parabolic-definitions}

Let $X$ be a connected, smooth projective complex algebraic  curve
of genus $g$ together with a finite (non-zero) number of {\em
marked\/} distinct points $p_1,\ldots, p_n$. We will denote the
divisor $D=p_1+ p_2 +\cdots +p_n$.

Let $E$ be a holomorphic bundle over $X$. A parabolic structure on
$E$ consists of weighted flags
 $$
 E_p=E_{p,1} \supset E_{p,2} \supset \cdots \supset E_{p,s_p}\supset E_{p,s_p+1}=0,
 $$
 $$
 0\leq \a_1(p)< \cdots <\a_{s_p}(p) < 1,
 $$
over each $p\in D$. A holomorphic map $\p: E\to F$ between
parabolic bundles is called parabolic if $\a_i^E(p) >\a_j^F(p)$
(the superindex indicating to which parabolic bundle the weight
corresponds) implies $\p(E_{p,i}) \subset F_{p,j+1}$ for all
$p\in D$. We call $\p$ strongly parabolic if $\a_i^E(p)\geq
\a_j^F(p)$ implies $\p(E_{p,i}) \subset F_{p, j+1}$ for all $p\in
D$.

We will abuse notation by simply writing  $E$ for a bundle with a
parabolic structure.

This notion  can be generalised to  higher dimensions. In this
case, a parabolic bundle  consists of a holomorphic vector  bundle
together with a weighted holomorphic filtration of the
restriction of the  bundle to a fixed  divisor. This
generalisation is also relevant for us in the particular case in
which the manifold is the product of the curve $X$ and a higher
dimensional manifold $Y$, and  the divisor is $D\times Y$.
Parabolic and strongly parabolic homomorphisms are defined in a
similar way.

Also $\PH(E,F)$ and $\SPH(E,F)$ will denote respectively the
sheaves  of  parabolic and strongly parabolic homomorphisms from $E$
to $F$.

Let $m_i(p)=\dim E_{p,i}/E_{p,i+1}$ be the multiplicity of
$\a_i(p)$. It will sometimes be convenient to repeat each weight
according to its multiplicity, i.e., we set $\tilde{\alpha}_1(p) =
\dotso = \tilde{\alpha}_{m_1(p)}(p) = \alpha_1(p)$, etc.  We then
have weights $0 \leq\tilde{\alpha}_1(p) \leq \cdots \leq
\tilde{\alpha}_r(p) < 1$, where $r= \rk (E)$ is the rank of $E$.
Define the parabolic degree and parabolic slope of $E$ by
$$
  \pdeg (E)= \deg( E )+ \sum_{p\in D} \sum_{i=1}^{s_p}
  m_i(p)\a_i(p)
  = \deg (E) + \sum_{p\in D} \sum_{i=1}^{r}
  \tilde{\alpha}_i(p),
$$
$$
  \pmu (E) =\frac{\pdeg (E)}{\rk (E)}.
$$

If $F$ is a subbundle of $E$, then $F$ inherits a parabolic structure
by setting $F_{p,i} =F_p \cap E_{p,i}$ and discarding those weights
of multiplicity zero. We call this the induced parabolic structure on
$F$. In a similar manner, one can give a parabolic structure to the
quotient $E/F$.

A parabolic bundle $E$ is said to be stable if
$\pmu(F)<\pmu(E)$ for all proper parabolic subbundles
$F\subset E$. Semistability is defined by replacing the strict
inequality by the weak inequality. For generic weights, stability
and semistability are equivalent.

For parabolic bundles $E$ there is a well-defined notion of dual
$E^*$. This is done by considering the bundle $\Hom(E,
{\mathcal O}(-D))$, and at each $p\in D$, defining the filtration
 $$
 E^*_p=E^*_{p,1} \supset \cdots \supset E^*_{p,s_p} \supset 0,
 $$
with $E^*_{p,i}=\Hom (E_p/ E_{p,s_p+2-i} , {\mathcal O}(-D)_p)$
and weights $1-\a_s(p) < \cdots <1-\a_1(p)$. It is easy to prove
that $E^{**}=E$, and  $\pdeg (E^*)= -\pdeg (E)$.

There is also a notion of tensor product $\otimes^P$ of two
parabolic bundles \cite{Y2}, which is best understood in terms of
$\RR$-filtered sheaves. Here we shall only use the case of
tensoring a parabolic bundle $E$ with a parabolic line bundle $L$.
Let $\alpha_i(p)$ be the weights of $E$ and $\beta(p)$ be the
weights of $L$. Then the parabolic bundle $F=E\otimes^P L$ is, as
a bundle, the kernel of
  $$
   E\otimes L(D) \surj \oplus_{p\in D} \left( (E_p/E_{p,i_p}) \otimes
   L(D)_p\right) ,
  $$
where $i_p=\min \{ s_p +1, i \, |\, \alpha_i(p) +\beta(p) \geq
1\}$, $p\in D$. The weights of $F$ are
  \begin{equation}\label{eq:weights-tensor}
  \alpha_{i_p}(p) +\beta(p)-1 < \cdots <\alpha_{s_p}(p) +\beta(p)-1
  < \alpha_{1}(p) +\beta(p) < \cdots <\alpha_{i_p-1}(p)
  +\beta(p)\, ,
   \end{equation}
with multiplicities $m_{i_p}(p), \ldots, m_{s_p}(p),m_1(p),\ldots,
m_{i_p-1}(p)$. It is then easy to see that $\pmu (E\otimes^P
L)=\pmu(E) +\pmu(L)$, whereas
 $$ 
  \deg(E\otimes^P L) = \deg(E) + \rk(E) \deg(L) + \sum_p \dim
  E_{p,i_p}\, .
 $$ 

We denote by $K$ the canonical bundle on $X$. A {\em parabolic
Higgs bundle}  is a pair $(E,\Phi)$, where $E$ is a parabolic
bundle and $\Phi \in H^0(\SPE(E)\otimes K(D))$, i.e.\ $\Phi$ is a
meromorphic endomorphism valued one-form with simple poles along
$D$ whose residue at $p$ is nilpotent with respect to the flag.
We shall sometimes denote a parabolic Higgs bundle by $\mathbf{E}
=(E,\Phi)$.

The notion of stability is extended to parabolic Higgs bundles in
the usual way:
  $$
  \pmu(F)<\pmu(E)
  $$
for all proper parabolic subbundles $F\subset E$ which are
preserved by $\Phi$. Semistability is defined by replacing the
strict inequality by the weak inequality.

The standard properties of stable bundles also apply to parabolic
Higgs bundles; for example, if $\mathbf{E}$ and $\mathbf{F}$ are
stable parabolic Higgs bundles of the same parabolic slope, then
there are no parabolic maps between them unless they are
isomorphic, and the only parabolic endomorphisms of a stable
parabolic Higgs bundle are the scalar multiples of the identity.

We shall say that the weights are {\em generic\/} when every
semistable parabolic Higgs bundle is automatically stable, i.e.,
when there are no properly semistable Higgs bundles. Let us fix
(generic) weights $\alpha_i(p)$ and topological invariants $\rk
(E)$ and $\deg (E)$.  The moduli space $\mathcal{M}$ of stable
parabolic Higgs bundles was constructed using Geometric Invariant
Theory by Yokogawa \cite{Y1,Y2}, who also showed that it is a
smooth irreducible complex variety. The moduli space $\mathcal{M}$
contains the cotangent bundle of the moduli space of stable
parabolic bundles.

The following result will facilitate the computation of the Betti
numbers of $\mathcal{M}$.

\begin{prop} \label{prop:deg-tensor}
  Fix the rank $r$. For different choices of degrees and generic
  weights, the moduli spaces of parabolic Higgs bundles have the same
  Betti numbers.
\end{prop}

\begin{proof}
  For fixed degree, it is a consequence of the results of Thaddeus
  \cite{T2} that the moduli spaces for different generic weights have
  the same Betti numbers, as we now explain.  The space of weights is
  divided into chambers by a finite number of hyperplanes, or
  \emph{walls}, and in each chamber the moduli spaces are isomorphic.
  Call the moduli spaces on each side of a wall $\mathcal{M}^+$ and
  $\mathcal{M}^-$, respectively (here, and in the following, we use
  the notation of \cite{T2}).  Thaddeus proves that
  $\mathcal{M}^+$ and $\mathcal{M}^-$ have a common blow-up with the
  same exceptional divisor.  The loci in $\mathcal{M}^\pm$ to be blown
  up (\emph{flip loci}, in our language) are isomorphic to projective
  bundles $\PP U^\pm$ over a product $\mathcal{N}^+ \times
  \mathcal{N}^-$ of moduli spaces of lower rank parabolic Higgs
  bundles.  This is similar\footnote{In these cases it is the
    stability parameter which is varied, rather than the parabolic
    weights.} to the situation in \cite{T1} (for ordinary triples) and
  our analysis in Section~\ref{sec:critical-values-flips} below (for
  parabolic triples), and shows that the difference between the
  Poincar\'e polynomials of $\mathcal{M}^+$ and $\mathcal{M}^-$ equals
  the difference between the Poincar\'e polynomials of the
  respective flip loci (cf.\ \cite[p.\ 605]{GH}).  However, in the
  case of parabolic Higgs bundles something special happens, namely
  the bundles $U^+$ and $U^-$ are dual to each other (see the
  paragraph of \cite{T2} preceding (5.6)). Hence $\PP U^+$ and $\PP
  U^-$ are projective bundles of the same rank, over the same base.
  But the Poincar\'e polynomial of a projective bundle splits as the
  product of the Poincar\'e polynomial of the base and the Poincar\'e
  polynomial of projective space. Thus the flip loci in
  $\mathcal{M}^+$ and $\mathcal{M}^-$ have the same Betti numbers and,
  therefore, the moduli spaces $\mathcal{M}^+$ and $\mathcal{M}^-$
  themselves have the same Betti numbers.

  To extend the result to moduli spaces of parabolic Higgs bundles
  with different degrees, we proceed as follows. Fix any parabolic
  line bundle $L$ with degree $d_L$ and weights $\beta(p)$. Then the
  map
  $$
   (E,\Phi) \mapsto (E\otimes^P L,\Phi)
   $$
   gives an {\em isomorphism\/} between the moduli space of parabolic Higgs
   bundles of rank $r$, degree $\Delta$ and weights $\alpha_i(p)$ and
   the moduli space of parabolic Higgs bundles of rank $r$, degree $\Delta + r
   d_L + \sum_p \sum_{i\geq i_p} m_i(p)$, and weights given by
   (\ref{eq:weights-tensor}). Choosing weights of multiplicity one and
   a suitable parabolic line bundle $L$ we see that moduli spaces of
   parabolic Higgs bundles for different degrees are isomorphic. Since
   we already know that the Betti numbers are independent of the
   (generic) weights for fixed degree, this concludes the proof.
\end{proof}

\subsection{Deformation theory}
\label{sec:deformation}

The deformation theory of parabolic Higgs bundles was worked out by
Yokogawa \cite{Y2}; see also Thaddeus \cite{T2} and Biswas and
Ramanan \cite{BR}.  Everything in this section is essentially
contained in these references but we shall find it convenient to give
an exposition tailored to our purposes.  Let $\mathbf{E} = (E,\Phi)$
and $\mathbf{F} = (F,\Psi)$ be parabolic Higgs bundles.  We define a
complex of sheaves
\begin{align*}
  C^\bullet(\mathbf{E},\mathbf{F}) : \quad
  \PH(E,F) &\to \SPH(E,F) \otimes K(D) \\
  f &\mapsto (f \otimes 1)\Phi - \Psi f,
\end{align*}
and write $C^\bullet(\mathbf{E}) = C^\bullet(\mathbf{E},\mathbf{E})$.
\begin{prop}\label{prop:deformation}
  \begin{enumerate}
  \item[(i)] The space of infinitesimal deformations of a parabolic Higgs
    bundle $\mathbf{E}$ is naturally  isomorphic to the first
    hypercohomology group of the complex
    \begin{align*}
      C^\bullet (\mathbf{E}) : \quad
      \PE(E) &\overset{[-,\Phi]}{\longrightarrow} \SPE(E) \otimes K(D) \\
      f &\longmapsto (f \otimes 1)\Phi - \Phi f.
    \end{align*}
    Thus the tangent space
    to $\mathcal{M}         $ at a point represented by a stable
    parabolic Higgs bundle $\mathbf{E}$ is isomorphic to
    $\mathbb{H}^1(C^\bullet(\mathbf{E}))$.
  \item[(ii)] The space of homomorphisms between parabolic Higgs bundles
    $\mathbf{E}$ and $\mathbf{F}$ is naturally isomorphic to the zeroth
    hypercohomology group
    $\mathbb{H}^0(C^\bullet(\mathbf{E},\mathbf{F}))$.
  \item[(iii)] The space of extensions $0 \to
    \mathbf{E}' \to \mathbf{E} \to \mathbf{E}'' \to 0$ of parabolic
    Higgs bundles $\mathbf{E}'$ and $ \mathbf{E}''$ is naturally
    isomorphic to the first hypercohomology
    $\mathbb{H}^1(C^\bullet(\mathbf{E}'',\mathbf{E}'))$.
  \item[(iv)] There is a long exact sequence
    \begin{equation}\label{eq:les-phiggs}
      \begin{split}
      0 &\to \mathbb{H}^0(C^\bullet(\mathbf{E},\mathbf{F})) \to
      H^0(\PH(E,F)) \to H^0(\SPH(E,F) \otimes K(D)) \\
      &\to
      \mathbb{H}^1(C^\bullet(\mathbf{E},\mathbf{F}))
      \to
      H^1(\PH(E,F)) \to H^1(\SPH(E,F) \otimes K(D)) \\
      &\to
      \mathbb{H}^2(C^\bullet(\mathbf{E},\mathbf{F})) \to 0.
      \end{split}
    \end{equation}
  \end{enumerate}
\end{prop}

\begin{proof}
  For proofs of (i) -- (iii) see Thaddeus \cite{T2}.  The proof of
  (iv) follows by a standard argument in Higgs bundle theory (see,
  e.g., Biswas and Ramanan \cite{BR}).
\end{proof}

As for ordinary Higgs bundles, duality plays an important role for
parabolic Higgs bundles.  The results of the following proposition are
consequences of the theory developed by Yokogawa (cf.\ (3.1) and
Proposition~3.7 of \cite{Y2}, see also \S~3 of Thaddeus \cite{T2} and
\S~5 of Bottacin \cite{Bo}).
\begin{prop}\label{prop:par-duality}
  \begin{enumerate}
  \item[(i)] Let $E$ and $F$ be parabolic bundles.  The sheaves $\PH(E,F)$
  and $\SPH(F,E(D))$ are naturally dual.
  \item[(ii)] Let $\mathbf{E}$ and $\mathbf{F}$ be parabolic Higgs bundles.
  Then there is a natural isomorphism
  \begin{displaymath}
    \mathbb{H}^{i}(C^\bullet(\mathbf{E},\mathbf{F})) \cong
    \mathbb{H}^{2-i}(C^\bullet(\mathbf{F},\mathbf{E}))^*.
  \end{displaymath}
  In particular we obtain a natural isomorphism
  $T_{\mathbf{E}}\mathcal{M} \cong T^*_{\mathbf{E}}\mathcal{M}$ for a
  stable parabolic Higgs bundle $\mathbf{E}$.
  \end{enumerate}
\end{prop}

\begin{proof}
  We just show how (i) implies (ii).  {}From (i) it follows that the dual
  complex $C^\bullet(\mathbf{E},\mathbf{F})^*$ is related to the
  original one by
\begin{displaymath}
  C^\bullet(\mathbf{E},\mathbf{F})^* \otimes K \cong
  C^\bullet(\mathbf{F},\mathbf{E}).
\end{displaymath}
Thus, Serre duality for hypercohomology and
Proposition~\ref{prop:deformation} (i) give statement (ii) of the
present proposition.
\end{proof}

Next we shall show how these results can be used to calculate the
dimension of the moduli space of parabolic Higgs bundles.  The
result is well known but it seems worthwhile to include the
calculation here, since we shall use the same ideas again below.
First we introduce some convenient notation for parabolic bundles
$E$ and $F$ as follows.  We denote by $P_p(E,F)$ the subspace of
$\Hom(E_{p},F_{p})$ consisting of parabolic maps and by $N_p(E,F)$
the subspace of strictly parabolic maps.  We also write $P_D(E,F)
= \bigoplus P_p(E,F)$ and $N_D(E,F) = \bigoplus N_p(E,F)$.  When
there is no risk of confusion we shall omit the parabolic bundles
$E$ and $F$ from the notation.  We then have short exact sequences
of sheaves
\begin{displaymath}
  0 \to \PH(E,F) \to \Hom(E,F) \to \Hom(E_D,F_D)/P_D(E,F) \to
  0,
\end{displaymath}
and
\begin{displaymath}
  0 \to \SPH(E,F) \to \Hom(E,F) \to \Hom(E_D,F_D)/N_D(E,F) \to
  0.
\end{displaymath}
Thus we can calculate the Euler characteristics of $\PH(E,F)$ and
$\SPH(E,F)$ as follows:
  \begin{equation}\label{eq:chi-phEF}
  \begin{split}
  \chi(\PH(E,F)) &= \chi(\Hom(E,F)) +  \sum_{p\in D}(\dim P_p -
  \rk(E)\rk(F)), \\
  \chi(\SPH(E,F)) &= \chi(\Hom(E,F)) + \sum_{p\in D}(\dim N_p -
  \rk(E)\rk(F)).
 \end{split}
 \end{equation}
With these preliminaries in place we can calculate the dimension of
the moduli space.

\begin{prop}\label{prop:phb-moduli-dim}
  The complex dimension of the moduli space $\mathcal{M}$ of
  stable rank $r$ parabolic Higgs bundles is
  \begin{displaymath}
    r^2(2g-2) + 2 + 2\sum_{p} f_p,
  \end{displaymath}
  where $r=\rk (E)$ and $f_p = \frac{1}{2}\big(r^2 - \sum_{i}m_{i}(p)^2\big)$.
\end{prop}

\begin{proof}
  Since $\mathbf{E}$ is stable, its only endomorphisms are the
  scalars.  Hence, using Proposition~\ref{prop:deformation} (ii)
  and the duality statement Proposition~\ref{prop:par-duality} (ii),
  we have that $\dim \mathbb{H}^0(C^\bullet(\mathbf{E}))
  = \dim
  \mathbb{H}^2(C^\bullet(\mathbf{E})) = 1$.  It follows that the dimension
  of the moduli space is
  \begin{align*}
    \dim \mathcal{M}
    &= \dim \mathbb{H}^1(C^\bullet(\mathbf{E})) \\
    &= 2 - \chi(C^\bullet(\mathbf{E}))  \\
    &= 2 - \chi(\PE(E)) + \chi(\SPE(E) \otimes K(D)),
  \end{align*}
  where in the last equality we have used the long exact sequence
  \eqref{eq:les-phiggs}.   {}From this we obtain the result by using
  equations \eqref{eq:chi-phEF},  
  the fact
  that $\dim P_p - \dim N_p = \sum_i m_i(p)^2$ and the Riemann--Roch
  formula.
\end{proof}

\subsection{Parabolic Higgs bundles and gauge theory}
\label{sec:moduli}

Our main goal is to study the topology of $\mathcal{M}$. To do
this we  need the gauge-theoretic interpretation of this moduli
space in terms of solutions to Hitchin's equations due to Simpson
\cite{S1}. The construction of the moduli space from this point of
view is due to Konno \cite{K}.  Let $E$ be a smooth parabolic
vector bundle of rank $r$ and fix a hermitian metric $h$ in $E$
which is smooth in $X\setminus D$ and whose (degenerate) behaviour
around the punctures is given as follows.  We say that a local
frame $\{e_1,\dotsc,e_r\}$ for $E$ around $p$ \emph{respects the
flag at $p$} if $E_{p,i}$ is spanned by the vectors
$\{e_{M_i+1}(p),\dotsc,e_r(p)\}$, where $M_i = \sum_{j \leq i}
m_j$. Let $z$ be a local coordinate around $p$ such that $z(p) =
0$.  We require that $h$ be of the form
\begin{displaymath}
  h =
  \begin{pmatrix}
    \abs{z}^{2\tilde{\alpha}_1} & & 0 \\
     & \ddots & \\
    0 & & \abs{z}^{2\tilde{\alpha}_r} \\
  \end{pmatrix}
\end{displaymath}
with respect to some local frame around $p$ which respects the
flag at $p$.  We denote the space of smooth $\dbar$-operators on
$E$ by $\CCC$ and the space of associated $h$-unitary connections
by $\AAA$.  Note that the unitary connection associated to a
smooth $\dbar_A$ via the hermitian metric $h$ is singular around
the punctures: if we write $z = \rho \exp(\ima \theta)$ and
$\{e_i\}$ is the local frame used in the definition of $h$, then
with respect to the local frame $\{\epsilon_i =
e_i/\abs{z}^{\tilde{\alpha}_i}\}$, the connection is of the form
 $$
  d_A = d + \ima \left(
  \begin{smallmatrix}
    \tilde{\alpha}_1 && 0 \\
    & \ddots & \\
    0 && \tilde{\alpha}_r
  \end{smallmatrix}
  \right)
  d\theta + A',
 $$
where  $A'$ is regular.

We denote the space of Higgs fields by $\boldsymbol{\Omega} =
\Omega^{1,0}(\SPE(E)\otimes\cO(D))$, the group of complex parabolic gauge
transformations by $\GGG_{\CC}$ and the subgroup of $h$-unitary
parabolic gauge transformations by $\GGG$.

Following Biquard \cite{Bi}, Konno introduces certain weighted Sobolev
norms; we denote the corresponding Sobolev completions of the spaces
defined above by ${\CCC}^{p}_{1}$, $\boldsymbol{\Omega}^p_1$,
${(\GGG_{\CC})}^p_2$ and ${\GGG}^p_2$ (the detailed
definitions are not important to us so we refer to \cite{K} for
them).  Let
  $$
  {\HHH}= \{ (\dbar_A,\Phi) \in {\CCC} \times {\boldsymbol{\Omega}}
  \ | \  \dbar_A \Phi = 0 \}
  $$
and let ${\HHH}^p_1$ be the corresponding subspace of ${\CCC}^p_1
\times {\boldsymbol{\Omega}}^p_1$.  Then ${\HHH}^p_1$ carries a
hyper-K\"ahler metric induced by $h$. Let $F(A)^{\perp}$ denote
the trace-free part of the curvature of the $h$-unitary connection
$A$ corresponding to $\dbar_A$ and let $\Phi^*$ be the adjoint
with respect to $h$.  One can then consider the moduli space $\cS$
defined by the subspace  of ${\HHH}^p_1$ satisfying
\emph{Hitchin's equation} (modulo ${\GGG}^p_2$),
\begin{displaymath}
  \cS =\{(\dbar_A,\Phi) \in {\HHH}_1^p \ | \  F(A)^{\perp} + [\Phi,\Phi^*] = 0
  \}/{\GGG}^p_2,
\end{displaymath}
where the equation is only defined on $X\setminus D$. Konno proves
that $\cS$ is a  hyper-K\"ahler quotient and that it can be
naturally identified with the moduli space
 $$
 \mathcal{M} = \HHH^p_1 / {(\GGG_{\CC})}^p_2\ .
 $$
Furthermore, Konno proves that the natural map $\HHH / \GGG_{\CC}
\to \HHH^p_1 / {(\GGG_{\CC})}^p_2$ is a diffeomorphism.

\section{Morse theory on the moduli space}
\label{sec:morse}

\subsection{The Morse function}
\label{sec:morse-function}

The non-zero complex numbers $\CC^*$ act on the moduli space
$\mathcal{M}$ via the map $\lambda\cdot(E,\Phi) = (E, \lambda
\Phi)$.  However, to have an action on the set of solutions to
Hitchin's equations, one must restrict to the action of $S^1
\subset \CC^*$.  Obviously the identification $\mathcal{S} \cong
\mathcal{M}$ respects the circle action and thus we have a circle
action on this hyper-K\"ahler manifold.  With respect to one of
the complex structures (coinciding with the one on $\mathcal{M}$)
this is a Hamiltonian action and the associated moment map is
\begin{equation}\label{eq:aaa3}
  [(A,\Phi)] \mapsto -\tfrac{1}{2}\norm{\Phi}^2
  = -\ima \int_X \Tr (\Phi\Phi^*).
\end{equation}
We shall, however, prefer to consider the positive function
\begin{equation}\label{eq:aaa4}
  f([A,\Phi]) = \tfrac{1}{2}\norm{\Phi}^2.
\end{equation}
In the case of non-parabolic Higgs bundles, Hitchin \cite{H}
proved that this is a proper map, using Uhlenbeck's compactness
theorem.  It was observed by Boden and Yokogawa \cite{BY} that the
same argument works in the parabolic case, by using the parabolic
analogue of Uhlenbeck's theorem, proved by Biquard \cite{Bi}. Thus
we have the following result.
\begin{prop}
  The map $f \colon \mathcal{M} \to \RR$ is proper. \qed
\end{prop}

Next we recall a general result of Frankel \cite{F}, which was first
used in the context of moduli spaces of Higgs bundles by Hitchin
\cite{H}.

\begin{thm} \label{thm:3.2}
  Let $\tilde{f}\colon M \to \RR$ be a proper moment map for a
  Hamiltonian circle action on a K\"ahler manifold $M$.  Then
  $\tilde{f}$ is a perfect Bott--Morse function.  \qed
\end{thm}

The following result on the Morse indices of such a Morse function is
implicit in Frankel's paper.

\begin{prop}\label{prop:critical=fixed}
  In the situation of Theorem \ref{thm:3.2}, the critical points of
  $\tilde{f}$ are exactly the fixed points of the circle action.
  Moreover, the eigenvalue $l$ subspace for the Hessian of $\tilde{f}$
  is the same as the weight $-l$ subspace for the infinitesimal circle
  action on the tangent space.
  In particular, the Morse index of $\tilde{f}$ at a critical point equals
  the dimension of the positive weight space of the circle action
  on the tangent space.
\end{prop}

\begin{proof}
The condition for $\tilde{f}$ to be a moment map is that
\begin{displaymath}
  \mathrm{grad}(\tilde{f}) = IX,
\end{displaymath}
where $X$ is the vector field generating the circle action and $I$ is
the complex structure on $M$.  Hence $p$ is a critical point of $\tilde{f}$
if and only if it is fixed under the circle action.  Let $\nabla$ be
the Levi-Civita connection on $M$, then the Hessian $H_{\tilde{f}}$ of
$\tilde{f}$ at $p$ is the quadratic form associated to the symmetric
endomorphism $\nabla(\mathrm{grad}(\tilde{f}))_p$ of $T_pM$.  Let
$Y_{p} \in T_{p}M$ and let the vector field $Y$ be an extension of
$Y_p$ around $p$.  Then we have
\begin{align*}
  H_{\tilde{f}}(Y_p) &= \nabla_{Y_p}(IX) \\
  &= \nabla_{IX_p}(Y) - [IX,Y]_p \\
  &= -[IX,Y]_p,
\end{align*}
where we have used that $X_p = 0$.  On the other hand it is easy to
see (cf.\ \cite{F}) that the infinitesimal circle action on $T_pM$ is given
by $Y_p \mapsto [Y,X]_p$.  It follows that the eigenvalues of
$H_{\tilde{f}}$ are exactly minus the weights of the circle action on $T_pM$.
\end{proof}

Thus we must identify the fixed point set of the action of $S^1
\subset \CC^*$ on $\mathcal{M}$. This was done by Simpson and is
analogous to what happens for ordinary Higgs bundles.

\subsection{Fixed points of the $S^1$ action on the moduli space}
\label{sec:fixed}

\begin{prop}[{\cite[Theorem 8]{S1}}]\label{prop:fixed=vhs}
  The equivalence class of a stable parabolic Higgs bundle
  $(E,\Phi)$ is fixed under the action of $S^1$ if
and only if it is
  a \emph{parabolic complex variation of Hodge structure}.  This means
  that $E$ has a direct sum decomposition
  \begin{displaymath}
    E = E_0 \oplus \cdots \oplus E_m
  \end{displaymath}
  as parabolic bundles, such that $\Phi$ is strongly parabolic and of
  degree one with respect to this decomposition, in other words the
  restriction $\Phi_l = \Phi_{\vert E_l}$ belongs to
  $$
  H^0 (\SPH(E_l,E_{l+1}) \otimes K(D)).
  $$
Furthermore, stability implies that $\Phi_l \neq 0$ for $l = 0,
\dotsc, m-1$. The \emph{type} of the parabolic complex variation
of Hodge structure is the vector $(\rk(E_0), \ldots, \rk(E_m))$.
\qed
\end{prop}

\begin{rem}
If $m=0$, then $E=E_0$ and $\Phi=0$, corresponding to the obvious
fixed points $(E,0)$, with $E$ a stable parabolic bundle.
\end{rem}

The following important fact was also noted by Simpson.  For a proof
see \cite[Proposition 3.11]{AG} (in fact, this deals with the
ordinary case but the argument can easily be adapted to the
parabolic case).

\begin{prop}\label{prop:vhs-stab}
  A parabolic complex variation of Hodge structure $(E = \bigoplus
  E_l, \Phi)$ is stable as a parabolic Higgs bundle if and only if
  the stability condition is satisfied for subbundles of $E$ which
  respect the decomposition $E = \bigoplus E_l$.
  \qed
\end{prop}

Next we need to calculate the weights of the circle action on the
tangent space to $\mathcal{M}$ at a critical point of $f$,
represented by $\mathbf{E} = (\bigoplus E_l,\Phi)$.  By the
characterization of the critical points provided by
Propositions~\ref{prop:critical=fixed} and \ref{prop:fixed=vhs},
we have decompositions
\begin{align*}
  \PE(E) &= \bigoplus_{l= -m}^{m} U_l, &
  \SPE(E) &= \bigoplus_{l= -m}^{m} \hat{U}_l,
\end{align*}
where we use the notation
\begin{align*}
  U_l &= \bigoplus_{j-i = l} \PH(E_i,E_j), &
  \hat{U}_l &= \bigoplus_{j-i = l} \SPH(E_i,E_j).
\end{align*}
We get a corresponding decomposition of the deformation complex
\begin{displaymath}
  C^\bullet(\mathbf{E}) = \bigoplus_{l=-m-1}^{m} C^\bullet(\mathbf{E})_l,
\end{displaymath}
where $C^\bullet(\mathbf{E})_l$ denotes the subcomplex
\begin{displaymath}
 C^\bullet(\mathbf{E})_l : \quad
  U_l \to \hat{U}_{l+1}\otimes K(D).
\end{displaymath}
With this notation we have the following result.

\begin{prop}\label{prop:weight-hh}
  Let $\mathbf{E} = (\bigoplus E_l,\Phi)$ represent a fixed point of
  the circle action on $\mathcal{M}$. Then the weight $l$ subspace of
  $T_{\mathbf{E}}\mathcal{M}$ is isomorphic to the first
  hypercohomology $\mathbb{H}^1(C^\bullet(\mathbf{E})_{-l})$.
\end{prop}

\begin{proof}
  It is clear that the derivative of the circle action at $\mathbf{E}
  = (E,\Phi)$
  is induced by the following map of deformation complexes
  $C^\bullet(E,\Phi) \to
  C^\bullet(E,e^{\ima\theta}\Phi)$:
  \begin{displaymath}
  \begin{CD}
    C^\bullet(E,\Phi):\quad
    @. \bigoplus U_l @>[-,\Phi]>> \bigoplus \hat{U}_{l+1} \otimes K(D) \\
   @VVe^{\ima\theta}V  @VV1V  @VVe^{\ima\theta}V \\
   C^\bullet(E,e^{\ima\theta}\Phi):\quad
   @. \bigoplus U_l @>[-,e^{\ima\theta}\Phi]>>
     \bigoplus \hat{U}_{l+1} \otimes K(D). \\
  \end{CD}
  \end{displaymath}
  In order to work out the circle action on
  $T_{\mathbf{E}}\mathcal{M}$ from this we need to determine
  the identification $\mathbb{H}^1(C^\bullet(E,\Phi)) \cong
  \mathbb{H}^1(C^\bullet(E,e^{\ima\theta}\Phi))$ induced by the
  isomorphism between $(E,\Phi)$ and $(E,e^{\ima\theta}\Phi)$.
  But it is easy to write down such an isomorphism $f_{\theta}$: with
  respect to the decomposition $E = \bigoplus E_l$ we can
  define $f_{\theta}$ to be multiplication by $e^{\ima l \theta}$ on
  $E_l$.  The corresponding isomorphism between the complexes
  $C^\bullet(E,\Phi)$ and $C^\bullet(E,e^{\ima\theta}\Phi)$ is given by
  the adjoint $\Ad(f_{\theta}) \colon \psi \mapsto f_{\theta}\psi
  f_{\theta}^{-1}$.  Note that $f_{\theta}$ is unique up to
  multiplication by scalars and hence $\Ad(f_{\theta})$ is unique.
  Since $\Ad(f_{\theta})$ is multiplication by $e^{\ima l \theta}$ on both
  $U_l$ and $\hat{U}_l$, we can write down the induced isomorphism of
  complexes; the piece in degree $l$ is given by
  \begin{displaymath}
  \begin{CD}
    C^\bullet(E,\Phi)_l:\quad
    @. U_l @>[-,\Phi]>>  \hat{U}_{l+1} \otimes K(D) \\
    @VV{\Ad(f_{\theta})}V @VVe^{\ima l \theta}V  @VVe^{\ima (l+1)\theta}V \\
      C^\bullet(E,e^{\ima\theta}\Phi)_l:\quad
    @. U_l @>[-,e^{\ima\theta}\Phi]>>  \hat{U}_{l+1} \otimes K(D). \\
  \end{CD}
  \end{displaymath}
  It follows that the derivative of the action of $e^{\ima\theta}$ is
  the endomorphism of $\HH^1(C^\bullet(E,\Phi))$ induced by the composite
  map of complexes
  \begin{displaymath}
    C^\bullet(E,\Phi)
    \overset{e^{\ima\theta}}{\longrightarrow}
    C^\bullet(E,e^{\ima\theta}\Phi)
    \overset{\Ad(f_{\theta})^{-1}}{\longrightarrow}
    C^\bullet(E,\Phi),
  \end{displaymath}
  whose degree $l$ piece is
  \begin{displaymath}
  \begin{CD}
    C^\bullet(E,\Phi)_l:\quad
    @. U_l @>[-,\Phi]>>  \hat{U}_{l+1} \otimes K(D) \\
    @VVV @VVe^{-\ima l \theta}V  @VVe^{-\ima l \theta}V \\
      C^\bullet(E,\Phi)_l:\quad
    @. U_l @>[-,\Phi]>>  \hat{U}_{l+1} \otimes K(D). \\
  \end{CD}
  \end{displaymath}
  Thus $\mathbb{H}^1(C^\bullet(E,\Phi)_l)$ is isomorphic to the weight
  $-l$ subspace of $\mathbb{H}^1(C^\bullet(E,\Phi)) \cong
  T_{\mathbf{E}}\mathcal{M}$.
\end{proof}

Summarizing the results of this section so far, we obtain the
following.

\begin{thm}\label{thm:morse-facts}
  The function $f\colon\mathcal{M} \to \RR$ defined by
  $f([A,\Phi]) = \frac12 \norm{\Phi}^2$ is a perfect Bott--Morse function.  A
  parabolic Higgs bundle $(E,\Phi)$ represents a critical point of
  $f$ if and only if it is a parabolic complex variation of Hodge
  structure, i.e., $E = \bigoplus_{l=0}^{m} E_l$ with $\Phi_l = \Phi_{\vert E_l}
  \colon
  E_l \to E_{l+1} \otimes K(D)$ strongly parabolic (where $\Phi=0$ if and
  only if $m=0$).  The tangent
  space to $\mathcal{M}$ at a critical point $\mathbf{E}$
  decomposes as
  \begin{displaymath}
    T_{\mathbf{E}}\mathcal{M}
    = \bigoplus_{l=-m}^{m+1} T_{\mathbf{E}}\mathcal{M}_l,
  \end{displaymath}
  where the eigenvalue $l$ subspace of the Hessian of $f$ is
  \begin{displaymath}
    T_{\mathbf{E}}\mathcal{M}_{l} \cong
    \mathbb{H}^1(C^\bullet(E,\Phi)_{-l}).
  \end{displaymath}
\end{thm}

\begin{proof}
  Immediate from Propositions \ref{prop:critical=fixed} and
  \ref{prop:weight-hh}.  Note that since our Morse function $f$ is
  \emph{minus} the moment map $\tilde{f}$ (cf.\ (\ref{eq:aaa3}) and
  (\ref{eq:aaa4})), the eigenvalue $l$ subspace of the Hessian
  coincides with the weight $l$ subspace for the circle action (with
  the \emph{same} sign).
\end{proof}

\subsection{Morse indices}
\label{sec:morse-indices}

\begin{prop}\label{prop:l-duality}
  \begin{enumerate}
  \item[(i)] There is a natural isomorphism
  \begin{displaymath}
    \mathbb{H}^1(C^\bullet(\mathbf{E})_{l}) \cong
    \mathbb{H}^1(C^\bullet(\mathbf{E})_{-l-1})^*
  \end{displaymath}
  and hence a natural isomorphism
  \begin{displaymath}
    T_{\mathbf{E}}\mathcal{M}_{l} \cong
    (T_{\mathbf{E}}\mathcal{M}_{1-l})^*.
  \end{displaymath}
  \item[(ii)] If $\mathbf{E}$ is stable, then we have
    \begin{align*}
      \mathbb{H}^0(C^\bullet(\mathbf{E})_{l}) &=
      \begin{cases}
        \CC &\text{if \ $l=0$,} \\
        0 &\text{otherwise,}
      \end{cases} \\
    \intertext{and}
      \mathbb{H}^2(C^\bullet(\mathbf{E})_{l}) &=
      \begin{cases}
        \CC &\text{if \ $l=-1$,} \\
        0 &\text{otherwise.}
      \end{cases} \\
    \end{align*}
  \end{enumerate}
\end{prop}

\begin{proof}
  (i) It follows from Proposition~\ref{prop:par-duality} (i) that
  there is an isomorphism of complexes
  \begin{displaymath}
    (C^\bullet(\mathbf{E})_{l})^* \otimes K \cong
    C^\bullet(\mathbf{E})_{-l-1}.
  \end{displaymath}
  Hence Serre duality for hypercohomology gives the first isomorphism
  of the statement.  The second isomorphism is now
  immediate from the last statement of Theorem~\ref{thm:morse-facts}.

  (ii) When $\mathbf{E}$ is
  stable we have that $\mathbb{H}^0(C^\bullet(\mathbf{E})) \cong \CC$,
  generated by the identity endomorphism of $\mathbf{E}$, and hence
  the first statement follows.
  For the same reason as in the proof of (i) we have the isomorphism
  $\mathbb{H}^0(C^\bullet(\mathbf{E})_{l}) \cong
  \mathbb{H}^2(C^\bullet(\mathbf{E})_{-l-1})^*$ and thus the second
  statement follows from the first.
\end{proof}

\begin{cor}\label{prop:down-morse-dim}
  Let $\mathbf{E}$ represent a critical point of $f$, let
  $T_{\mathbf{E}}\mathcal{M}_{\leq 0}$ be the subspace of the
  tangent space on which the Hessian of $f$ has eigenvalues less
  than or equal to zero and let
  $T_{\mathbf{E}}\mathcal{M}_{>0}$ be the subspace on which
  the Hessian of $f$ has eigenvalues greater than zero.  Then
  \begin{displaymath}
    T_{\mathbf{E}}\mathcal{M}_{\leq 0} \cong
    (T_{\mathbf{E}}\mathcal{M}_{>0})^*
  \end{displaymath}
  under the isomorphism of Proposition~\ref{prop:par-duality} {\rm (ii)}.  It
  follows that the dimension of
  $T_{\mathbf{E}}\mathcal{M}_{\leq 0}$ is half the dimension
  of the moduli space, i.e.,
  \begin{displaymath}
    \dim T_{\mathbf{E}}\mathcal{M}_{\leq 0}
    = r^2(g-1) + 1 + \sum_p f_p.
  \end{displaymath}
\end{cor}

\begin{proof}
  Immediate from Propositions~\ref{prop:l-duality}
  and~\ref{prop:phb-moduli-dim}.
\end{proof}

\begin{prop}\label{prop-index}
  Let the parabolic Higgs bundle $\mathbf{E}=(E,\Phi)$ represent
a critical point of $f$.
Then the Morse index of $f$ at this point is
  \begin{align*}
    \lambda_{\mathbf{E}}
    &= r^2(2g-2) + 2\sum_p f_p + 2\chi(C^\bullet(\mathbf{E})_0) \\
    &= r^2(2g-2) + 2\sum_p f_p + 2\sum_{l=0}^{m} \chi\bigl(\PE(E_l)\bigr)
    - 2\sum_{l=0}^{m-1}\chi\bigl(\SPH(E_l,E_{l+1}) \otimes K(D)\bigr)\\
    & = r^2(2g-2) + 2\sum_p f_p + 2\sum_{l=0}^{m}\Bigl[
    (1-g-n)\rk(E_l)^2 +\sum_p \dim P_p(E_l,E_l)\Bigl]\\
     &\qquad+2\sum_{l=0}^{m-1}\Bigl[(1-g)\rk(E_l)\rk(E_{l+1})
    - \rk(E_{l})\deg(E_{l+1}) + \rk(E_{l+1})\deg(E_{l}) \\
     &\qquad -\sum_p\dim N_p(E_l,E_{l+1})\Bigr],
  \end{align*}
where   $E = \bigoplus_{l=0}^{m} E_l$ with $\Phi\in
  H^0(\SPH(E_l,E_{l+1}) \otimes K(D))$.
\end{prop}

\begin{proof}
  Since we are calculating real dimensions, the Morse index is
  twice the dimension of $T_{\mathbf{E}}\mathcal{M}_{< 0}$,
  the subspace on which the Hessian of $f$ has negative eigenvalues.
  Hence Corollary~\ref{prop:down-morse-dim} shows that
  \begin{align*}
    \frac{1}{2} \lambda_{\mathbf{E}} &=
    \dim T_{\mathbf{E}}\mathcal{M}_{< 0} \\
    &= \dim T_{\mathbf{E}}\mathcal{M}_{\leq 0}
    - \dim T_{\mathbf{E}}\mathcal{M}_{0} \\
    &= r^2(g-1) + 1 + \sum_p f_p
    - \dim T_{\mathbf{E}}\mathcal{M}_{0}.
  \end{align*}
  On the other hand from Proposition~\ref{prop:l-duality} (ii) we
  have that $\mathbb{H}^0(C^\bullet(\mathbf{E})_0) = \CC$, while
  $\mathbb{H}^2(C^\bullet(\mathbf{E})_0) = 0$.  Hence
  Theorem~\ref{thm:morse-facts} shows that we have
  \begin{align*}
    \dim T_{\mathbf{E}}\mathcal{M}_{0}
    &= \dim \mathbb{H}^1(C^\bullet(\mathbf{E})_0)  \\
    &= 1 - \chi(C^\bullet(\mathbf{E})_0),
  \end{align*}
  and this finishes the proof of the first identity of the statement
  of the Proposition.  The rest  can be deduced from the
  long exact sequence in hypercohomology for the complex
  $C^\bullet(\mathbf{E})_0$, analogous to \eqref{eq:les-phiggs},
and using the same method as in the proof of
Proposition~\ref{prop:phb-moduli-dim}.
\end{proof}

\begin{rem}
Obviously, the absolute minima  is for $m=0$, for which the
computation in Proposition~\ref{prop-index} naturally gives
  \begin{align*}
    \lambda_{(E,0)}
    &= r^2(2g-2) + 2\sum_p f_p +  2\chi(\PE(E)) \\
    &= r^2(2g-2) + 2\sum_p f_p +2 r^2 (1-g)+ 2 \sum_p(r^2 -f_p -r^2) \\
    & =0.
  \end{align*}
\end{rem}

\subsection{Rank three  parabolic Higgs bundles}
\label{sec:betti-higgs}

Now we turn our attention to the moduli space $\cM$ of
parabolic Higgs bundles of rank three. Let $(E,\Phi)$  be a
critical point of $f$. By
Theorem \ref{thm:morse-facts}, the only
possibilities that we have in this situation are:

\begin{enumerate}
 \item[(a)] $E$ is a stable rank three  parabolic Higgs bundle and
$\Phi=0$.
 \item[(b)] $E=E_0 \oplus E_1\oplus E_2$ where $E_l$ are
parabolic line bundles. These are line bundles with weights at
each  $p\in D$. The map $\Phi$ decomposes as strongly parabolic
maps $\Phi_0:E_0\to E_1\otimes K(D)$ and $\Phi_1:E_1\to E_2\otimes
K(D)$.
 \item[(c)] $E=E_0\oplus E_1$ where $E_0$ is a parabolic
line bundle $L$ and $E_1$ is a rank $2$ parabolic bundle. Here
$\Phi$ gives a strongly parabolic map $\Phi_0:L\to E_1\otimes
K(D)$.
 \item[(d)] $E=E_0\oplus E_1$ where $E_0$ is a rank $2$ parabolic
bundle and $E_1$ is a parabolic line bundle $L$. Here $\Phi$ gives
a strongly parabolic map $\Phi_0:E_0\to L\otimes K(D)$.
\end{enumerate}

In case (a) the corresponding critical subvariety can obviously be
identified with the moduli space of ordinary parabolic bundles.
Its Betti numbers can be computed from a formula given  by Nitsure
\cite{Ni} and Holla \cite{Ho}.  In Section~\ref{par-rank3} below
we work out explicitly what their formula gives for the Poincar\'e
polynomial in our situation of rank three parabolic bundles. Case
(b) involves basically line bundles and divisors and can be dealt
with easily \cite{McD}. The other two cases, (c) and (d), are more
involved. They are particular cases of objects called {\em
parabolic triples}, which have been introduced and studied from a
gauge-theoretic point of view in \cite{BiG}, and will be studied
in Section~\ref{sec:stable-triples} below.

\subsection{Laumon's Theorem for parabolic Higgs bundles}
\label{sec:laumon}

At this point we make a small digression in order to deduce,
following Hausel, a parabolic version of a Theorem of Laumon from
the analysis leading to our calculation of the Morse indices.

As in the non-parabolic case studied by Hitchin \cite{H, H2},
there is a \emph{Hitchin map}
 $$
  \chi\colon \mathcal{M} \to B=\bigoplus_{i=1}^{r} H^0(K(D)^i)
 $$
defined by taking the parabolic Higgs bundle $(E,\Phi)$ to the
characteristic polynomial of $\Phi$.  Since the Higgs field is
strictly parabolic, this map takes values in a subspace of $B$ of
dimension $r^2(g-1) + 1 + \sum_p f_p$, i.e., half the dimension of
$\mathcal{M}$. The Hitchin map defines an algebraic completely
integrable system. This means that the $r^2(g-1) + 1 + \sum_p f_p$
functions defined by $\chi$ Poisson commute, their differentials are
linearly independent and the generic fibre of $\chi$ is an open set in
an abelian variety.  In fact $\mathcal{M}$ is a symplectic leaf of a
Poisson manifold equipped with the structure of a generalized
integrable system, see Bottacin \cite{Bo} and Markman
\cite{markman:1994}.

The pre-image of $0$ under the Hitchin map,
\begin{displaymath}
  N = \chi^{-1}(0),
\end{displaymath}
is called the \emph{nilpotent cone}.  The main result of Laumon
\cite{La}, proved for the moduli stack of Higgs bundles, is that
the nilpotent cone is Lagrangian.

In the non-parabolic case, Hausel \cite[Theorem~5.2]{Ha} proved that
the downwards Morse flow on the moduli space of Higgs bundles
coincides with the nilpotent cone.  His proof goes over word by word
to the parabolic case, so we have the following Theorem.
\begin{thm}
  The downwards Morse flow on the moduli space of parabolic Higgs
  bundles coincides with the
  nilpotent cone $N$.
  \qed
\end{thm}
As pointed out by Hausel, the nilpotent cone is isotropic because the
Hitchin map is a completely integrable system. (To be precise, the
nilpotent cone being isotropic means that its tangent space at any
non-singular point of the nilpotent cone is an isotropic subspace of
the tangent space to $\mathcal{M}$, cf.\ Ginzburg \cite{Gi}.)  Hence
the nilpotent cone is Lagrangian if its dimension equals half that of
the moduli space $\mathcal{M}$.  But this fact follows at once from
our Corollary~\ref{prop:down-morse-dim}.  Thus we have the following
version of Laumon's theorem for the moduli space of parabolic Higgs
bundles.

\begin{thm}\label{thm:laumon}
  The nilpotent cone $N$ is a Lagrangian subvariety of the moduli space of
  parabolic Higgs bundles.
  \qed
\end{thm}

\section{Parabolic triples}
\label{sec:stable-triples}

\subsection{Definitions and basic facts}
\label{sec:triples-definitions}

A \emph{parabolic triple}  $T = (E_{1},E_{2},\phi)$ on $X$ consists
of two parabolic vector bundles $E_{1}$ and $E_{2}$ on $X$ and a
$\phi\in H^0(\SPH( E_{2}, E_{1}(D)))$. A homomorphism from $T' =
(E_1',E_2',\phi')$ to $T = (E_1,E_2,\phi)$   is a commutative
diagram
\begin{displaymath}
  \begin{CD}
    E_2' @>\phi'>> E_1'(D) \\
    @VVV @VVV  \\
    E_2 @>\phi>> E_1(D),
  \end{CD}
\end{displaymath}
where the vertical arrows are parabolic sheaf homomorphisms. A
triple $T'=(E_1',E_2',\phi')$ is a subtriple of $T = (E_1,E_2,\phi)$
if the sheaf homomorphims  $E_1'\to E_1$ and $E_2'\to E_2$ are
injective. A subtriple $T'\subset T$  is called {\em proper} if
$T'\neq 0 $ and $T'\neq T$.

\begin{defn}
For any $\sigma \in \RR$ the \emph{$\sigma$-degree} and
\emph{$\sigma$-slope} of $T$ are defined to be
\begin{align*}
  \deg_{\sigma}(T)
  &= \pdeg(E_{1}) + \pdeg(E_{2}) + \sigma \
  \rk(E_{2}), \\
  \mu_{\sigma}(T)
  &=
  \frac{\deg_{\sigma}(T)}
  {\rk(E_{1})+\rk(E_{2})} \\
  &= \pmu(E_{1} \oplus E_{2}) +
  \sigma\frac{\rk(E_{2})}{\rk(E_{1})+
    \rk(E_{2})}.
\end{align*}

We say  $T = (E_{1},E_{2},\phi)$ is \emph{$\sigma$-stable} if
 $$
 \mu_{\sigma}(T')  < \mu_{\sigma}(T),
 $$
for any proper subtriple $T' = (E_{1}',E_{2}',\phi')$. We define
\emph{$\sigma$-semistability} by replacing the above strict
inequality with a weak inequality. A triple is called
\emph{$\sigma$-polystable} if it is the direct sum of
$\sigma$-stable triples of the same $\sigma$-slope.
\end{defn}

Let us fix the topological and parabolic types  of
$E_1$ and  $E_2$. We denote by   $\mathcal{N}_\sigma$
the moduli space of $\sigma$-stable triples
$T = (E_{1},E_{2},\phi)$ of  the given type.

Given a triple  $T=(E_1,E_2,\phi)$ one has the dual triple
$T^*=(E_2^*,E_1^*,\phi^*)$, where $E_i^*$ is the parabolic dual of $E_i$ and
$\phi^*$ is the transpose of $\phi$. The following  is not difficult
to prove.
\begin{prop}
\label{prop:duality}
The  $\sigma$-(semi)stability of $T$   is equivalent to the
$\sigma$-(semi)stability of $T^*$.
The map $T\mapsto T^*$ defines an isomorphism of moduli
spaces. \qed
\end{prop}

This can be used to restrict  our study to $\rk(E_1)\geq \rk(E_2)$
and appeal to duality to deal with the case $\rk(E_1)<\rk(E_2)$.

There are certain necessary conditions in order for
$\sigma$-semistable triples to exist. Let $r_1=\rk (E_1)$,
$r_2=\rk (E_2)$, $\pmu_1=\pmu(E_1)$ and $\pmu_2=\pmu(E_2)$ be the
ranks and parabolic degrees of $E_1$ and $E_2$, and define
\begin{align}
   \sigma_m= &\pmu_1-\pmu_2, \label{alpha-bounds-m} \\
      \sigma_M = & \left(1+ \frac{r_1+r_2}{|r_1 - r_2|}\right)(\pmu_1 - \pmu_2)
                   + n \ \frac{r_1+r_2}{|r_1 - r_2|}\ ,
     \;\;\; \text{if }\ \ r_1\neq r_2. \label{alpha-bounds-bigM}
\end{align}

\begin{prop} \label{prop:alpha-range} A necessary condition for
$\mathcal{N}_\sigma$ to be non-empty is
\begin{enumerate}
 \item[(i)] $\sigma_m \leq \sigma \leq \sigma_M$, \  if  \  $r_1\neq
r_2$,
 \item[(ii)] $\sigma_m \leq \sigma$, \    if   \ $r_1= r_2$.
\end{enumerate}
\end{prop}

\begin{proof}
The proof is similar to the one given  in \cite[Proposition
3.18]{BG} for ordinary triples.
\end{proof}

\begin{rem}
The upper bound given for $\sigma$ is not optimal. A better one
can be found, as will be seen later in Section
\ref{sec:thaddeus-program}.
\end{rem}

Using the  dimensional reduction construction given in \cite{BiG}, the
moduli space $\mathcal{N}_\sigma$ can be realised as a  subvariety
of a certain moduli space of parabolic bundles on $X\times \PP^1$.
Such moduli spaces have been constructed by
Maruyama and Yokogawa \cite{MY} in arbitrary dimensions using  GIT methods.

Another important aspect that follows also from the dimensional
reduction point of view is the existence of  a correspondence
between stability and the existence of solutions to certain
gauge-theoretic equations on a parabolic triple
$T=(E_1,E_2,\phi)$, known as the {\em parabolic vortex equations}
\cite{BiG}. The parabolic vortex equations
  \begin{equation}
    \label{eq:coupled-vortex}
  \begin{split}
    \ima\Lambda F(E_1) + \phi\phi^{*}
      &= \tau_{1} \, \Id_{E_1},  \\
    \ima\Lambda F(E_2)  - \phi^{*}\phi
      &= \tau_{2} \, \Id_{E_2},
  \end{split}
  \end{equation}
are equations for Hermitian metrics  on $E_1$ and $E_2$ adapted to
the parabolic structure. Here $\Lambda$ is contraction by the
K\"ahler form of a metric on $X$ (normalized so that
$\vol(X)=2\pi$), $F(E_i)$ is the curvature of the unique connection
on $E_i$ compatible with the Hermitian metric and the holomorphic
structure of $E_i$, and $\tau_1$ and $\tau_2$ are real parameters
satisfying $\pdeg(E_1)+\pdeg(E_2)=r_1\tau_1+r_2\tau_2$. Also, here
$\phi^\ast$ is the adjoint of $\phi$ with respect to the Hermitian
metrics. One has the following.

\begin{thm}\cite[Theorem 3.4]{BiG}
\label{thm:vortices-hitchin-kobayashi}
A  solution to \eqref{eq:coupled-vortex} exists if and only if $T$ is
$\sigma$-polystable for $\sigma=\tau_1-\tau_2$.
\qed
\end{thm}

\subsection{Parabolic Higgs bundles and parabolic triples}
\label{sec:higgs-bundles-triples}

The relation between parabolic Higgs bundles and parabolic triples
is given by the following.

\begin{prop}\label{prop:august}
Suppose that $(E,\Phi)$ is a stable parabolic Higgs bundle such
that $E=E_1\oplus E_2$ and
$$
 \Phi=
  \begin{pmatrix}
    0& \phi \\
    0 & 0
  \end{pmatrix}
$$
with $\phi:E_2\to E_1\otimes K(D)$ a strongly parabolic map. Then
$(E,\Phi)$ is stable if and only if the parabolic triple
$(E_1\otimes K, E_2,\phi)$ is $\sigma$-stable for $\sigma=2g-2$.
\end{prop}

\begin{proof}
Take a sub-object $E'\subset E$ with $\Phi(E')\subset E'\ox K(D)$.
This can be assumed to be of the form $E'= E'_1\oplus E'_2$ and
hence it defines a subtriple $(E_1'\ox K,E_2',\phi')$ where
$\phi'=\phi|_{E_2'}$. The result follows now from the equivalence
between
 $$
\frac{\pdeg(E_1')+\pdeg(E_2')}{r_1'+r_2'} <
\frac{\pdeg(E_1)+\pdeg(E_2)}{r_1+r_2}, \;\; \mbox{and}
 $$
 $$
 \frac{\pdeg(E_1') + r_1'(2g-2)+\pdeg(E_2')}{r_1'+r_2'}+
\s \frac{r_2'}{r_1'+r_2'} < \qquad \qquad
 $$
 $$
 \qquad \qquad \frac{\pdeg(E_1)+r_1(2g-2)+ \pdeg(E_2)}{r_1+r_2}
 +\s\frac{r_2}{r_1+r_2},
 $$
which is the $\s$-stability of the  triple $(E_1\ox K, E_2, \Phi)$,
for $\s=2g-2$.
\end{proof}

\subsection{Extensions and deformations of parabolic triples}
\label{sec:extensions-of-triples}

In order to analyse the differences between the  moduli spaces
$\mathcal{N}_\sigma$ as  $\sigma$ changes, as well as the
smoothness properties of the moduli space for a given value of
$\sigma$, we need to study the homological algebra of parabolic
triples. This is done by considering the hypercohomology of  a
certain complex of sheaves, in an analogous way to the case of
holomorphic triples studied in \cite{BGG}, and the parabolic Higgs
bundle case studied in Subsection \ref{sec:deformation}.

Let $T'=(E'_1,E'_2,\phi')$ and $T''=(E''_1,E''_2,\phi'')$ be two parabolic
triples. Let $\Hom(T'',T')$ denote the linear
space of homomorphisms from  $T''$ to
$T'$, and let $\Ext^1(T'',T')$  denote the linear space of equivalence
classes of extensions of the form
\begin{displaymath}
  0 \too T' \too T \too T'' \too 0,
\end{displaymath}
where by this we mean a commutative  diagram
  $$
  \begin{CD}
  0@>>>E_2'@>>>E_2@>>> E_2''@>>>0\\
  @.@V\phi' VV@V \phi VV@V \phi'' VV\\
  0@>>>E'_1(D)@>>>E_1(D)@>>>E_1''(D)@>>>0.
  \end{CD}
  $$
Hence, to analyse $\Ext^1(T'',T')$ one considers the
complex of sheaves
\begin{equation}
  \label{eq:extension-complex}
    C^{\bullet}(T'',T'):
  \PH(E_{1}'', E_{1}') \oplus  \PH(E_{2}'', E_{2}')
  \overset{c}{\too}
  \SPH(E_{2}'', E_{1}'(D)),
\end{equation}
where the map $c$ is defined by
\begin{displaymath}
  c(\psi_{1},\psi_{2}) = \phi'\psi_{2} - \psi_{1}\phi''.
\end{displaymath}

\begin{prop}
  \label{prop:hyper-equals-hom}
  There are natural isomorphisms
  \begin{align*}
    \Hom(T'',T') &\cong \HH^{0}(C^{\bullet}(T'',T')), \\
    \mbox{\em Ext}^{1}(T'',T') &\cong \HH^{1}(C^{\bullet}(T'',T')),
  \end{align*}
and a long exact sequence associated to the complex
$C^{\bullet}(T'',T')$:
\begin{equation}
  \label{eq:long-exact-extension-complex}
\begin{array}{ll}
  0 &\ar \mathbb{H}^0
  \ar  H^0(\PH(E_{1}'', E_{1}') \oplus \PH(E_{2}'', E_{2}'))
   \ar   H^0(\SPH(E_{2}'',  E_{1}'(D))) \\
    & \ar \mathbb{H}^1
  \ar   H^1(\PH(E_{1}'', E_{1}') \oplus \PH(E_{2}'', E_{2}'))
 \ar  H^1(\SPH(E_{2}'',  E_{1}'(D))) \\
   & \ar \mathbb{H}^2 \ar  0,
\end{array}
\end{equation}
where $\mathbb{H}^i=\mathbb{H}^i(C^{\bullet}(T'',T'))$.
\qed
\end{prop}

We introduce the following notation:
\begin{align}
  h^{i}(T'',T') &= \dim\HH^{i}(C^{\bullet}(T'',T')), \notag \\
  \chi(T'',T') &= h^0(T'',T') - h^1(T'',T') + h^2(T'',T'). \notag
\end{align}

\begin{prop}
  \label{prop:chi(T'',T')}
  For any parabolic  triples $T'$ and $T''$ we have
  \begin{align}\notag
    \chi(T'',T') &= \chi(\PH(E_{1}'', E_{1}'))
    + \chi(\PH(E_{2}'',  E_{2}'))
    - \chi(\SPH(E_{2}'', E_{1}'(D)))
  \end{align}
where $\chi(E)=\dim H^0(E) - \dim H^1(E)$ is the Euler characteristic
of $E$.
\end{prop}

\begin{proof}
  Immediate from the long exact sequence
  \eqref{eq:long-exact-extension-complex} and the Riemann--Roch formula.
\end{proof}
\begin{cor}
  \label{cor:chi-relation}
  For any extension $0 \to T' \to T \to T'' \to 0$ of parabolic triples,
  \begin{displaymath}
    \chi(T,T) = \chi(T',T') + \chi(T'',T'') + \chi(T'',T') +
    \chi(T',T'').
  \end{displaymath}
  \qed
\end{cor}

\begin{prop}
   \label{prop:h0-vanishing}
   Suppose that $T'$ and $T''$ are $\sigma$-semistable.
   \begin{enumerate}
 \item[(i)]
   If $\mu_\sigma(T')<\mu_\sigma (T'')$ then
 $\HH^{0}(C^{\bullet}(T'',T')) \cong 0$.
 \item[(ii)]
 If $\mu_\sigma(T')=\mu_\sigma (T'')$ and $T'$ and $T''$ are both
   $\sigma$-stable, then
  $$
     \HH^{0}(C^{\bullet}(T'',T')) \cong
     \begin{cases}
       \CC\, , \quad &\text{if \ $T' \cong T''$}, \\
       0 \, ,\quad &\text{if \ $T' \not\cong T''$}.
     \end{cases}
 $$
 \end{enumerate}
\qed
   \end{prop}

\begin{cor}\label{dimension-ext1}
Let  $T'$ and $T''$ be  $\sigma$-semistable parabolic triples
with  $\mu_\sigma(T')=\mu_\sigma (T'')$, and suppose that
 $\HH^{2}(C^{\bullet}(T'',T')) =0$. Then
 $$
 \dim \mbox{\em Ext}^1(T'',T')=    h^0(T'',T')  - \chi(T'',T').
 $$
\qed
\end{cor}

Since the  space of infinitesimal deformations of $T$ is
isomorphic to $\HH^{1}(C^{\bullet}(T,T))$, the considerations of
the previous sections also apply to studying deformations of a
parabolic triple $T$ (the proofs are analogous to the
non-parabolic case \cite{BGG}). To be precise, one has the
following.

\begin{thm}\label{thm:smoothdim}
Let $T=(E_1,E_2,\phi)$ be a $\sigma$-stable parabolic triple.
\begin{itemize}
\item[(i)] The Zariski tangent space at the point defined by $T$
in the moduli space of stable triples  is isomorphic
to $\HH^{1}(C^{\bullet}(T,T))$.
\item[(ii)]
If $\HH^{2}(C^{\bullet}(T,T))= 0$, then the moduli space of $\sigma$-stable
parabolic triples is smooth in  a neighbourhood of the point defined by $T$.
\item[(iii)]
$\HH^{2}(C^{\bullet}(T,T))= 0$ if and only if the homomorphism
 $$
  H^1(\PE(E_1)) \oplus H^1(\PE(E_2))
  \too  H^1(\SPH(E_2, E_1(D)))
 $$
in the corresponding long exact sequence is surjective.
\item[(iv)] At a smooth point $T\in \mathcal{N}_\sigma$
the dimension of the moduli space of $\sigma$-stable parabolic triples is
\begin{equation}\label{dimension-triples}
  \begin{split}
  \dim \mathcal{N}_\sigma
  &= h^{1}(T,T) = 1 - \chi(T,T) \\
  &= \chi(\PE(E_{1}, E_{1}))
    + \chi(\PE(E_{2},  E_{2}))
    - \chi(\SPH(E_{2}, E_{1}(D)))
\end{split}
\end{equation}
\item[(v)] If $\phi$ is injective or surjective then $T=(E_1,E_2,\phi)$
defines a smooth point in the moduli space.
\end{itemize}
\qed
\end{thm}

\section{Critical values and flips}
\label{sec:critical-values-flips}

\subsection{Critical values}
\label{sec:critical-values}

A parabolic triple $T=(E_1,E_2,\phi)$ of fixed topological  and
parabolic type is strictly $\sigma$-semistable if and only if it has
a proper subtriple $T'=(E_1',E_2',\phi')$ such that
$\mu_{\sigma}(T')= \mu_{\sigma}(T)$, i.e.
\begin{equation}
  \label{eq:strict-alpha-ss}
  \pmu(E'_1 \oplus E'_2) + \sigma \frac{r'_2}{r_1'+r_2'}
  =  \pmu(E_1 \oplus E_2) + \sigma \frac{r_2}{r_1+r_2},
\end{equation}
where $r_1, r_2,r_1', r_2'$ are the ranks of $E_1, E_2,E_1',
E_2'$. There are two ways in which this can happen. The first one
is if there exists a subtriple $T'$ such that
\begin{align*}
  \frac{r'_2}{r_1'+r_2'} &= \frac{r_2}{r_1+r_2},\;\; \mbox{and} \\
  \pmu(E'_1 \oplus E'_2) &= \pmu(E_1 \oplus E_2).
\end{align*}
In this case the terms containing $\sigma$ drop from
\eqref{eq:strict-alpha-ss} and $T$ is strictly $\sigma$-semistable
for all values of $\sigma$.  We refer to this phenomenon as
\emph{$\sigma$-independent semistability}.

The other way in which strict $\sigma$-semistability can happen is
if equality holds in \eqref{eq:strict-alpha-ss} but
 $$
  \frac{r'_2}{r_1'+r_2'} \neq \frac{r_2}{r_1+r_2}.
 $$
The values of $\sigma$ for which this happens are called critical
values.

{}From now on we shall make the following assumption on the
weights.
\begin{assumption}\label{assumption}
  Let $\alpha_i(p)$ be the collection of all the weights of $E_1$ and
  $E_2$ together. We assume that they are all of multiplicity one
  and that, for a large integer $N$ depending only on the ranks,
  they satisfy the following property:
 $$
 \sum_{1\leq i\leq r, \, p\in D} n_{i,p}\, \alpha_i(p)
 \in \ZZ, \;\; n_{i,p} \in \ZZ, |n_{i,p}| \leq N
 \Longrightarrow n_{i,p}=0, \, \text{for all } p, i \, .
 $$
  The weights failing this genericity condition are a finite union of hyperplanes
  in $[0,1)^{nr}$, where $nr$ is the total number of weights, $r=r_1+r_2$.
\end{assumption}

\begin{prop}
\label{triples-critical-range}
\begin{itemize}
 \item[(i)] Under Assumption \ref{assumption}, there are no
$\sigma$-independent semistable triples (by taking $N$ larger than
$r_1+r_2$).
 \item[(ii)] The critical values of $\sigma$ form a discrete subset
of $[\sigma_m,\infty)$, where $\sigma_m$ is as in
(\ref{alpha-bounds-m}).
 \item[(iii)] If $r_1\neq r_2$ the number
of critical values is  finite and they  lie in the interval
$[\sigma_m,\sigma_M]$, where $\sigma_M$ is as in
(\ref{alpha-bounds-bigM}).
 \item[(iv)] The stability criteria for
two values of  $\sigma$ lying between two consecutive critical
values are equivalent; thus the corresponding moduli spaces
coincide.
 \item[(v)] If $\sigma_c$ is a critical value and $T'$ is a subtriple of a
$\sigma_c$-semistable triple $T$ such that
$\mu_{\sigma_c}(T')=\mu_{\sigma_c}(T)$, then $T'$ and the quotient
triple $T''=T/T'$ are $\sigma_c$-stable (for this, it may be
necessary to take a larger value of $N$ in Assumption
\ref{assumption}).
\end{itemize}
\qed
\end{prop}

\subsection{Crossing critical values and universal extensions}
\label{sec:crossing-critical-values}

In this section we study the differences between the  moduli spaces
$\mathcal{N}_\sigma$, for fixed type but different values of $\sigma$.

We begin with a set theoretic description of the differences
between two spaces $\mathcal{N}_{\sigma}$  and
$\mathcal{N}_{\sigma'}$ when $\sigma$ and $\sigma'$ are separated
by a single critical value (as defined in Subsection
\ref{sec:triples-definitions}). For the rest of this section we
adopt the following notation: when $r_1\ne r_2$ the bounds
$\sigma_m$ and $\sigma_M$ are as in (\ref{alpha-bounds-m}) and
(\ref{alpha-bounds-bigM}). When $r_1=r_2$ we adopt the convention
that $\sigma_M=\infty$. Let $\sigma_c\in \RR $ be a critical value
such that
 $$
 \sigma_m \leq \sigma_c \leq \sigma_M.
 $$
Set
 $$
\sigma_c^+ = \sigma_c + \epsilon,\quad \sigma_c^- = \sigma_c -
\epsilon,
 $$
where $\epsilon > 0$ is small enough so that $\sigma_c$ is the
only critical value in the interval $(\sigma_c^-,\sigma_c^+)$.

\begin{defn}Let $\sigma_c$ be a critical value. We define the {\it flip loci}
$\mathcal{S}_{\sigma_c^{\pm}}\subset\mathcal{N}_{\sigma_c^{\pm}}$ by
the conditions that the points in $\mathcal{S}_{\sigma_c^+}$
represent triples which are $\sigma_c^+$-stable  but
$\sigma_c^-$-unstable, while the points in
$\mathcal{S}_{\sigma_c^-}$ represent triples which are
$\sigma_c^-$-stable  but $\sigma_c^+$-unstable.
\end{defn}


\begin{lem}\label{lemma:fliploci}
In the above notation,
 $$
 \mathcal{N}_{\sigma_c^+}-\mathcal{S}_{\sigma_c^+}=
\mathcal{N}_{\sigma_c}=
\mathcal{N}_{\sigma_c^-}-\mathcal{S}_{\sigma_c^-}.
 $$ \qed
\end{lem}

As a consequence of  Proposition \ref{triples-critical-range} (v)
we have the following.

\begin{prop}\label{prop:vicente}
Let $\sigma_c$  be a critical value. Let $T=(E_1,E_2,\phi)$ be a
triple of this type which is $\sigma_c$-semistable. Then $T$ has a
(unique) description as the middle term in an extension
 \begin{equation}\label{destab}
 0\to T'\to T\to T'' \to 0
 \end{equation}
 in which $T'$  and $T''$ are  $\sigma_c$-stable and
 $\mu_{\sigma_c}(T')=\mu_{\sigma_c}(T)=\mu_{\sigma_c}(T'')$.
\qed
\end{prop}

We thus have the following.

\begin{prop}\label{prop:critical-loci}
The set $\mathcal{S}_{\sigma_c^+}$ coincides with the set of equivalence
classes of extensions (\ref{destab}), in which $T'$  and $T''$ are
$\sigma_c$-stable,
$\mu_{\sigma_c}(T')=\mu_{\sigma_c}(T)=\mu_{\sigma_c}(T'')$,
and $r_2'/r'< r_2''/r''$.

Similarly, $\mathcal{S}_{\sigma_c^-}$ coincides with the set of equivalence
classes of extensions (\ref{destab}), in which $T'$  and $T''$ are
$\sigma_c$-stable,
$\mu_{\sigma_c}(T')=\mu_{\sigma_c}(T)=\mu_{\sigma_c}(T'')$,
and $r_2'/r'> r_2''/r''$; or, equivalently, extensions
$$
 0\to T''\to T\to T' \to 0
$$
where $T'$ and $T''$ are as above, but $r_2'/r'< r_2''/r''$.
\qed
\end{prop}

To construct the locus $\mathcal{S}_{\sigma_c^\pm}$, we first
observe that, by the genericity of the weights, the moduli spaces
$\mathcal{N}_{\sigma_c}'$ and $\mathcal{N}_{\sigma_c}''$ are fine
moduli spaces (cf.\ \cite{Y1}), i.e.,  there are universal
parabolic triples $\cT'=(\cE_1',\cE_2',\Phi')$ and
$\cT''=(\cE_1'', \cE_2'',\Phi'')$  over
$\mathcal{N}_{\sigma_c}'\times X$ and
$\mathcal{N}_{\sigma_c}''\times X$ respectively. Let
$B=\cN_{\sigma_c} '\x \cN_{\sigma_c}''$ and let  pull back $\cT'$
and $\cT''$ to $B\x X$. Considering  the complex
$C^\bullet(\cT'',\cT')$ as defined in
\eqref{eq:extension-complex}, taking relative hypercohomology
$\HH_\pi(C^\bullet(\cT'',\cT'))$ with respect to the projection
$\pi:B \x X\to B$, and putting
 \begin{equation}\label{def-W+}
 W^+:=\HH^1_\pi(C^\bullet(\cT'',\cT')),
 \end{equation}
we have the following exact sequence of sheaves over $B$:
 \begin{equation}\label{eq:relative-hyper}
      \begin{split}
  0 & \to  \HH^0_\pi(C^\bullet(\cT'',\cT'))\to\pi_*\PH(\cE_2'', \cE_2') \oplus \pi_*\PH(\cE_1'', \cE_1') \to
  \pi_*\SPH(\cE_2'', \cE_1'(D)) \\
   & \to  W^+ \to R^1\pi_*\PH(\cE_2'', \cE_2') \oplus R^1\pi_*\PH(\cE_1'', \cE_1') \to
  R^1\pi_*\SPH(\cE_2'', \cE_1'(D))\\
  &\to
 \HH^2_\pi(C^\bullet(\cT'',\cT'))
 \to 0.
 \end{split}
  \end{equation}

Analogously, we can   consider  the complex
$C^\bullet(\cT',\cT'')$ and define
$$
W^-:=\HH^1_\pi(C^\bullet(\cT',\cT'')).
$$

\begin{prop}\label{vanishing-h2}
If $\HH^2(C^\bullet(T'',T'))=0$ for every $(T',T'')\in
\mathcal{N}_{\sigma_c}'\times  \mathcal{N}_{\sigma_c}''$, then
$W^+$ defined in (\ref{def-W+}) is locally free. Similarly for
$W^-$.
\end{prop}

\begin{proof}
By Proposition \ref{prop:h0-vanishing},
$\HH^0_\pi(C^\bullet(\cT'',\cT'))=0$ for every $(T',T'')\in
\mathcal{N}_{\sigma_c}'\times  \mathcal{N}_{\sigma_c}''$
and hence  $\HH^0_\pi(C^\bullet(\cT'',\cT'))=0$.
By assumption $\HH^2_\pi(C^\bullet(\cT'',\cT'))=0$ and the result
thus follows from  (\ref{eq:relative-hyper}).
\end{proof}

\begin{rem}
In our applications, the vanishing assumption in Proposition
\ref{vanishing-h2}  will always be satisfied due to the small rank
of the bundles involved. In fact, the vanishing is probably true
in general for $\sigma\geq 2g-2$, as in the non-parabolic case
\cite{BGG}.
\end{rem}

Clearly, from  Proposition \ref{prop:critical-loci}, we have the
following.
\begin{prop}\label{prop:5.9}
If $\HH^2(C^\bullet(T'',T'))=0$ and  $\HH^2(C^\bullet(T',T''))=0$
for every $(T',T'')\in \mathcal{N}_{\sigma_c}'\times
\mathcal{N}_{\sigma_c}''$, then
 $$
 \cS_{\s_c^\pm}= \PP W^\pm.
 $$ \qed
\end{prop}

The following will be important to study the relation between
$\cN_{\s_c^-}$ and $\cN_{\s_c^+}$.

\begin{prop} \label{prop:5.6}
Over  $\PP W^+\x X$ there is a universal extension
 \begin{equation}\label{universal-ext+}
 0\to \cT'\otimes \cO_{\PP W^+}(1)\to \cT^+\to \cT'' \to 0,
 \end{equation}
where $\cT'\otimes \cO_{\PP W^+}(1):=(\cE_1'\otimes \cO_{\PP
W^+}(1),\cE_2'\otimes \cO_{\PP W^+}(1),\Phi')$ (we omit pull-backs
for clarity). Similarly, on  $\PP W^-\x X$ there is a universal
extension
 $$
  0\to \cT''\otimes \cO_{\PP W^-}(1)\to \cT^-\to \cT' \to 0,
 $$
where $\cT''\otimes \cO_{\PP W^-}(1):=(\cE_1''\otimes \cO_{\PP
W^-}(1),\cE_2''\otimes \cO_{\PP W^-}(1), \Phi'')$.
\end{prop}

\begin{proof}
The proof is analogous to the one given by Lange \cite{L} for
extensions of sheaves. In fact, our result could be derived from
that one by  making use of the correspondence between parabolic
triples over $X$ and $\Sl(2,\CC)$-invariant parabolic vector bundles
over $X\x \PP^1$ (cf. \cite{BiG}). Hence we only give the main
ingredients of the proof.

Let $(T',T'')\in B=\mathcal{N}_{\sigma_c}'\times
\mathcal{N}_{\sigma_c}''$. Let $\WW=\HH^1(X,C^\bullet(T'',T'))$,
and let $\PP=\PP(\WW)$. Over $\PP\x X$ there is a universal
extension
  \begin{equation}\label{universal-ext}
  0\to T'(1)\to T\to T'' \to 0,
  \end{equation}
where $T'(1):=(E_1'\otimes \cO_{\PP}(1), E_2'\otimes
\cO_{\PP}(1),\phi')$ and we are omitting pull-backs. By the
universal property of this extension we  mean that $T$ restricted
to $\{p\}\x X$ is a triple whose corresponding equivalence class
is precisely $p\in \PP$. Extensions like (\ref{universal-ext}) are
parametrised by $\HH^1(\PP\x X,C^\bullet(T'',T'(1)))$ which by the
K\"unneth formula is isomorphic to
 $$
 \HH^1(X,C^\bullet(T'',T'))\otimes H^0(\PP,\cO_{\PP}(1))\cong \WW\otimes \WW^*
 \cong \End(\WW).
 $$
One can show that the identity element in $\End(\WW)$ defines the
universal extension.

To prove the relative version stated in the proposition,
we consider the spectral sequence
 $$
 H^p(B,\HH_\pi^q(C^\bullet(\cT'',\cT')))\Rightarrow
 \HH^{p+q}(B\x X,C^\bullet(\cT'',\cT'))
 $$
relating relative and global hypercohomology groups. Since
$\HH_\pi^0(C^\bullet(\cT'',\cT'))=0$, the induced map
 $$
 \HH^{1}(B\x X,C^\bullet(\cT'',\cT')) \to
 H^0(B,\HH_\pi^1(C^\bullet(\cT'',\cT')))
 $$
is an isomorphism. Similarly, if $P:=\PP W^+$, we have an
isomorphism
 $$
 \HH^{1}( P\x X,C^\bullet(\cT'',\cT'\otimes \cO_P(1))) \cong
 H^0(P,\HH_\nu^1(C^\bullet(\cT'',\cT'\otimes \cO_P(1)))),
 $$
where $\nu: P\x X \to P$ is the canonical projection.

Now, write $p:P \to B$ for the projection. Again omitting
pull-backs when convenient, we have that the  image of the
identity of $W^+$ under the canonical isomorphisms
 \begin{align*}
 H^0(B,\End \, W^+)=H^0(B,W^+\otimes (W^+)^*) &= H^0(B,W^+\otimes p_*\cO_P(1))\\
                                    &= H^0(B,p_* (p^*W^+\otimes \cO_P(1)))\\
 &= H^0(P, p^*(\HH_\pi^1(C^\bullet(\cT'',\cT')))\otimes \cO_P(1))\\
 &= H^0(P,\HH_\nu^1(C^\bullet(\cT'',\cT'\otimes \cO_P(1))))
 \end{align*}
is a nonvanishing section defining the universal extension
(\ref{universal-ext+}). A technical ingredient in proving the
universal property is the commutation of $\HH^{1}(B\x
X,C^\bullet(\cT'',\cT'))$ with base change (see \cite{L} for
details on the analogous situation of extensions of sheaves).
\end{proof}

\subsection{Flips} \label{sec:flips}
Now we assume that $\cN_{\s_c^+}$, $\cN_{\s_c^-}$, $\cN'_{\s_c}$
and $\cN''_{\s_c}$ are smooth. Also assume that
$\HH^2(C^\bullet(T'',T'))=0$ and $\HH^2(C^\bullet(T',T''))=0$ for
every $(T',T'')\in \mathcal{N}_{\sigma_c}'\times
\mathcal{N}_{\sigma_c}''$. This will always be the case in our
applications. In order to relate $\cN_{\s_c^-}$ and $\cN_{\s_c^+}$
we have to blow up $\mathcal{N}_{\s_c^\pm}$ along
$\mathcal{S}_{\s_c^\pm}$. For this it is necessary to study the
normal bundle to $\mathcal{S}_{\s_c^\pm}= \PP W^\pm$ in
$\mathcal{N}_{\s_c^\pm}$.

\begin{prop} \label{prop:normal-bundle}
Let  $p:\PP W^\pm \to B$ be the natural  projection and  $j:\PP
W^\pm\inc \mathcal{N}_{\s_c^\pm}$ be the natural inclusion. Then
there is an  exact sequence
  $$
   0\too  T\PP W^\pm
  \too j^\ast T\cN_{\s_c^\pm} \too p^* W^\mp \ox
  \mathcal{O}_{\PP W^\pm}(-1)\too 0,
  $$
and hence, the normal bundle to $\mathcal{S}_{\s_c^\pm}= \PP
W^\pm$ in $\mathcal{N}_{\s_c^\pm}$ is isomorphic to $p^* W^\mp \ox
\mathcal{O}_{\PP W^\pm}(-1)$.
\end{prop}

\begin{proof}
We consider the case of $\cN_{\s_c^+}$ --- the  case of
$\cN_{\s_c^-}$ is analogous. Over $\cN_{\s_c^+}\x X$ there is a
universal triple $\cT=(\cE_1, \cE_2,\Phi)$, whose restriction to
$\cS_{\s_c^+}=\PP W^+$ is the universal extension $\cT^+$ in
(\ref{universal-ext+}). The tangent bundle of  $\cN_{\s_c^+}$ is
given by the  relative $\HH^1$  of the complex
 $$
  \PH (\cE_1,\cE_1) \oplus \PH(\cE_2,\cE_2) \too
  \SPH(\cE_2,\cE_1(D))
 $$
over $\mathcal{N}_{\s_c^+}\x X$ with respect to the natural
projection $\mathcal{N}_{\s_c^+} \x X \to \mathcal{N}_{\s_c^+}$.

Denote $\cE_i'(1)= \cE_i' \otimes \cO_{\PP W^+}(1)$
and define
 $$
 \PH_U(\cE_i,\cE_i):=\ker \left(\PH(\cE_i,\cE_i)\to
 \PH(\cE_i'(1),\cE_i'') \right)
 $$
and
 $$
 \SPH_U(\cE_2,\cE_1(D)):=\ker \left(\SPH(\cE_2,\cE_1(D))\to
 \SPH(\cE_2'(1),\cE_1''(D))\right).
 $$
The tangent bundle of  $\PP W^+$ is the relative $\HH^1$  with
respect to the projection  $\PP W^+ \x X\to \PP W^+$ of the middle
complex in the following exact sequence of complexes
 $$
 \begin{array}{ccccc}
  \PH (\cE_1'',\cE_1'(1)) \oplus
  \PH(\cE_2'',\cE_2'(1))  &\too& \SPH(\cE_2'',\cE_1'(1)(D)) \\
  \downarrow & & \downarrow \\
 \PH_U(\cE_1,\cE_1) \oplus \PH_U(\cE_2,\cE_2) &\too& \SPH_U(\cE_2,\cE_1(D)) \\
  \downarrow & & \downarrow \\
 \PH (\cE_1'(1),\cE_1'(1))
  \oplus \PH(\cE_2'(1),\cE_2'(1))\oplus &  & \SPH(\cE_2'(1),\cE_1'(1)(D)) \oplus \\
  \oplus \PH (\cE_1'',\cE_1'') \oplus \PH(\cE_2'',\cE_2'') &\too&
  \ \oplus \SPH(\cE_2'',\cE_1''(D)).
\end{array}
 $$
Note that when passing to cohomology, this gives us the exact
sequence
 $$
  0\longrightarrow T_V \PP W^+ \longrightarrow T \PP W^+ \longrightarrow T
  (\cN_{\sigma_c}' \times \cN_{\sigma_c}'') \longrightarrow 0,
 $$
where $T_V \PP W^+ \cong p^* W^+ (1) / \cO_{\PP W^+}$ is the
vertical tangent bundle.

Therefore the normal bundle to $\PP W^+$ is the relative $\HH^1$
of the quotient complex in the following exact sequence of
complexes
 $$
 \begin{array}{ccccc}
 \PH_U (\cE_1,\cE_1) \oplus \PH_U(\cE_2,\cE_2) &\too& \SPH_U(\cE_2,\cE_1(D)) \\
  \downarrow & & \downarrow \\
 \PH (\cE_1,\cE_1) \oplus \PH(\cE_2,\cE_2) &\too& \SPH(\cE_2,\cE_1(D)) \\
  \downarrow & & \downarrow \\
 \PH (\cE_1' (1),\cE_1'')
\oplus \PH(\cE_2'(1),\cE_2'') &\too&
\SPH(\cE_2'(1),\cE_1''(D)),
\end{array}
 $$
and it is hence isomorphic to $p^*W^-\otimes \cO_{\PP W^+}(-1)$.
\end{proof}

In particular, we conclude  that the embedding $\PP W^\pm\inc
\cN_{\s_c^\pm}$ is smooth. So we can blow-up $\cN_{\s_c^ \pm}$ along
$\PP W^ \pm$ to get $\widetilde{\cN}_{\s_c^\pm}$ with exceptional
divisor $E_\pm\subset \widetilde{\cN}_{\s_c^\pm}$ such that
 $$
 E_\pm = \PP W^\pm \x_B \PP W^\mp \, .
 $$
Note that ${\mathcal{O}}_{E_\pm}(E_\pm)={\mathcal{O}}_{\PP
W^+}(-1)\otimes {\mathcal{O}}_{\PP W^-}(-1)$, by Proposition
\ref{prop:normal-bundle}.

\begin{prop}\label{prop:5.12}
 There is a natural isomorphism
 $\widetilde{\cN}_{\s_c^+}\cong \widetilde{\cN}_{\s_c^-}$.
\end{prop}

\begin{proof}
Let $\cT$ be the universal triple over $\cN_{\s_c^+} \times X$. By
Proposition \ref{prop:5.6}, the restriction of $\cT$ to $\PP W^+
\times X$ lies in the universal extension
 \begin{equation}\label{eqn:exacta0}
  0\to \cT'\otimes \cO_{\PP W^+}(1)\to \cT|_{\PP W^+
 \times X}\to \cT'' \to 0.
 \end{equation}

Now pull back $\cT$ by the blow-up map $q:\widetilde{\cN}_{\s_c^+}
\to \cN_{\s_c^+}$. Consider the composition $q^*\cT \to q^*
\cT|_{E_+ \times X} \to q^*i_* \cT''$, where $i:\PP W^+ \x X\inc
\cN_{\sigma_c^+}\x X$. Since $q^*i_* \cT''$ is a triple formed by
two bundles supported on a divisor, the kernel is a triple (i.e.\
it is formed by two bundles, and not just coherent sheaves).
Define the triple $\hat{\cT}$ on $\widetilde{\cN}_{\s_c^+} \times
X$ by the exact sequence
  \begin{equation}\label{eqn:exacta}
  0 \longrightarrow \hat{\cT}\otimes \cO_{\PP W^+}(1) \longrightarrow
  q^*\cT \longrightarrow q^*i_*\cT'' \longrightarrow 0.
  \end{equation}
(This is called an {\em elementary transformation}.)

Let us see that all the triples in the family $\hat{\cT}$ are
$\s_c^-$-stable. Therefore this defines a map
$\widetilde{\cN}_{\s_c^+} \to \cN_{\s_c^-}$. Obviously, off
$E_+\times X$, $\hat{\cT} \cong \cT$ is the family parametrising
triples which are $\s_c^+$-stable and $\s_c^-$-stable at the same
time. Tensoring the  exact sequence \eqref{eqn:exacta} with
$\cO_{E_+\times X}$, we obtain, over $E_+\times X$,
 $$
  0 \to \underline{\Tor}(q^*i_*\cT'', \cO_{E_+\times X})
  \to \hat\cT|_{E_+\times X} \otimes \cO_{\PP W^+}(1)
  \to q^*\cT|_{E_+\times X} \to q^*i_* \cT'' \to  0 .
 $$
Since $\underline{\Tor}(q^*i_*\cT'', \cO_{E_+\times X})=\cT''\ox
\cO_{E_+\times X}(-E_+\times X)=\cT''\ox \cO_{\PP W^-}(1)\ox
\cO_{\PP W^+}(1)$, and also using (\ref{eqn:exacta0}), we get a
triple
 \begin{equation}\label{eqn:exacta2}
  0 \to \cT'' \ox \cO_{\PP W^-}(1) \to \hat\cT|_{E_+\times X} \to
  \cT' \to 0\, .
 \end{equation}
We have to check that all extensions in the family
(\ref{eqn:exacta2}) are non-trivial. For this, restrict to the fibre $\PP
W^+\x \PP W^-$ over a point $b\in B$. This corresponds to fixing
some specific triples $T'$ and $T''$. We have an exact sequence
 $$
   0 \to T'' \ox \cO_{\PP W^-_b}(1) \to \hat\cT|_{\PP
 W^+_b\x \PP W^-_b \times X} \to T' \to 0\, .
 $$
This extension class is parametrised by
  $$
   \Ext^1(T', T'' \ox \cO_{\PP W^-_b}(1)) =  W^-_b \ox H^0(\cO_{\PP
   W^-_b}(1)) = \End(W^-_b).
  $$
Moreover the linear group $GL (W^-_b)$ acts on $W^-_b$. The
extension class is invariant by this action, therefore it is a
linear multiple of the identity. Letting $b$ move in $B$ we have a
section of $\End (W^-)$. Since this is a multiple of the identity,
it lives in $\cO \cdot \Id \subset \End(W)$, therefore it is a
{\em constant\/} multiple of the identity. This cannot be
constantly zero for, otherwise, it would be $\hat\cT|_{E_+\times
X}=\cT' \oplus (\cT'' \ox \cO_{\PP W^-}(1))$. Then $\Hom
(\hat\cT,q^*i_*\cT'' \ox \cO_{\PP W^-}(1)) \neq 0$. Hence
(\ref{eqn:exacta}) would imply that the map
  \begin{equation}\label{eqn:inj-not}
  \Ext^1(\cT''\ox \cO_{\PP W^+}(-1),\cT'' \ox \cO_{\PP W^-}(1))  \to
  \Ext^1(\cT\ox \cO_{\PP W^+}(-1),\cT'' \ox \cO_{\PP W^-}(1))
  \end{equation}
is not injective. On the other hand, using the proyection $\pi:
\cN_{\sigma_c^+}\times X \to \cN_{\sigma_c^+}$ in
(\ref{eqn:exacta0}) we have  that
  \begin{align*}
  & \Hom_{\pi}(\cT'',\cT'')  \cong \Hom_{\pi}(\cT,\cT'') \, ,\\
  & \Ext^1_\pi (\cT'',\cT'') \hookrightarrow \Ext^1_\pi (\cT,\cT'') \,  ,
  \end{align*}
as bundles over $\cN_{\sigma_c^+}$. Twisting by $\cO_{\PP
W^+}(1)\ox \cO_{\PP W^-}(1)$ and using the spectral sequence
 $$
 H^p(\cN_{\sigma_c^+},\Ext^q_\pi ( \cdot,\cdot ))\Rightarrow
 \HH^{p+q}(\cN_{\sigma_c^+}\x X,C^\bullet(\cdot,\cdot))
 $$
we have that (\ref{eqn:inj-not}) is injective, giving a
contradiction.

Hence the extension class of (\ref{eqn:exacta2}) is a non-zero
multiple of the identity. This gives a map $\tilde\cN_{\s_c^+}\to
\cN_{\s_c^-}$ which restricts to $E_+$ as the natural projection
on the second factor $E_+\iso \PP W^+\x_B \PP W^- \to \PP W^-$.
Analogously we obtain a map $\tilde\cN_{\s_c^-} \to \cN_{\s_c^+}$.
So there are two {\em injective} maps $\tilde\cN_{\s_c^\pm} \to
\cN_{\s_c^+} \x\cN_{\s_c^-}$. Their images are both the closures
of the image of $\tilde\cN_{\s_c^\pm} \setminus E_\pm$, which are
the same. So they coincide.
\end{proof}

\begin{rem} \label{rem:5.8}
Let $\sigma_c$ be a critical value. If $w^+=\rk (W^+)>0$ and
$w^-=\rk (W^-)>0$ then the moduli spaces $\cN_{\sigma^-_c}$ and
$\cN_{\sigma^+_c}$ are birational. This is true because $w^+ + w^-
+ \dim B-1=\dim \cN_{\sigma_c^\pm}$ by Corollary
\ref{cor:chi-relation}, and the flip loci $\cS_{\s_c^\pm}$ have
dimension $w^{\pm} + \dim B -1$.
\end{rem}

\section{Parabolic triples with $r_1=2$ and $r_2=1$}
\label{sec:thaddeus-program}

In this section we use the results of Section
\ref{sec:critical-values-flips} to compute the Poincar\'e
polynomial of the moduli space of parabolic triples $\cN_{\sigma}$
for the case $r_1=2$ and $r_2=1$ and for non-critical values of
$\s$. We are studying triples of the form $\phi:L \to E_1(D)$,
where $L$ is a parabolic line bundle with $\deg (L)=d_2$ and
weights $\alpha(p)$, and $E_1$ is a parabolic rank $2$ bundle with
$\deg (E_1)=d_1$ and weights $\beta_1(p)< \beta_2(p)$. By Theorem
\ref{thm:smoothdim} (v), the moduli space of stable elements in
$\cN_{\sigma}$ is smooth. Moreover, applying the exact sequence
(\ref{eq:long-exact-extension-complex}), one can easily show that
we are in the situation given in Proposition \ref{vanishing-h2}.

\subsection{Flips}\label{subsec:thaddeus-program.1}

By Subsection \ref{sec:critical-values}, there are the following
three possibilities for the existence of critical values:

\begin{itemize}
 \item $r_1'=1$ and $r_2'=0$. Then the subtriple $T'$ is of the form
 $0\to M(D)$, where $M$ is a line bundle of degree $d_M$.
 Since $M$ inherits weights from $E_1$,
 there is a function $\e=\{\e(p)\}_{p\in D}$, which assigns to each
 $p\in D$ a number $\e(p)\in \{1,2\}$ such that the weight of $M$ at $p$
 is $\b_{\e(p)}(p)$. We have an exact sequence of triples
 $$
  \begin{array}{ccccc}
  0 &\too& L &\too& L \\
  \downarrow & & \downarrow & & \downarrow \\
  M(D) & \too & E_1(D) & \too & F(D) \, .
  \end{array}
  $$
 The quotient triple is of the form $L\to F(D)$, where $F$ is a
  parabolic line bundle of degree $d_1-d_M$ and weights
  $\b_{\varsigma(p)}(p)$, with $\varsigma(p)=3-\e(p)$. Note that
  $\{\b_{\e(p)}(p), \b_{\varsigma(p)}(p)\}= \{\b_1(p), \b_2(p)\}$.
  By (\ref{eq:strict-alpha-ss}), the critical value is
 \begin{equation}\label{eqn:sc-july}
 \s_c= 3d_M -d_1 -d_2 + \sum_p \left(2\b_{\e(p)}(p)-\a(p) -\b_{\varsigma(p)}(p)\right).
 \end{equation}
  As described in Subsection \ref{sec:flips},
  this defines the subspace $\mathcal{S}_{{\sigma_c}^+}=\PP W^+_{\s_c}$,
  where
   $$
   W^+_{\s_c} \longrightarrow B_{\s_c} =\cN'_{\s_c} \times \cN''_{\s_c}
  $$
(we shall make the dependence on $\s_c$ explicit in this section,
since we shall be working with various flip loci simultaneously).
The moduli space parametrizing the possible parabolic line bundles
$M$ with fixed weights $\beta_{\e(p)}(p)$ is $\cN'_{\s_c} =
\Jac^{d_M} X$. The moduli space parametrizing triples of the form
$L\to F(D)$, which are parabolic line bundles with fixed weights
is $\Jac^{d_2} X \x S^{N}X$, where $N=\deg \SPH(L, F(D))$. To
compute this we use the following.

\begin{lem}\label{lem:sphom-phi_l}
Let $L_1, L_2$ be two parabolic line bundles with weights
$\alpha_{L_1}(p)$ and $\alpha_{L_2}(p)$, respectively. Then
 $$
  \SPH(L_1,L_2 \otimes K(D))  \cong \Hom(L_1,L_2 \otimes K(S))\, ,
 $$
where $S = \{p \in D \;|\; \alpha_{L_1}(p) <\alpha_{L_2}(p) \}$.
\end{lem}

\begin{proof}
By definition, a strongly parabolic map $\Phi: L_1\to L_2$
satisfies
\begin{displaymath}
  \mathrm{Res}_{p}\Phi = 0 \;\iff\; \alpha_{L_1}(p) \geq
  \alpha_{L_2}(p).
\end{displaymath}
{}From this the result is clear.
\end{proof}

In our case, we introduce the following notations:
 \begin{equation}\label{eqn:6.2.5}
 \begin{split}
 S_1 &= \{ p \in D \, | \, \alpha(p)< \beta_{\varsigma(p)}(p) \}\, ,\\
 S_2 &= \{ p \in D \, | \, \alpha(p)< \beta_{\e(p)}(p) \}\, ,\\
 S_3 &= \{ p \in D \, | \,  \beta_{\e(p)}(p) <\beta_{\varsigma(p)}(p)
 \}\, , \\
 s_1 &=\# S_1, \, s_2=\# S_2,\,  s_3=\# S_3\, .
 \end{split}
 \end{equation}
Then
 \begin{equation} \label{eqn:N-july}
 \begin{split}
 N &= \deg \SPH(L, F(D))= \deg \Hom(L,F(S_1)) = \\ &= \deg (F)-\deg (L)  + s_1  = d_1- d_2-d_M  +
 s_1.
 \end{split}
 \end{equation}

Now $\PP W^+_{\s_c}$ is a projective fibration over $B_{\s_c}$
with fibres projective spaces of dimension $w^+_{\s_c} - 1$. By
Proposition \ref{prop:chi(T'',T')} and Corollary
\ref{dimension-ext1},
 \begin{equation} \label{eqn:D+july}
 \begin{split}
 w^+_{\s_c} &= \dim \Ext^1(T'',T') =  - \chi(T'',T') \\
 &=-\chi(\PH(F,M)) +\chi(\SPH(L,M(D))) \\ &=
 -\chi(\Hom(F,M(-S_3))) +\chi(\Hom(L,M(S_2))) \\
 &=  d_1-d_2 -d_M  + s_2 + s_3 \, .
 \end{split}
 \end{equation}

\item $r_1'=1$, $r'_2=1$. Then the subtriple $T'$ is of the form
 $L\to F(D)$ and the quotient triple is of the form
 $0\to M(D)$, yielding an exact sequence
  $$
  \begin{array}{ccccc}
  L &\too& L &\too& 0 \\
  \downarrow & & \downarrow & & \downarrow \\
  F(D) & \too & E(D) & \too & M(D) \, ,
  \end{array}
  $$
 where $M$ is a line bundle of degree $d_M$
 and weights $\b_{\e(p)}(p)$, for some $\e=\{\e(p)\}_{p\in D}$,
 and $F$ is a
  parabolic line bundle of degree $d_1-d_M$ and weights
  $\b_{\varsigma(p)}(p)$, with $\varsigma(p)=3-\e(p)$.
  The critical value is again given by (\ref{eqn:sc-july}).
  These extensions define the subspace
  $\mathcal{S}_{{\sigma_c}^-}=\PP W^-_{\s_c}$,
  where
   $$
   W^-_{\s_c} \longrightarrow B_{\s_c}=\cN'_{\s_c} \times \cN''_{\s_c}=
   \Jac^{d_M} X \times  \Jac^{d_2} X \times S^N X\, ,
   $$
  with $N$ as in (\ref{eqn:N-july}). Now $\PP W^-_{\s_c}$ is a projective
  fibration over $B_{\s_c}$ with fibres projective spaces of dimension $w^-_{\s_c} - 1$
  where
\begin{equation} \label{eqn:D-july}
 \begin{split}
 w^-_{\s_c} &=   h^1(T',T'')= - \chi(T',T'') = -\chi(\PH(M,F))= \\
 &= -\chi (\Hom(M,F (-(D-S_3)))) = 2d_M-d_1 + g-1 + n-s_3 \, .
 \end{split}
 \end{equation}

 \item $r'_1=2$, $r_2'=0$. Then
the triples $T$ are extensions of $0\to E_1$ by $L\to 0$. The
critical value is $\s_c=\pmu_1-\pmu_2=\s_m$, which is the minimum
possible value for the parameter $\s$. At this value, the moduli
space $\cN_{\s_m^-}=\emptyset$ and $\cN_{\s_m^+}=\cS_{\s_m^+}$.
This can be described explicitly as a projective fibration over a
product of a Jacobian and a moduli space of rank $2$ stable
parabolic bundles, but we will not go into this since we shall not
use it.
\end{itemize}

\begin{rem}
Since we are taking generic values for the weights, the values of
$\s_c$ are distinct for the different choices of $d_M$ and $\e$.
The genericity condition was necessary in Section
\ref{sec:critical-values-flips} to have smooth flip loci, and this
is essentially due to the fact that at a critical value, the
Jordan-H\"older filtrations are of length at most two. In the case
we treat here, $r_1=2$ and $r_2=1$, the Jordan-H\"older
filtrations are of length at most two even with non-generic
weights, because the ranks are too small. Therefore the
computations of this section work as well for the case of {\em
distinct\/} non-generic weights. Of course, we shall need the
genericity of weights at many other places in the coming sections.
\end{rem}

\begin{rem} \label{rem:N_L}
Let  $\sigma_L$ be the largest critical value. This means that the
moduli space $\cN_{\sigma_L^+}=\emptyset$, i.e.,
$\cS_{\sigma_L^-}=\cN_{\sigma_L^-}$. Since
 $$
 \dim B_{\s_c}+ w^+_{\s_c} +w^-_{\s_c} -1=\dim \cN_{\s_c^\pm},
 $$
by Corollary \ref{cor:chi-relation}, we have that $w^+_{\s_L}=0$.
By (\ref{eqn:D+july}), this corresponds to the case $s_3=0$ and
$d_1-d_2-d_M+s_2=0$, i.e., when $L\cong F(S_2)\subset F(D)$. In
this case $\cN_{\s_L^-}$ equals $\PP W^-_{\s_L} \to B_{\s_L}$.

Also $\s_L<\s_M$; in general they are not equal. The value of
$\s_L$ obtained in (\ref{eqn:sc-july}) is always slightly smaller
than that of $\s_M$ in (\ref{alpha-bounds-bigM}).
\end{rem}

\subsection{Poincar{\'e} polynomial of moduli of triples} \label{sec:betti-triples}

Let $\s_c$ be a critical value as in Subsection
\ref{subsec:thaddeus-program.1} with the only condition $\s_c\neq
\s_m$. Then we have that
 \begin{equation}\label{eqn:Pt-a}
  P_t(\cN_{\s_c^-}) -P_t(\cN_{\s_c^+}) =P_t( \PP W^-_{\s_c}) - P_t(\PP
  W^+_{\s_c})\, .
 \end{equation}
Note that this formula also holds when $w^+_{\s_c}=0$ or
$w^-_{\s_c}=0$. For instance, if $w^+_{\s_c}=0$ then
$\cS_{\s_c^+}=\emptyset$ and $\cS_{\s_c^-}=\PP W^-_{\s_c}$ is of
the same dimension as $\cN_{\s_c^-}$, hence it is a component of
it. So \eqref{eqn:Pt-a} holds. In particular we can use
\eqref{eqn:Pt-a} for $\s_c=\s_L$ (see Remark \ref{rem:N_L}). But
we cannot use it for $\s_c=\s_m$ (see Remark \ref{rem:Pt(Ns)}).

\begin{thm} \label{thm:Pt(Ns)} Let $\s>\s_m$ be a non-critical value.
For any $\e=\{\e(p)\}_{p\in
D}$, $\e(p)\in \{1,2\}$, let $s_1,s_2,s_3$ be given by
(\ref{eqn:6.2.5}) and
 $$
 \bdm=\left[\frac13 \left(d_1+d_2 +\sum \left(
 \a(p) +\b_{\varsigma(p)}(p) -
 2\b_{\e(p)}(p)\right)+\s\right)\right]+1\, ,
 $$
 where $\varsigma(p)=3-\e(p)$ and $[x]$ is the integer part of
 $x$. Then $P_t(\cN_\sigma)$ equals
 $$
 \begin{aligned}
 \sum_{\e} \Coeff_{x^0} \bigg( & \frac{(1+t)^{4g}(1+xt)^{2g}
 t^{2d_1-2d_2+2s_2+2s_3-2\bdm}x^{\bdm}}{(1-t^2)(1-x)(1-xt^2)(1-t^{-2}x)
 x^{d_1-d_2+s_1}} \\
 &-\frac{(1+t)^{4g}(1+xt)^{2g}
 t^{-2d_1+2g-2 +2n-2s_3+4\bdm}x^{\bdm} }{(1-t^2)(1-x)(1-xt^2)(1-t^4x)
 x^{d_1-d_2+s_1}}\bigg)\, .
 \end{aligned}
 $$
\end{thm}

\begin{proof}
{}From (\ref{eqn:Pt-a}), we have that
 \begin{eqnarray*}
  P_t(\cN_{\s}) &=& \sum_{\s_c>\s} \left( P_t( \PP W^-_{\s_c})
 - P_t(\PP W^+_{\s_c}) \right) \\
 &=&  \sum_{\s_c>\s} \left(P_t(\PP^{w^-_{\s_c}-1}) -
 P_t(\PP^{w^+_{\s_c}-1}) \right) P_t(B_{\s_c}) \\
 &=&  \sum_{\s_c>\s} \bigg(\frac{1-t^{2w^-_{\s_c}}}{1-t^2} -
 \frac{1-t^{2w^+_{\s_c}}}{1-t^2}\bigg) P_t(\Jac \,X)^2 P_t(\Sym^N X) \\
 &=& \sum_{\s_c>\s} \frac{t^{2w^+_{\s_c}}-t^{2w^-_{\s_c}}}{1-t^2} \ (1+t)^{4g}
 \Coeff_{x^0}\left(
 \frac{(1+xt)^{2g}}{(1-x)(1-xt^2)x^{N}}
 \right)  \qquad \qquad \hbox{by \cite{McD}}\\
 &=& \sum_{\s_c>\s} \frac{t^{2d_1-2d_2-2d_M+2s_2+2s_3}-
 t^{4d_M-2d_1+2g-2 +2n-2s_3}}{1-t^2} \ (1+t)^{4g} \cdot \\
  & & \qquad \cdot \Coeff_{x^0}\left(
 \frac{(1+xt)^{2g}}{(1-x)(1-xt^2)x^{d_1-d_2-d_M+s_1}}
 \right) \qquad \hbox{by (\ref{eqn:N-july}), (\ref{eqn:D+july})
 and (\ref{eqn:D-july})} \\
 &=& \sum_{\e} \Coeff_{x^0} \bigg(  \frac{(1+t)^{4g}(1+xt)^{2g}
 t^{2d_1-2d_2+2s_2+2s_3}}{(1-t^2)(1-x)(1-xt^2)x^{d_1-d_2+s_1}}
 \sum_{d_M | \s_c>\s} t^{-2d_M}x^{d_M}   \\ & &  \qquad -
 \frac{(1+t)^{4g}(1+xt)^{2g}
 t^{-2d_1+2g-2 +2n-2s_3}}{(1-t^2)(1-x)(1-xt^2)x^{d_1-d_2+s_1}}
 \sum_{d_M | \s_c>\s} t^{4d_M}x^{d_M}
 \bigg)\, .
 \end{eqnarray*}

The condition for $d_M$ is
 $$
 \s_c= 3d_M -d_1 -d_2 + \sum \left( 2\b_{\e(p)}(p)-\a(p) -\b_{\varsigma(p)}(p)
 \right) >\s\, ,
 $$
which translates into
  $$
  d_M > \frac13 \left(d_1+d_2 +\sum \left( \a(p) +\b_{\varsigma(p)}(p)
  - 2\b_{\e(p)}(p)\right)+\s\right).
  $$
Since $\s$ is not a critical value, we cannot have equality, so
the right hand side is not an integer. The inequality becomes
$d_M\geq \bdm$, with $\bdm$ as in the statement. Now
 \begin{eqnarray*}
 \sum_{d_M=\bdm}^\infty t^{-2d_M}x^{d_M}&=&
 \frac{t^{-2\bdm}x^{\bdm}}{1-t^{-2}x}\, , \\
  \sum_{d_M=\bdm}^\infty t^{4d_M}x^{d_M} &=&
 \frac{t^{4\bdm}x^{\bdm}}{1-t^4x}\,  .
 \end{eqnarray*}

So finally
 \begin{eqnarray*}
  P_t(\cN_{\s}) &=& \sum_{\e} \Coeff_{x^0} \bigg(  \frac{(1+t)^{4g}(1+xt)^{2g}
 t^{2d_1-2d_2+2s_2+2s_3}t^{-2\bdm}x^{\bdm}}{(1-t^2)(1-x)(1-xt^2)(1-t^{-2}x)
 x^{d_1-d_2+s_1}} \\
 & & -\frac{(1+t)^{4g}(1+xt)^{2g}
 t^{-2d_1+2g-2 +2n-2s_3} t^{4\bdm}x^{\bdm} }{(1-t^2)(1-x)(1-xt^2)(1-t^4x)
 x^{d_1-d_2+s_1}} \bigg)\, .
 \end{eqnarray*}
\end{proof}

\begin{rem} \label{rem:Pt(Ns)}
The formula in this theorem only works for $\sigma>\sigma_m$. For
$\sigma<\sigma_m$, $\cN_\s$ is empty, but the formula above does
not give zero for such values.
\end{rem}

\section{Critical submanifolds of type $(1,1,1)$}\label{sec:(1,1,1)}

\subsection{Description of the critical submanifolds}
\label{sec:(1,1,1)-description}

In this section we consider the critical points of the Bott--Morse
function $f$ represented by parabolic Higgs bundles $(E,\Phi)$ of
type $(1,1,1)$, i.e., of the form $E=L_1 \oplus L_2 \oplus L_3$
where $L_l$ are parabolic line bundles, i.e., line bundles with
weights at the points $p\in D$.  We denote the (fixed) weights of
$(E,\Phi)$ at $p \in D$ by $0 \leq \alpha_1(p) < \alpha_2(p) <
\alpha_3(p) < 1$.  Each possible choice of the distribution of
these weights among the line bundles $L_l$ is given by a
permutation $\varpi_p \in S_3$ such that the weight on the fibre
$L_{l,p}$ at $p$ is $\alpha_{\varpi_p(l)}(p) =
\alpha_{\varpi(l)}(p)$ for $l=1,2,3$. The map $\Phi$ decomposes as
strongly parabolic maps $\Phi_1:L_1\to L_2\otimes K(D)$ and
$\Phi_2:L_2\to L_3\otimes K(D)$.

We define
\begin{align*}
   d_l &= \deg (L_l) \quad \text{for $\, l=1,2,3$}, \\
   m &= d_{1} + d_2.
\end{align*}
We shall choose to describe the topological data $(d_1,d_2,d_3)$
using the parameters $(d_1,m,\Delta)$, where $\Delta = d_1 + d_2 +
d_3$. In terms of this data we have $d_2=m-d_1$ and
$d_3=\Delta-m$.
We also introduce the notation
\begin{equation} \label{eq:def-FG}
\begin{split}
   F(\alpha,\varpi) &= \sum_{p \in D}
   (\alpha_1(p)+\alpha_2(p)+\alpha_3(p)-3\alpha_{\varpi(3)}(p)), \\
   G(\alpha,\varpi) &=  \sum_{p \in D}
   (2\alpha_1(p)+2\alpha_2(p)+2\alpha_3(p)
   - 3\alpha_{\varpi(2)}(p) -3\alpha_{\varpi(3)}(p)).
\end{split}
\end{equation}

\begin{prop}\label{prop:par-stability111}
  A parabolic Higgs bundle $(L_1 \oplus L_2 \oplus L_3,\Phi)$ of type
  $(1,1,1)$ is stable if and only if the maps $\Phi_1$ and $\Phi_2$
  are non-zero and, furthermore,
  \begin{align*}
    3m     &> 2\Delta - F(\alpha,\varpi),\\
    3d_{1} &> \Delta - G(\alpha,\varpi).
  \end{align*}
\end{prop}

\begin{proof}
  It is immediate from Propositions
  \ref{prop:fixed=vhs} and \ref{prop:vhs-stab} that a
  parabolic Higgs bundle of type $(1,1,1)$ is stable if and only if  the
  conditions $\pmu (L_3) < \pmu(E)$ and $\pmu (L_2 \oplus L_3) <
  \pmu(E)$ hold and $\Phi_l \neq 0$ for $l=1,2$.  {}From this
  we obtain the characterization given in the Proposition by
  calculating the relevant parabolic degrees:
  \begin{align*}
    \pdeg(E) &= \Delta + \sum_{p\in
      D}(\alpha_1(p)+\alpha_2(p)+\alpha_3(p))\ , \\
    \pdeg(L_3) &= d_3 + \sum_{p \in D} \alpha_{\varpi(3)}(p) \\
    &= \Delta -m
    + \sum_{p \in D} \alpha_{\varpi(3)}(p)\ , \\
    \pdeg(L_2 \oplus L_3) &= d_2 + d_3
    + \sum_{p \in D} (\alpha_{\varpi(2)}(p)+\alpha_{\varpi(3)}(p)) \\
    &= \Delta - d_1
    + \sum_{p \in D} (\alpha_{\varpi(2)}(p)+\alpha_{\varpi(3)}(p))\ .
  \end{align*}
  Thus we get
  \begin{align*}
    && \pmu (L_3) &< \pmu(E) \\
    &\iff&  3\Delta - 3m +3\sum_{p \in D}
    \alpha_{\varpi(3)}(p)
    &< \Delta + \sum_{p\in
      D}(\alpha_1(p)+\alpha_2(p)+\alpha_3(p)) \\
    &\iff& 3m & > 2\Delta - F(\alpha,\varpi)\ ,
  \end{align*}
 while
  \begin{align*}
  && \pmu (L_2 \oplus L_3) &< \pmu(E) \\
  \iff& & 3\Delta - 3d_1
  +3\sum_{p \in D} (\alpha_{\varpi(2)}(p)+\alpha_{\varpi(3)}(p))
  &< 2\Delta + 2\sum_{p\in
    D}(\alpha_1(p)+\alpha_2(p)+\alpha_3(p)) \\
  \iff& & 3d_{1} &> \Delta - G(\alpha,\varpi)\ .
  \end{align*}
\end{proof}

Denote by $\mathcal{N}_{(1,1,1)}(d_1,m,\varpi)$ the critical
submanifold of parabolic Higgs bundles of type $(1,1,1)$ with
invariants $(d_1,m)$ and weights given by $\varpi =
\{\varpi_p\}_{p \in D}$. Introduce the following notation
\begin{equation}\label{eqn:S,1,1,1}
\begin{split}
  S_1 &= \{p \in D \;|\; \alpha_{\varpi(1)}(p) >
  \alpha_{\varpi(2)}(p) \}, \\
  S_2 &= \{p \in D \;|\; \alpha_{\varpi(2)}(p) >
  \alpha_{\varpi(3)}(p) \}, \\
  s_1 &= \# S_1, \  s_2 = \# S_2\, .
\end{split}
\end{equation}
There is no risk of confusion with the sets $S_1$, $S_2$ defined
in Section \ref{sec:thaddeus-program}, since this notation will
only apply to this section. By Lemma \ref{lem:sphom-phi_l},
\begin{displaymath}
  \SPH(L_l,L_{l+1} \otimes K(D))
  \cong \Hom(L_l,L_{l+1} \otimes K(D-S_l))
\end{displaymath}
for $l=1,2$. We let for $l=1,2$,
\begin{displaymath}
  m_{l} = \deg(\Hom(L_{l},L_{l+1} \otimes K(D-S_l))).
\end{displaymath}
Then
\begin{equation}\label{eq:def-m_i}
\begin{split}
  m_{1} &= d_{2} - d_{1} + 2g-2 + n - s_{1} \\
        &= m - 2d_{1} + n - s_{1} + 2g-2,  \\
  m_{2} &= d_{3} - d_{2} + 2g-2 + n - s_{2} \\
        &= \Delta - 2m + d_1 + n - s_2 + 2g-2.
\end{split}
\end{equation}

\begin{prop} \label{prop:boundsmm}
The critical submanifold $\mathcal{N}_{(1,1,1)}(d_1,m,\varpi)$ is
non-empty if and only if the following conditions are satisfied
\begin{align*}
    &&&& 3m &> 2\Delta - F(\alpha,\varpi),\\
    &&&& 3d_{1} &> \Delta - G(\alpha,\varpi), \\
    m_{1} &\geq 0 & &\iff&  2d_{1} -m &\leq n - s_{1} + 2g-2,\\
    m_{2}&\geq 0 & &\iff&   2m - d_1 &\leq \Delta  + n - s_2 + 2g-2,
\end{align*}
where $F$ and $G$ were defined in \eqref{eq:def-FG}. Moreover, the
map
\begin{align*}
  \mathcal{N}_{(1,1,1)}(d_1,m,\varpi)
  &\to \mathrm{Jac}^{d_1}(X) \times S^{{m}_1} X \times
  S^{{m}_2} X \\
  (L_1 \oplus L_2 \oplus L_3, \Phi_1, \Phi_2)
  &\mapsto (L_1,\mathrm{div}(\Phi_1),\mathrm{div}(\Phi_2))
\end{align*}
is an isomorphism.
\end{prop}

\begin{proof}
  Proposition~\ref{prop:par-stability111} shows that the conditions
  in the statement are necessary
  for $\mathcal{N}_{(1,1,1)}(d_1,m,\varpi)$ to be non-empty.

  Assume now that we are given $(d_1,m,\varpi)$ satisfying
  this conditions.  For any line bundle $L_1$ in
  $\mathrm{Jac}^{d_1}(X)$ and effective divisors $D_1 \in S^{{m}_1} X$
  and $D_2 \in S^{{m}_2} X$ we
  get line bundles $M_l = \mathcal{O}(D_l)$ with non-zero sections
  $\Phi_l$ determined up to multiplication by nonzero scalars for
  $l=1,2$.  We then obtain a parabolic Higgs bundle of type $(1,1,1)$
  by letting
  \begin{align*}
    L_2 &= L_1\otimes K^{-1}(S_1-D)\otimes M_1, \\
    L_3 &= L_2\otimes K^{-1}(S_2-D)\otimes M_2,
  \end{align*}
  and defining $\Phi$ to have components $\Phi_{1}$ and $\Phi_{2}$.
  Clearly this parabolic Higgs bundle has the desired invariants
  $(d_1,m,\varpi)$  and, if the conditions
  in the statement are satisfied, then
  Proposition \ref{prop:par-stability111} shows that it is indeed
  stable.

  It follows from this construction that the map given in the
  statement of the Proposition is surjective.  To see that it is
  injective, we note that taking non-zero scalar multiples of the Higgs
  fields $\Phi_1 \in H^0(L_1^{-1}\otimes L_2\otimes K(D-S_1))$ and
  $\Phi_2 \in H^0(L_2^{-1}\otimes L_3\otimes K(D-S_2))$ gives rise to
  isomorphic parabolic Higgs
  bundles of type $(1,1,1)$.  Thus the map given is, in fact, an
  isomorphism.
\end{proof}

\begin{cor}\label{cor:P-t-N-111}
  The Poincar\'e polynomial of the critical submanifold
  $\mathcal{N}_{(1,1,1)}(d_1,m,\varpi)$ is
  \begin{displaymath}
    P_{t}(\mathcal{N}_{(1,1,1)}(d_1,m,\varpi)) = (1+t)^{2g}
    \Coeff_{x^0y^0}\left(
    \frac{(1+xt)^{2g}}{(1-x)(1-xt^2)x^{m_{1}}}
    \cdot
    \frac{(1+yt)^{2g}}{(1-y)(1-yt^2)y^{m_{2}}}
    \right) .
  \end{displaymath}
\end{cor}

\begin{proof}
  Immediate from MacDonald's formula \cite{McD}
  for the Poincar\'e polynomial of
  a symmetric product of $X$.
\end{proof}

The total contribution to the Poincar\'e polynomial of
$\mathcal{M}$ from submanifolds of type $(1,1,1)$ is
  \begin{equation}\label{eq:contrib-11}
  P_t(\Delta,(1,1,1)) : =
  \sum_{d_{1},m,\varpi} t^{\lambda_{(d_1,m,\varpi)}}
    P_{t}(\mathcal{N}_{(1,1,1)}(d_1,m,\varpi)),
  \end{equation}
where $\lambda_{(d_1,m,\varpi)}$ is the index of the critical
submanifold $\mathcal{N}_{(1,1,1)}(d_1,m,\varpi)$ and the sum is
over all permutations $\varpi = \{\varpi_{p}\}_{p\in D} \in
(S_{3})^n$ and pairs of integers $(d_{1},m)$ such that the bounds
of Proposition \ref{prop:boundsmm} are satisfied.

\begin{lem}\label{lem:index-111}
  The index of the critical submanifold $\mathcal{N}_{(1,1,1)}(d_1,m,\varpi)$
  is
  \begin{displaymath}
    \lambda_{(d_1,m,\varpi)} =
    2\bigl(
    4g-4+n+s_{1}+s_{2}-\Delta+d_{1}+m
    \bigr).
  \end{displaymath}
\end{lem}

\begin{proof}
  The formula for the Morse index is given in Proposition
  \ref{prop-index}.

  We need to calculate, for each $p\in D$ the numbers $f_{p}$ and the
  dimensions of the spaces $P_{p}$ and $N_{p}$, which enter in this
  formula.

  Recall (from
  Proposition~\ref{prop:phb-moduli-dim}) that $f_{p} = (1/2)(r^2 -
  \sum_{i}m_{i}(p)^2)$.  Since the multiplicities are all $1$ and
  $r=3$ we get $f_{p} = (1/2)(9-(1+1+1)) = 3$ and hence
  \begin{displaymath}
    \sum_{p} f_{p} = 3n\ .
  \end{displaymath}

  For $p \in D$, the space $P_{p}(L_{l},L_{l})$ consists of the
  parabolic endomorphisms of $L_{l,p}$, so $\dim P_{p}(L_{l},L_{l}) =
  1$.  The space $N_{p}(L_{l},L_{l+1})$ is the space of strictly
  parabolic maps from $L_{l,p}$ to $L_{l+1,p}$ and hence
  \begin{displaymath}
    N_{p}(L_{l},L_{l+1}) =
    \begin{cases}
      0 &\text{if $\alpha_{\varpi(l)}(p)
         > \alpha_{\varpi(l+1)}(p)$,} \\
      \Hom(L_{l,p},L_{l+1,p}) &\text{otherwise.}
    \end{cases}
  \end{displaymath}
  Recalling from \eqref{eqn:S,1,1,1} the definition of
  $S_{l}$, it follows that
  \begin{displaymath}
    \dim N_{p}(L_{l},L_{l+1}) =
    \begin{cases}
      0 &\text{if $p \in S_{l}$,} \\
      1 &\text{if $p \in D-S_{l}$,}
    \end{cases}
  \end{displaymath}
  and thus
  \begin{displaymath}
    \sum_{p} \dim N_{p}(L_{l},L_{l+1}) = n - s_{l}\ .
  \end{displaymath}

  Substituting this in the formula for the Morse index we get
  \begin{align*}
    \lambda_{(d_1,m,\varpi)} &=
    r^2(2g-2) + 2\sum_{p}f_{p}
    + 2\sum_{l=1}^3 \Bigl(
    (1-g-n)\rk(L_{l})^2 + \sum_{p}\dim P_{p}(L_{l},L_{l})\Bigr)\\
    &\quad + 2\sum_{l=1}^2 \Bigl((1-g)\rk(L_{l})\rk(L_{l+1}) -
    \rk(L_{l})\deg(L_{l+1})
    + \rk(L_{l+1})\deg(L_{l}) \\
    &\quad - \sum_{p}\dim
    N_{p}(L_{l},L_{l+1})
    \Bigr) \\
    &= 9(2g-2) + 2\cdot 3n
    +2\bigl(3(1-g-n)+3n\bigr) \\
    &\quad +2\bigl(2(1-g)-d_{2} -d_{3} +d_{1} +d_{2}
    -(n-s_{1}+n-s_{2})\bigr) \\
    &=2(4g-4+n+s_{1}+s_{2}-\Delta+d_{1}+m)\ .
  \end{align*}

\end{proof}

\subsection{The sum for fixed $\varpi$.}

We shall now calculate the total contribution
\eqref{eq:contrib-11} to the Poincar\'e polynomial of
$\mathcal{M}$ from submanifolds of type $(1,1,1)$ in several
stages.  We begin by doing the sum over $(d_1,m)$ for a fixed
permutation $\varpi$.

\begin{lem}
  Let $\varpi = \{\varpi_{p}\}_{p\in D} \in (S_{3})^n$ be fixed.  Then
  $$
      \sum_{d_{1},m} t^{\lambda_{(d_1,m,\varpi)}}
      P_{t}(\mathcal{N}_{(1,1,1)}(d_1,m,\varpi))
      = \Coeff_{x^0y^0}\Psi(\varpi),
  $$
  where we have defined
  $$
      \Psi(\varpi) = \sum_{
        \substack{
          d_{1} \geq \bd \\ m \geq \bm}}
    t^{2(4g-4+n+s_{1}+s_{2}-\Delta+d_{1}+m)}
    (1+t)^{2g}
    \frac{(1+xt)^{2g}}{(1-x)(1-xt^2)x^{m_{1}}}\cdot
    \frac{(1+yt)^{2g}}{(1-y)(1-yt^2)y^{m_{2}}}
  $$
  with
\begin{equation}\label{eq:bd1m}
\begin{split}
  \bm &= [(2/3)\Delta -(1/3)F(\alpha,\varpi) +1], \\
  \bd &= [(1/3)\Delta - (1/3)G(\alpha,\varpi) + 1].
  \end{split}
\end{equation}
\end{lem}

\begin{proof}
  The identity would be clear from Corollary~\ref{cor:P-t-N-111} and
  Lemma~\ref{lem:index-111} if the latter sum were over $(d_{1},m)$
  satisfying the conditions of Proposition \ref{prop:boundsmm}.
  Now, from these
  equations we see that we need to sum over the closed region in the
  $(m,d_1)$-plane bounded by the lines
  \begin{align*}
    m &= \bm, \\
    d_{1} &= \bd, \\
    2d_{1} -m &= n - s_{1} + 2g-2, \\
    2m - d_1 &= \Delta  + n - s_2 + 2g-2.
  \end{align*}
  Thus, summing over the semi-infinite region defined by
  $d_1 \geq \bd$ and $m \geq \bm$, we introduce in the sum
  extra terms.  But, since the lines given by the third and fourth equations
  represent the conditions $m_1 \geq 0$ and $m_2 \geq
  0$, these extra terms have strictly positive powers of $x$ or $y$
  and hence this does not change the coefficient of $x^0y^0$.
\end{proof}

Using \eqref{eq:def-m_i} we have  that $x^{m_{1}} =
x^{n-s_{1}+2g-2}x^{-2d_{1}}x^{m}$ and $y^{m_{2}} =
y^{\Delta+n-s_{2}+2g-2}y^{d_{1}}y^{-2m}$, hence
\begin{equation}\label{eq:BS2}
  \begin{aligned}
  \Psi(\varpi) &=
  t^{2(4g-4+n+s_{1}+s_{2}-\Delta)}
  (1+t)^{2g}
  \frac{(1+xt)^{2g}}{(1-x)(1-xt^2)}\cdot
  \frac{(1+yt)^{2g}}{(1-y)(1-yt^2)} \\
  &\quad\cdot
  \frac{1}{x^{n-s_{1}+2g-2}y^{\Delta+n-s_{2}+2g-2}}
  \sum_{d_{1}=\bd}^{\infty}\frac{t^{2d_{1}}x^{2d_{1}}}{y^{d_{1}}}
  \sum_{m=\bm}^{\infty}
  \frac{t^{2m}y^{2m}}{x^{m}} \\
  &=
  t^{2(4g-4+n+s_{1}+s_{2}-\Delta)}
  (1+t)^{2g}
  \frac{(1+xt)^{2g}}{(1-x)(1-xt^2)}\cdot
  \frac{(1+yt)^{2g}}{(1-y)(1-yt^2)} \\
  &\quad\cdot
  \frac{(x^2y^{-1}t^2)^{\bd}(x^{-1}y^2t^2)^{\bm}}
  {x^{n-s_{1}+2g-2}y^{\Delta+n-s_{2}+2g-2}
    (1-x^2y^{-1}t^2)(1-x^{-1}y^2t^2)} \\
  &=
  t^{2(s_1+s_2)}x^{s_1}y^{s_2}(x^2y^{-1}t^2)^{\bd}(x^{-1}y^2t^2)^{\bm}
  \\
  &\quad\cdot
  \frac{t^{2(4g-4+n-\Delta)}(1+t)^{2g}(1+xt)^{2g}(1+yt)^{2g}}
  {x^{n+2g-2}y^{\Delta+n+2g-2}(1-x)(1-xt^2)(1-y)(1-yt^2)
    (1-x^2y^{-1}t^2)(1-x^{-1}y^2t^2)},
  \end{aligned}
\end{equation}
where we have separated powers of $x$, $y$ and $t$ which
potentially depend on $\varpi$.

\subsection{The sum over $\varpi$.}

In order to proceed with the calculation we need to sum the
contribution (\ref{eq:BS2}) over all permutations $\varpi
=(\varpi_p)\in (S_3)^n$.  For this we need to understand the
dependence of $s_1$, $s_2$, $\bm$ and $\bd$ on $\varpi$.  Now,
looking at the definitions (\ref{eq:bd1m}), we see that $\bd$ and
$\bm$ also depend on the weights. In order to deal with this
dependence, we shall take advantage of Proposition
\ref{prop:deg-tensor} which allows us to do the computation in the
case where the degree $\Delta$ satisfies that $\Delta\not \equiv
0\pmod 3$.  Hence we can choose the weights so as to facilitate
the computations, as long as we keep the same choice throughout.
Now, if $\Delta \not\equiv 0 \pmod 3$ and we choose the weights
$\a_i(p)$ sufficiently small, then $\bd$ and $\bm$ are independent
of $\varpi$. For future reference, we state here our assumptions.

\begin{assumption}
  \label{ass:small-weights}
  Write $D = p_1 + \cdots + p_n$.  In addition to
  Assumption~\ref{assumption}, we shall from now on assume that
  $\Delta\not\equiv 0 \pmod 3$ and that the
  weights $\a_i(p)$ are chosen to satisfy
   $$
    \a_i(p_j)\ll 1 \qquad\text{for all $i,j$.}
   $$
\end{assumption}

Next we consider the dependence of $s_1$ and $s_2$ on $\varpi$.
We can write
\begin{align*}
  s_1 &= \sum_{p \in D} s_1(p), \\
  s_2 &= \sum_{p \in D} s_2(p),
\end{align*}
where $s_1(p)$ and $s_2(p)$ are defined in the obvious way:
\begin{align*}
  s_1(p) &=
  \begin{cases}
    1 &\text{if $\varpi_p(1) > \varpi_p(2)$,} \nonumber\\
    0 &\text{otherwise,}
  \end{cases}
  \intertext{and}
  s_2(p) &=
  \begin{cases}
    1 &\text{if $\varpi_p(2) > \varpi_p(3)$,} \nonumber\\
    0 &\text{otherwise.}
  \end{cases}
\end{align*}
We give the values of $s_1(p)$, $s_2(p)$ and $s_1(p)+s_2(p)$ as a
function of $\varpi_p$ in Table \ref{tab:si-sigma}, using the
notation $\varpi = (\varpi(1)\,\varpi(2)\,\varpi(3))$ for a
permutation $\varpi \in S_3$.
\begin{table}[htbp]
  \centering
  \caption{$s_1(p)$ and $s_2(p)$ as a function of $\varpi_p$}
\begin{tabular}{|c|c|c|c|c|c|c|}
\hline
  $\varpi_p$ & $(123)$ & $(231)$ & $(312)$ & $(213)$ & $(132)$ & $(321)$ \\
\hline
  $s_1(p)$          &   0   &   0   &   1   &   1   &   0   &   1   \\
\hline
  $s_2(p)$          &   0   &   1   &   0   &   0   &   1   &   1   \\
\hline
  $s_1(p)+s_2(p)$   &   0   &   1   &   1   &   1   &   1   &   2   \\
\hline
\end{tabular}
  \label{tab:si-sigma}
\end{table}

Under Assumption~\ref{ass:small-weights}, $\bm$ and $\bd$ are
independent of $\varpi$: in fact we have from the definitions
(\ref{eq:bd1m}) of $\bm$ and $\bd$ that
\begin{equation}\label{eq:bdm-sigma-fix}
\begin{split}
  \bm &= \left[\tfrac{2\Delta}{3}\right]+1, \\
  \bd &= \left[\tfrac{\Delta}{3}\right]+1.
\end{split}
\end{equation}
Therefore to do the sum $\sum_{\varpi}\Psi(\varpi)$, we only need
to do
\begin{math}
  \sum_{\varpi \in (S_3)^n}t^{2(s_1+s_2)}x^{s_1}y^{s_2}
\end{math}
in (\ref{eq:BS2}).  For this we use Table \ref{tab:si-sigma} and
obtain:
\begin{equation}
  \label{eq:sum-sigma-1}
  \begin{aligned}
    \sum_{\varpi \in (S_3)^n}t^{2(s_1+s_2)}x^{s_1}y^{s_2}
    &=
    \prod_{p \in D}\sum_{\varpi_p\in S_3}
      t^{2(s_1(p)+s_2(p))}x^{s_1(p)}y^{s_2(p)} \\
    &=
    \prod_{p \in D}(1+2t^2x+2t^2y+t^4xy) \\
    &= (1+2t^2x+2t^2y+t^4xy)^n.
  \end{aligned}
\end{equation}
Combining (\ref{eq:sum-sigma-1}) with (\ref{eq:BS2}) we finally
obtain:
\begin{equation}\label{eq:sum-sigma-Delta-not-0}
\begin{split}
 \sum_{\varpi}\Psi(\varpi)
   &= (1+2t^2x+2t^2y+t^4xy)^n \, \cdot\\
  &\cdot \, \frac{t^{2(4g-4+n-\Delta)}(1+t)^{2g}
    (1+xt)^{2g}(1+yt)^{2g}
  (x^2y^{-1}t^2)^{\bd}(x^{-1}y^2t^2)^{\bm}}
  {x^{n+2g-2}y^{\Delta+n+2g-2}(1-x)(1-xt^2)(1-y)(1-yt^2)
    (1-x^2y^{-1}t^2)(1-x^{-1}y^2t^2)}.
\end{split}
\end{equation}
Note that this expression has arbitrarily large positive and
negative powers of $x$ and $y$.  Therefore it is not very suitable
for extracting the coefficient to $x^0y^0$.  However, to
facilitate this task we can make the following change of variable:
 $$
  x=u^2v, \qquad  y =uv^2.
 $$
Then we have
\begin{displaymath}
  x^2y^{-1} = u^3, \qquad x^{-1}y^2 = v^3 \qquad\text{and}\qquad xy = u^3v^3.
\end{displaymath}
Substituting in (\ref{eq:sum-sigma-Delta-not-0}), and
using (\ref{eq:bdm-sigma-fix}), we finally obtain the formula for
the contribution to the Poincar\'e polynomial from critical
submanifolds of type $(1,1,1)$:

\begin{prop}\label{prop:P_t111_nonzero}
  Under Assumption \ref{ass:small-weights}, let $\Delta_0 \in
  \{1,2\}$ be the remainder modulo $3$ of $\Delta$. Then
\begin{multline*}
 P_t(\Delta,(1,1,1)) =
\Coeff_{u^0v^0}  \bigg( (1+2u^2vt^2+2uv^2t^2+u^3v^3t^4)^n  \,
 \cdot \\ \cdot \,
  \frac{t^{2(4g-3+n)}(1+t)^{2g}
  (1+u^2vt)^{2g}(1+uv^2t)^{2g}}
  {u^{3n+6g-9+\Delta_0}v^{3n+6g-6-\Delta_0}
    (1-u^2v)(1-uv^2)(1-u^2vt^2)(1-uv^2t^2)
    (1-v^3t^2)(1-u^3t^2)} \bigg).
\end{multline*}
\qed
\end{prop}

\section{Critical submanifolds of type $(1,2)$} \label{sec:(1,2)}

\subsection{Description of the critical submanifolds}
\label{sec:(1,2)-description}
In this section, we consider the critical points of the
Bott--Morse function $f$ represented by Higgs bundles $(E,\Phi)$
which are of the form $E=E_0\oplus E_1$ where $E_0=L$ is a
parabolic line bundle, $E_1$ is a rank $2$ parabolic bundle and
$\Phi:L\to E_1\otimes K(D)$ is a strongly parabolic homomorphism. This
defines a parabolic triple $(E_1\ox K, L)$ of type $(1,2)$. By
Proposition \ref{prop:august}, the triple is $\s$-stable exactly
for the value $\s=2g-2$.

In order to do the computations, let us introduce some notation.
Recall that we keep our assumption of generic weights. The (fixed)
weights of $E$ at each $p\in D$ are $0<\a_1(p)<\a_2(p)<\a_3(p)<1$.
Each possible choice of distribution of these weights is given by
a function $\varpi=\{ \varpi_p\}_{p\in D}$ that assigns to each
$p\in D$ a number $\varpi(p) \in \{1,2,3\}$ such that the weight
of $L$ is $\a(p)=\a_{\varpi(p)}(p)$. The weights of $E_1$ are
$\b_1(p)<\b_2(p)$, so that $\{\b_1(p),\b_2(p),\a(p)\}=
\{\a_1(p),\a_2(p),\a_3(p)\}$. In the decomposition $E=L \oplus
E_1$, we define
  \begin{eqnarray*}
  d_1 &=& \deg(E_1\ox K)=\deg(E_1) + 4g-4\, , \\
  d_2 &=& \deg(L)\, ,
  \end{eqnarray*}
so that $\Delta=\deg (E)= d_1+d_2 +4-4g$.

Denote by $\cN_{(1,2)}(d_1,\varpi)$ the critical submanifold of
parabolic Higgs bundles of type $(1,2)$ with topological
invariants given by $(d_1, d_2=\Delta-d_1+4g-4)$ and weights
determined by $\varpi$. The contribution of all critical
submanifolds of type $(1,2)$ is given as
   $$
   P_t(\Delta, (1,2)):=\sum_{d_1,\varpi} t^{\lambda_{(d_1,\varpi)}}
   P_t(\cN_{(1,2)}(d_1,\varpi)),
   $$
where $\lambda_{(d_1,\varpi)}$ is the index of
$\cN_{(1,2)}(d_1,\varpi)$.

\begin{lem} \label{lem:ind(1,2)}
 The index of the critical submanifold $\cN_{(1,2)}(d_1,\varpi)$ is
 $$
 \lambda_{(d_1,\varpi)} = 12g-12 + 4 n -2d_1 + 4d_2-2s_0,
 $$
 where $s_0=\#\{\b_i(p)\, | \, \b_i(p)>\a(p) \}$.
\end{lem}

\begin{proof}
{}From Proposition \ref{prop:phb-moduli-dim}, $f_p=\frac12
(r^2-\sum_i m_i(p)^2)=3$, since the multiplicities are all $1$.
For $p\in D$, the space $P_p(L,L)$ consists of endomorphisms of
$L_p$, so $\dim P_p(L,L)=1$ and $P_p(E_1\ox K, E_1\ox K)$ consists
of parabolic endomorphisms of $(E_1\ox K)_p$, so it has dimension
$3$. The dimension of the space of strongly parabolic
homomorphisms from $L_p$ to $(E_1\ox K)_p$ is given by
  $$
  \dim N_p(L,E_1\ox K)=\left\{ \begin{array}{ll} 2 \qquad &
  \mathrm{if}\,\, \a(p)<\b_1(p), \\
  1 \qquad & \mathrm{if}\,\, \b_1(p)<\a(p)<\b_2(p),\\
  0 \qquad & \mathrm{if}\,\, \b_2(p)<\a(p). \end{array}\right.
  $$
Therefore
 $$
 s_0=\sum_p \dim N_p(L, E_1\ox K) =\#\{\b_i(p)\, |\, \b_i(p)>\a(p) \}.
 $$
Substituting this in the formula for the Morse index in
Proposition~\ref{prop-index}, we have
  \begin{eqnarray*}
  \lambda_{(d,\varpi)}
  &=& 9 (2g-2)+ 6n + 2\big( 5(1-g-n) + 2(1-g)
  - (d_1-4g+4)+ 2d_2 + 4n - s_0 \big) \\
   &=& 12g-12 + 4n -2d_1 + 4d_2-2s_0.
  \end{eqnarray*}
\end{proof}

By Proposition \ref{prop:august}, $\cN_{(1,2)}(d_1,\varpi)$ is
isomorphic to the moduli space of $\s$-stable triples of the
appropriate type with $\s=2g-2$. By the genericity of weights,
such $\sigma$ is not a critical value. Its Poincar{\'e} polynomial is
given by Theorem \ref{thm:Pt(Ns)}. So, for each
$\e=\{\e(p)\}_{p\in D}$, let $s_1,s_2,s_3$ be defined by
(\ref{eqn:6.2.5}), $\varsigma(p)=3-\e(p)$ and
 \begin{eqnarray*}
 \bdm &=& \left[\frac13 \left(d_1+d_2 +\sum \big(\a(p) +\b_{\varsigma(p)}(p)
 -2\b_{\e(p)}(p)\big)
 +2g-2\right) \right]+1 \\ &=&
   \left[\frac13\left(\Delta +\sum \big(\a(p) +\b_{\varsigma(p)}(p)
 -2\b_{\e(p)}(p) \big) \right)\right]  +2g-1.
 \end{eqnarray*}
Then Theorem \ref{thm:Pt(Ns)}, Lemma \ref{lem:ind(1,2)} and the
fact that  $s_0=s_1+s_2$ yield that
 \begin{multline} \label{eqn:Pt(N(1,2))}
 P_t(\Delta,(1,2))= \\
 \sum_{d_1,\varpi} t^{12g-12 + 4 n -2d_1 + 4d_2-2s_1-2s_2}
 \sum_{\e} \Coeff_{x^0}\bigg(
 \frac{(1+t)^{4g}(1+xt)^{2g}
 t^{2d_1-2d_2+2s_2+2s_3-2\bdm}x^{\bdm}}{(1-t^2)(1-x)(1-xt^2)(1-t^{-2}x)
 x^{d_1-d_2+s_1}} \\
 -
\frac{(1+t)^{4g}(1+xt)^{2g}
 t^{-2d_1+2g-2 +2n-2s_3+4\bdm}x^{\bdm} }{(1-t^2)(1-x)(1-xt^2)(1-t^4x)
 x^{d_1-d_2+s_1}} \bigg) \, .
 \end{multline}

\subsection{The sum for fixed $(\varpi,\e)$.}

We shall compute the contribution (\ref{eqn:Pt(N(1,2))}) to the
Poincar{\'e} polynomial of $\mathcal{M}$ from submanifolds of type
$(1,2)$. We start by performing the sum over all possibilities of
$d_1$ for each choice of $(\varpi,\e)$. The condition that the
moduli space $\cN_{(1,2)}(d_1,\varpi)$ be non-empty is
$2g-2>\pmu_1-\pmu_2$ (see Remark \ref{rem:Pt(Ns)}). This means
that
 $$
 2g-2> d_1/2-d_2 + \sum \big(\b_1(p) +\b_2(p)-2\a(p)\big)/2.
 $$
Using that $\Delta=d_1+d_2+4-4g$, this is translated into
 $$
 d_2-d_1 >4-4g - \frac\Delta3 +\frac23 \sum \big(\b_1(p) +\b_2(p)-2\a(p)\big).
 $$
But  $d_2-d_1\equiv \Delta \pmod 2$. So
$d_2-d_1=\bdd+2k$, $k\geq 0$, and
\begin{equation}\nonumber
 \bdd=4 -4g + 2\left[  \frac12 \left( \left[ -\frac{\Delta}{3} +
 \frac23 \sum \big(\b_1(p) +\b_2(p)-2\a(p)\big) \right] + \Delta \right)\right]-\Delta +2.
 \end{equation}
This gives the range for the summation in (\ref{eqn:Pt(N(1,2))})
for $d_1$ for fixed $(\varpi,\e)$. Looking at
(\ref{eqn:Pt(N(1,2))}), one sees that we need to compute
 \begin{eqnarray*}
 \sum t^{2d_2}  x^{d_2-d_1} &=&t^{\Delta+4g-4} \sum t^{d_2-d_1} x^{d_2-d_1}=
 \frac{t^{\bdd+\Delta+4g-4}x^{\bdd}}{1-t^2x^2}\, ,\\
 \sum t^{-4d_1+4d_2}  x^{d_2-d_1} &=&\frac{t^{4\bdd}x^{\bdd}}{1-t^8x^2} \, .
 \end{eqnarray*}

Substituting into (\ref{eqn:Pt(N(1,2))}), we get that
 \begin{multline}
 \label{eqn:contrib_2}
 P_t(\Delta,(1,2))=\sum_{\varpi,\e}
 \Coeff_{x^0}\bigg(  \frac{(1+t)^{4g}(1+xt)^{2g}
 t^{16g-16+4n-2s_1+2s_3-2\bdm+\bdd+\Delta}
 x^{\bdm+\bdd-s_1}}{(1-t^2)(1-x)(1-xt^2)(1-t^{-2}x)(1-t^2x^2)} \\
  -\frac{(1+t)^{4g}(1+xt)^{2g}
 t^{14g-14 +6n-2s_1-2s_2-2s_3+ 4\bdm + 4\bdd}
 x^{\bdm+\bdd-s_1} }{(1-t^2)(1-x)(1-xt^2)(1-t^4x) (1-t^8x^2)}
 \bigg)\, .
 \end{multline}

\subsection{The sum over $\varpi$ and $\e$.}

Now we need to perform the sum in (\ref{eqn:contrib_2}) for all
choices of $(\varpi,\e)$. For this we arrange the degree and the
weights to satisfy Assumption \ref{ass:small-weights}. Write
$\Delta=3k+\Delta_0$, $\Delta_0\in\{1,2\}$. Since $\alpha_i(p)$
are sufficiently small, we have that
 \begin{equation} \label{eqn:bddm-au}
 \begin{split}
  \bdm &= \left[\frac\Delta3 \right]  +2g-1 = k+2g-1\, , \\
  \bdd &= 4 -4g + 2\left[  \frac12 \left( \left[ -\frac{\Delta}{3}
  \right] + \Delta \right)\right]-\Delta+2 = 6-4g-k-\Delta_0\,
  \end{split}
 \end{equation}
are independent of $(\varpi,\e)$. Therefore to do the sum
\eqref{eqn:contrib_2}, we only need to do
 $$
  \sum_{\varpi,\e} t^{-2s_1+2s_3}x^{-s_1} \qquad
  \hbox{ and }\qquad \sum_{\varpi,\e} t^{-2s_1-2s_2-2s_3}x^{-s_1}.
 $$
We have to write down the dependence of $s_1,s_2,s_3$ on
$(\varpi,\e)$. Note that we can write
\begin{align*}
  s_1 &= \sum_{p \in D} s_1(p), \\
  s_2 &= \sum_{p \in D} s_2(p), \\
  s_3 &= \sum_{p \in D} s_3(p),
\end{align*}
where $s_1(p)$, $s_2(p)$ and $s_3(p)$ are defined in the obvious
way:
\begin{align*}
  s_1(p) &=
  \begin{cases}
    1 &\text{if $\alpha(p) < \beta_{\varsigma(p)}(p)$,} \nonumber\\
    0 &\text{otherwise,}
  \end{cases}
  \intertext{ }
  s_2(p) &=
  \begin{cases}
    1 &\text{if $\alpha(p) < \beta_{\e(p)}(p)$,} \nonumber \\
    0 &\text{otherwise,}
  \end{cases}
  \intertext{and}
  s_3(p) &=
  \begin{cases}
    1 &\text{if $\beta_{\e(p)}(p)<\beta_{\varsigma(p)}(p)$,}\nonumber \\
    0 &\text{otherwise.}
  \end{cases}
\end{align*}
We give the values of $s_1(p)$, $s_2(p)$ and $s_3(p)$ as a
function of $(\varpi(p),\e(p)) \in \{1,2,3\}\x \{1,2\}$ in Table
\ref{tab:si-varpi},
\begin{table}[htbp]
  \centering
  \caption{$s_1(p)$, $s_2(p)$ and $s_3(p)$ as a function of $(\varpi(p),\e(p))$}
\begin{tabular}{|c|c|c|c|c|c|c|}
\hline
 $(\varpi(p),\e(p))$ & $(1,1)$ & $(1,2)$ & $(2,1)$ & $(2,2)$ & $(3,1)$ & $(3,2)$ \\ \hline
  $s_1(p)$           &   1   &   1   &   1   &   0   &   0   &   0   \\
\hline
  $s_2(p)$           &   1   &   1   &   0   &   1   &   0   &   0   \\
\hline
  $s_3(p)$           &   1   &   0   &   1   &   0   &   1   &   0   \\
\hline
\end{tabular}
  \label{tab:si-varpi}
\end{table}

We obtain
 \begin{align*}
   \sum_{\varpi,\e} t^{-2s_1+2s_3}x^{-s_1}  &= \prod_{p\in D}
   \sum_{\varpi(p),\e(p)}
   t^{-2s_1(p)+2s_3(p)}x^{-s_1(p)} \\ &=
   \prod_{p\in D} (t^{-2}x^{-1}+2+2x^{-1}+t^2) \\
   &=(t^{-2}x^{-1}+2+2x^{-1}+t^2)^n\\
   &= t^{-2n}x^{-n}(1+2t^2+2t^2x+t^4x)^n\, , \\
\intertext{and}
   \sum_{\varpi,\e} t^{-2s_1-2s_2-2s_3}x^{-s_1} &= \prod_{p\in D}
   \sum_{\varpi(p),\e(p)}
   t^{-2s_1(p)-2s_2(p)-2s_3(p)}x^{-s_1(p)}  \\  &= (2 t^{-4}x^{-1}
   + x^{-1}t^{-6} +1 +2 t^{-2} )^n \\
   &= t^{-6n}x^{-n}(1 + 2t^2  +2t^{4}x + t^{6}x )^n\, .
  \end{align*}
Combining this with \eqref{eqn:contrib_2}
we get that
$P_t(\Delta,(1,2))$ equals
 $$
 \Coeff_{x^0}\bigg( \frac{(1+t)^{4g}(1+xt)^{2g}
 t^{16g-16+2n -2\bdm+\bdd+\Delta}
 x^{\bdm+\bdd-n} (1+2t^2+2t^{2}x+t^4x)^n
 }{(1-t^2)(1-x)(1-xt^2)(1-t^{-2}x)(1-t^2x^2)}
 $$
 $$
 \qquad \qquad - \frac{(1+t)^{4g}(1+xt)^{2g}
 t^{14g-14 + 4\bdm + 4\bdd}
 x^{\bdm+\bdd-n} (1+ 2 t^{2} +2 t^4x +t^6x )^n}
{(1-t^2)(1-x)(1-xt^2)(1-t^4x) (1-t^8x^2)}
 \bigg)\, .
 $$
Now, using that
$\bdm+\bdd=5-2g-\Delta_0$ and $-3\bdm+\Delta= 3-6g+\Delta_0$,
which follow from \eqref{eqn:bddm-au}, we have the following.

\begin{prop} \label{prop:P_t(1,2)_nonzero}
Under Assumption \ref{ass:small-weights}, let $\Delta_0\in
\{1,2\}$ be the remainder modulo $3$ of $\Delta$. Then
 \begin{multline*}
   P_t(\Delta,(1,2)) =\Coeff_{x^0}\bigg( \frac{(1+t)^{4g}(1+xt)^{2g}
 t^{8g-8+2n}x^{5-2g-\Delta_0-n}
 (1+2t^2+2t^{2}x+t^4x)^n}{(1-t^2)(1-x)(1-xt^2)(1-t^{-2}x)(1-t^2x^2)}  \\
 - \frac{(1+t)^{4g}(1+xt)^{2g}
 t^{6g +6 -4\Delta_0} x^{5-2g-\Delta_0-n} (1+ 2 t^{2} +2 t^4x +t^6x )^n}{(1-t^2)
 (1-x)(1-xt^2)(1-t^4x)(1-t^8x^2)} \bigg) \, .
  \end{multline*}
\qed
\end{prop}

\section{Critical submanifolds of type $(2,1)$} \label{sec:(2,1)}

\subsection{Description of the critical submanifolds}
In this section, we consider the critical points of the
Bott--Morse function $f$ represented by Higgs bundles $(E,\Phi)$
which are of the form $E=E_0\oplus E_1$ where $E_1=L$ is a
parabolic line bundle, $E_0$ is a rank $2$ parabolic bundle and
$\Phi:E_0\to L\otimes K(D)$ is a strongly parabolic homomorphism. This
defines a parabolic triple $(L\ox K,E_0)$ of type $(2,1)$. By
Proposition \ref{prop:august}, the triple is $\s$-stable exactly
for the value $\s=2g-2$.

As in Section \ref{sec:(1,2)}, the (fixed) weights of $E$ at each
$p\in D$ are $0<\a_1(p)<\a_2(p)<\a_3(p)<1$. Each possible choice
of distribution of these weights is given by some
$\varpi=\{\varpi(p)\}_{p\in D}$ where $\varpi(p) \in \{1,2,3\}$,
$p\in D$ such that the weight of $L$ is $\a(p)=\a_{\varpi(p)}(p)$.
The weights of $E_0$ are $\b_1(p)<\b_2(p)$, so that
$\{\b_1(p),\b_2(p),\a(p)\}= \{\a_1(p),\a_2(p),\a_3(p)\}$. In the
decomposition $E=E_0 \oplus L$, we define
  \begin{eqnarray*}
  d_1 &=& \deg(L\ox K)=\deg(L) + 2g-2, \\
  d_2 &=& \deg(E_0),
  \end{eqnarray*}
so that $\Delta=\deg (E)= d_1+d_2 +2-2g$.

Denote by $\cN_{(2,1)}(d_1,\varpi)$ the critical submanifold of
parabolic Higgs bundles of type $(2,1)$ with topological
invariants given by $(d_1, d_2=\Delta-d_1+2g-2)$ and where the
weights are determined by $\varpi$.

\begin{lem}\label{lem:9.1}
 $\cN_{(2,1)}(d_1,\varpi)$ is isomorphic to the moduli space of $\s$-stable
 parabolic triples of type $(1,2)$ with degrees $d_2'=-n-d_1, d_1'=-2n-d_2$
 and weights
 $1-\alpha(p)$ for the line bundle and $1-\beta_2(p)<1-\beta_1(p)$
 for the rank $2$-bundle, for $\s=2g-2$.
\end{lem}

\begin{proof}
 The result follows by dualizing and applying Proposition
 \ref{prop:duality}. And by the definition of
 dual of a parabolic bundle.
\end{proof}

Note that by the genericity of weights, the value $\sigma=2g-2$ is
not a critical value for the moduli space of parabolic triples.
Now let $\e=\{\e(p)\}_{p\in D}$, $\e(p)\in \{1,2\}$, and
$\varsigma(p)=3-\e(p)$. We introduce the following sets:
 \begin{eqnarray*}
 S_1 &=& \{ p \in D \, | \, 1-\alpha(p)< 1-\beta_{\varsigma(p)}(p) \}=
 \{ p \in D \, | \, \alpha(p)> \beta_{\varsigma(p)}(p) \},\\
 S_2 &=& \{ p \in D \, | \, 1-\alpha(p)< 1-\beta_{\e(p)}(p) \}=
 \{ p \in D \, | \, \alpha(p)> \beta_{\e(p)}(p) \}, \\
 S_3 &=& \{ p \in D \, | \,  1-\beta_{\e(p)}(p) <1-\beta_{\varsigma(p)}(p)
 \}=\{ p \in D \, | \,  \beta_{\e(p)}(p) >\beta_{\varsigma(p)}(p)
 \}.
 \end{eqnarray*}
and denote
\begin{displaymath}
  s_1=\# S_1, \qquad s_2=\# S_2 \qquad\text{and} \qquad s_3=\# S_3.
\end{displaymath}
Applying
Theorem \ref{thm:Pt(Ns)}, we have
 $$
 \begin{aligned}
   \bdm &= \bigg[\frac13 \Big(-n-d_1-2n-d_2 \\
       &\qquad +\sum \left(1-\a(p) +
 1-\b_{\varsigma(p)}(p) -2+2\b_{\e(p)}(p)\right)
 +2g-2\Big) \bigg]+1 \\ &=
  -n +  \left[\frac13\left(-\Delta -\sum \big(\a(p) +\b_{\varsigma(p)}(p)
 -2\b_{\e(p)}(p) \big) \right)\right]  +1.
 \end{aligned}
 $$
Then
 \begin{multline} \label{eqn:Pt(N(2,1))}
 P_t(\cN_{(2,1)}(d_1,\varpi)) = \sum_{\e} \Coeff_{x^0} \bigg( \frac{(1+t)^{4g}(1+xt)^{2g}
 t^{2d_1-2d_2-2n+2s_2+2s_3-2\bdm}x^{\bdm}}{(1-t^2)(1-x)(1-xt^2)(1-t^{-2}x)
 x^{d_1-d_2-n+s_1}}  \\ - \frac{(1+t)^{4g}(1+xt)^{2g}
 t^{2d_2+2g-2 +6n-2s_3+4\bdm}x^{\bdm} }{(1-t^2)(1-x)(1-xt^2)(1-t^4x)
 x^{d_1-d_2-n+s_1}} \bigg) \, .
 \end{multline}

Similarly to Lemma \ref{lem:ind(1,2)}, we can prove

\begin{lem} \label{lem:ind(2,1)}
 The index of the critical submanifold $\cN_{(2,1)}(d_1,\varpi)$ is
 $$
 \lambda_{(d_1,\varpi)} = 12g-12 + 4 n -4d_1 + 2d_2-2s_0,
 $$
 where $s_0=\#\{\b_i(p) \ | \  \b_i(p)<\a(p) \} =s_1+s_2$.
\end{lem}

\begin{proof}
  {}From Proposition \ref{prop:phb-moduli-dim}, $f_p=\frac12
(r^2-\sum m_j(p)^2)=3$, since the multiplicities are all $1$. For
$p\in D$, the space $P_p(L\ox K,L\ox K)$ consists of endomorphisms
of $L_p$, so $\dim P_p(L\ox K,L\ox K)=1$ and $P_p(E_0, E_0)$
consists of parabolic endomorphisms of $(E_0)_p$, so it has
dimension $3$. On the other hand, $\dim N_p(E_0, L\ox K)$ is the
dimension of the space of strongly parabolic homomorphisms from
$(E_0)_p$ to $(L\ox K)_p$, so
  $$
  \dim N_p(E_0,L\ox K)=\left\{ \begin{array}{ll} 0 \qquad &
  \mathrm{if}\,\, 1-\a(p)>1-\b_1(p), \\
  1 \qquad & \mathrm{if}\,\, 1-\b_1(p)>1-\a(p)>1-\b_2(p),\\
  2 \qquad & \mathrm{if}\,\, 1-\b_2(p)>1-\a(p). \end{array}\right.
  $$
Therefore
 $$
 s_0=\sum_p \dim N_p(E_0, L\ox K) =\#\{\b_i(p) | \b_i(p)<\a(p) \}.
 $$
Substituting this in the formula for the Morse index
of Proposition~\ref{prop-index}, we have
  \begin{eqnarray*}
  \lambda_{(d_1,\varpi)} &=&
  9 (2g-2)+ 6n + 2\left( 5(1-g-n) +4n \right) \\
  && + 2\left( 2(1-g)
  - 2(d_1-2g+2)+ d_2 - s_0 \right) \\
   &=& 12g-12 + 4 n -4d_1 + 2d_2-2s_0.
  \end{eqnarray*}
\end{proof}

Therefore the contribution of all critical submanifolds of type
$(2,1)$ is given as
 \begin{equation}
 \label{eqn:contrib(2,1)}
   \begin{split}
   P_t(\Delta, (2,1)) & : = \sum_{d_1,\varpi} t^{\lambda_{(d_1,\varpi)}}
   P_t(\cN_{(2,1)}(d_1,\varpi)) \\
  &= \sum_{d_1,\varpi}
 \sum_{\e} \Coeff_{x^0}\bigg(  \frac{(1+t)^{4g}(1+xt)^{2g}
  t^{12g-12 + 2 n -2d_1 -2s_1+2s_3-2\bdm}
 x^{\bdm}}{(1-t^2)(1-x)(1-xt^2)(1-t^{-2}x)
 x^{d_1-d_2-n+s_1}} \\
  & \qquad -\frac{(1+t)^{4g}(1+xt)^{2g}
 t^{14g-14+10n-4d_1+4d_2  -2s_1-2s_2-2s_3+4\bdm}x^{\bdm} }{(1-t^2)(1-x)(1-xt^2)(1-t^4x)
 x^{d_1-d_2-n+s_1}} \bigg)\, ,
 \end{split}
 \end{equation}
obtained by using Lemma \ref{lem:ind(2,1)} and
(\ref{eqn:Pt(N(2,1))}).

\subsection{The sum for fixed $(\varpi,\e)$.}

Now we shall compute the contribution (\ref{eqn:contrib(2,1)}) to
the Poincar{\'e} polynomial of $\mathcal{M}$ from submanifolds of type
$(2,1)$. As before, we do first the sum over all possibilities of
$d_1$ for each choice of $(\varpi,\e)$. The condition that the
moduli space $\cN_{(2,1)}(d_1,\varpi)$ be non-empty is
$2g-2>\pmu_1-\pmu_2$. This means that
 $$
 2g-2> d_1-d_2/2 + \sum \big(2\a(p)-\b_1(p) -\b_2(p)\big)/2,
 $$
Using that $\Delta=d_1+d_2+2-2g$, this is translated into
 $$
 d_2-d_1 >2-2g + \frac\Delta3 + \frac43 \sum \big( 2\a(p)-
 \b_1(p) -\b_2(p)\big).
 $$
But  $d_2-d_1\equiv \Delta \pmod 2$. So
$d_2-d_1=\bdd+2k$, $k\geq 0$, and
 $$
 \bdd=2-2g + 2\left[  \frac12 \left( \left[ \frac{\Delta}{3} +
 \frac43 \sum \big(2\a(p)-\b_1(p) -\b_2(p)\big) \right]-\Delta
 \right)\right]+\Delta.
 $$

In (\ref{eqn:contrib(2,1)}) we need to compute the terms
 \begin{eqnarray*}
  \sum t^{-2d_1}x^{d_2-d_1} &=&t^{-\Delta-2g+2} \sum t^{d_2-d_1} x^{d_2-d_1}=
 \frac{t^{\bdd-\Delta-2g+2}x^{\bdd}}{1-t^{2}x^{2}} \, ,\\
 \sum t^{-4d_1+4d_2}x^{d_2-d_1} &=&\frac{t^{4\bdd}x^{\bdd}}{1-t^8x^2}\, .
 \end{eqnarray*}
Plugging this into \eqref{eqn:contrib(2,1)} we get
  \begin{multline} \label{eqn:au-au}
 P_t(\Delta, (2,1))= \sum_{\varpi,\e}
 \Coeff_{x^0}\bigg( \frac{(1+t)^{4g}(1+xt)^{2g}
 t^{10g-10 + 2 n-2s_1+2s_3-2\bdm+\bdd-\Delta}
 x^{\bdm+\bdd+n-s_1}}{(1-t^2)(1-x)(1-xt^2)(1-t^{-2}x)(1-t^{2}x^{2})}
 \\ -\frac{(1+t)^{4g}(1+xt)^{2g}
 t^{14g-14 +10n-2s_1-2s_2-2s_3+ 4\bdm + 4\bdd}
 x^{\bdm+\bdd+n-s_1} }{(1-t^2)(1-x)(1-xt^2)(1-t^4x) (1-t^8x^{2})}
 \bigg)\, .
 \end{multline}

\subsection{The sum over $\varpi$ and $\e$.}

To perform the sum in (\ref{eqn:au-au}) for all choices of
$(\varpi,\e)$, we arrange the degree and the weights to satisfy
Assumption \ref{ass:small-weights}. Write $\Delta=3k+\Delta_0$,
$\Delta_0\in\{1,2\}$. Since $\alpha_i(p)$ are sufficiently small,
we have that
 \begin{eqnarray*}
  \bdm &=& -n+\left[-\frac\Delta3 \right]  +1 =-n-k\, , \\
  \bdd &=& 2 -2g + 2\left[  \frac12 \left( \left[ \frac{\Delta}{3}
  \right] - \Delta \right)\right]+\Delta+2=2-2g+k+\Delta_0\, ,
 \end{eqnarray*}
are independent of $(\varpi,\e)$. Therefore to do the sum in
\eqref{eqn:au-au}, we only need to compute
 $$
 \sum_{\varpi,\e} t^{-2s_1+2s_3}x^{-s_1} \qquad \hbox{and}\qquad \sum_{\varpi,\e}
 t^{-2s_1-2s_2-2s_3}x^{-s_1}\, .
 $$
As before,
\begin{align*}
  s_1 &= \sum_{p \in D} s_1(p)\ , \\
  s_2 &= \sum_{p \in D} s_2(p)\ , \\
  s_3 &= \sum_{p \in D} s_3(p)\ ,
\end{align*}
where $s_1(p)$, $s_2(p)$ and $s_3(p)$ are defined in the obvious
way:
\begin{align*}
  s_1(p) &=
  \begin{cases}
    1 &\text{if $\alpha(p) > \beta_{\varsigma(p)}(p)$,}  \\
    0 &\text{otherwise,}
  \end{cases}
  \intertext{ }
  s_2(p) &=
  \begin{cases}
    1 &\text{if $\alpha(p) > \beta_{\e(p)}(p)$,}\\
    0 &\text{otherwise,}
  \end{cases}
  \intertext{and}
  s_3(p) &=
  \begin{cases}
    1 &\text{if $\beta_{\e(p)}(p)>\beta_{\varsigma(p)}(p)$,}  \\
    0 &\text{otherwise.}
  \end{cases}
\end{align*}
The values of $s_1(p)$, $s_2(p)$ and $s_3(p)$ as a function of
$(\varpi(p),\e(p))$ are in Table \ref{tab:si-varpi2},
\begin{table}[htbp]
  \centering
  \caption{$s_1(p)$, $s_2(p)$ and $s_3(p)$ as a function of $(\varpi(p),\e(p))$}
\begin{tabular}{|c|c|c|c|c|c|c|}
\hline
 $(\varpi(p),\e(p))$ & $(1,1)$ & $(1,2)$ & $(2,1)$ & $(2,2)$ & $(3,1)$ & $(3,2)$ \\ \hline
  $s_1(p)$           &   0   &   0   &   0   &   1   &   1   &   1   \\
\hline
  $s_2(p)$           &   0   &   0   &   1   &   0   &   1   &   1   \\
\hline
  $s_3(p)$           &   0   &   1   &   0   &   1   &   0   &   1   \\
\hline
\end{tabular}
  \label{tab:si-varpi2}
\end{table}

We obtain
\begin{align*}
   \sum_{\varpi,\e} t^{-2s_1+2s_3}x^{-s_1}  &= \prod_{p\in D}
   \sum_{\varpi(p),\e(p)}
   t^{-2s_1(p)+2s_3(p)}x^{-s_1(p)} \\ &= ( t^{-2}x^{-1}+2+2x^{-1}+t^2)^n \\
   &= t^{-2n}x^{-n}( 1+2t^{2}+2t^2x+t^4x)^n \, ,\\
\intertext{and}
   \sum_{\varpi,\e} t^{-2s_1-2s_2-2s_3}x^{-s_1} &=
   (2 t^{-4}x^{-1} + t^{-6}x^{-1} +1 +2 t^{-2} )^n \\
   &= t^{-6n}x^{-n}(1+2t^{2} + 2t^4x +t^6 x)^n\, .
\end{align*}
Plugging this into \eqref{eqn:au-au} we get that
$P_t(\Delta,(2,1))$ equals
 $$
 \Coeff_{x^0}\bigg(  \frac{(1+t)^{4g}(1+xt)^{2g}
 t^{10g-10 -2\bdm+\bdd-\Delta}
 x^{\bdm+\bdd} ( 1+2t^2+2t^{2}x+t^4x)^n
 }{(1-t^2)(1-x)(1-xt^2)(1-t^{-2}x)(1-t^2x^2)}
 $$
 $$
 \qquad \qquad -\frac{(1+t)^{4g}(1+xt)^{2g}
 t^{14g-14 +4n+ 4\bdm + 4\bdd}
 x^{\bdm+\bdd} (1 + 2t^{2} +2t^4x +t^6x )^n }
 {(1-t^2)(1-x)(1-xt^2)(1-t^4x) (1-t^8x^2)}
 \bigg)\, .
 $$
Now, using that
$\bdm+\bdd=2-2g-n+\Delta_0$ and $-2\bdm+\bdd-\Delta= 2-2g+2n$, we
obtain the following.

\begin{prop} \label{prop:P_t(2,1)_nonzero}
 Under Assumption \ref{ass:small-weights}, let $\Delta_0\in
 \{1,2\}$ be the remainder modulo $3$ of $\Delta$. Then
  \begin{multline*}
  P_t(\Delta,(2,1)) =\Coeff_{x^0} \bigg( \frac{(1+t)^{4g}(1+xt)^{2g}
 t^{8g-8+2n}x^{2-2g+\Delta_0-n}( 1+2t^2+2t^{2}x+t^4x)^n}{(1-t^2)
 (1-x)(1-xt^2)(1-t^{-2}x)(1-t^2x^2) }  \\
 -  \frac{(1+t)^{4g}(1+xt)^{2g}
 t^{6g-6 +4\Delta_0 } x^{2-2g+\Delta_0- n}(1 + 2t^{2} +2t^4x +t^6x )^n}{(1-t^2)
 (1-x)(1-xt^2)(1-t^4x)(1-t^8x^2)}
 \bigg) \, .
  \end{multline*}
\qed
\end{prop}

\section{Betti numbers of the  moduli space
of rank three parabolic bundles}\label{par-rank3}

The Betti numbers of the moduli space of parabolic vector
bundles were computed  by Nitsure \cite{Ni} and Holla \cite{Ho}.
Here we work out Holla's formula for the special case when the
rank is $3$ and all flags at the parabolic points are full.  We
also continue to work with the choice of weights made in
Assumptions \ref{assumption} and \ref{ass:small-weights}.

\subsection{Notation}
\label{sec:hilla-notation}

Given a parabolic bundle $E$, the corresponding
\emph{quasi-parabolic
  data}, $R$, gives the multiplicity of each step of the flag at the
parabolic points:
\begin{displaymath}
  R^p_i = \dim E_{p,i} - \dim E_{p,i+1}
\end{displaymath}
where $E_p = E_{p,1} \supset \cdots \supset E_{p,s_p+1} = 0$ is
the parabolic filtration at $p$.  Thus $R^p_i = m_i(p)$ in the
notation of Subsection~\ref{sec:parabolic-definitions}. We choose
to keep Holla's original notation in this section because it is
better suited for the calculations to be carried out.  The
\emph{rank} of $R$ is just the rank of $E$,
\begin{math}
  n(R) = \sum_i R^p_i
\end{math}.
One defines
\begin{math}
  \alpha(R) = \sum_{p, i}\alpha_i(p)R^p_i
\end{math},
so that
\begin{displaymath}
  \pdeg(E) = \deg(E) + \alpha(R).
\end{displaymath}
Given a parabolic bundle $E$ with  Harder--Narasimhan filtration
\begin{displaymath}
  0 = G_0 \subsetneq G_1 \subsetneq \cdots \subsetneq G_r = E,
\end{displaymath}
each subbundle $G_j$ is a parabolic bundle with the induced
parabolic structure.  The induced quasi-parabolic data $R^I_{\leq
k}$ is defined by
\begin{displaymath}
  (R^I_{\leq k})^p_i = \dim (G_{k,p} \cap E_{p,i})
  - \dim (G_{k,p} \cap E_{p,i+1}).
\end{displaymath}
Thus $(R^I_{\leq
  k})^p_i$ is the multiplicity of the $i$-th step of the induced
parabolic structure on $G_k$ at $p$ (note that this may be zero).
Each subquotient $G_k/G_{k-1}$ is also a parabolic bundle and the
corresponding parabolic data is $R^I_k$, given by
\begin{displaymath}
  (R^I_k)^p_i = (R^I_{\leq k})^p_i - (R^I_{\leq k-1})^p_i.
\end{displaymath}
The \emph{intersection matrix} $I$ is defined by letting
\begin{displaymath}
  I^p_{i,k} = (R^I_k)^p_i,
\end{displaymath}
in other words, $I^p_{i,k}$ is the multiplicity of the $i$-th step
of the induced parabolic structure on $G_k/G_{k-1}$ at $p$. The
rank of the subquotient $G_k/G_{k-1}$ is
\begin{math}
  n(R^I_k) = \sum_i I^p_{i,k}
\end{math}
and hence the Harder--Narasimhan type of $E$ can be written as
\begin{displaymath}
  \mathbf{n}=(n_1,\ldots,n_r)
   = \bigl(n(R^I_1),\ldots,n(R^I_r)\bigr).
\end{displaymath}

\subsection{Holla's formula}
The formula \cite[Theorem~5.23]{Ho} for the Poincar\'e polynomial
of the moduli space of parabolic bundles of degree $\Delta$ and
rank $n(R)$ is
\begin{equation}
  \label{eq:1}
  P_t(\Delta,n(R))
  = (1-t^2)\sum_{r=1}^{n(R)}\sum_I
  \frac{t^{2\{\sigma'(I) - \Delta(n(R)-n(R_r^I)) + M_g(I,\alpha)\}}}
  {(t^{2n(R^I_1) + 2n(R^I_2)}-1) \cdots (t^{2n(R^I_{r-1}) + 2n(R^I_r)}-1)}
  \prod_{k=1}^rP_{R^I_k(t)},
\end{equation}
where the sum is over intersection matrices $I$ of all possible
Harder--Narasimhan filtrations of parabolic bundles. Here $r$ is
the length of the Harder--Narasimhan filtration corresponding to
$I$, the number $\sigma'(I)$ is defined by $\sigma'(I) = \sum_{p
\in D} \sigma'_p(I)$, where
\begin{displaymath}
  \sigma'_p(I) = \sum_{k>l,\ i<j} I_{i,k}^pI_{j,l}^p,
\end{displaymath}
the number $M_g(I,\alpha)$ is defined by
\begin{multline*}
  M_g(I,\alpha) =
  \sum_{k=1}^{r-1}
  \big(n(R^I_k)+n(R^I_{k+1})\big)
  \left(\left[n(R^I_{\leq k})
    \frac{\Delta+\alpha(R)}{n(R)}
    - \alpha(R^I_{\leq k})\right]+1\right) \\
    + (g-1)\sum_{i<j}n(R^I_i)n(R^I_j)
\end{multline*}
and
\begin{multline*}
  P_R(t) = \left(
    \frac{\prod_{i=1}^{n(R)}(1-t^{2i})^n}
      {\prod_{p \in D}\prod_{\{i\,|\, R^p_i \neq0 \}}
        \prod_{l=1}^{R^p_i}(1-t^{2l})}\right)
    \left(\frac{\prod_{i=1}^{n(R)}(1+t^{2i-1})^{2g}}
      {(1-t^{2n(R)})\prod_{i=1}^{n(R)-1}(1-t^{2i})^2}\right).
\end{multline*}
This formula is valid for all choices of weights such that a
parabolically semistable bundle is automatically parabolically
stable. In particular, it is valid under our
Assumption~\ref{assumption} on genericity of the weights.

\begin{rem}\label{rem:holla-non-fixed}
  This is the formula for the non-fixed determinant case, whereas
  Holla states the formula for the fixed determinant case.  The two
  formulas differ by a factor of $(1+t)^{2g}$ coming from the
  Poincar\'e polynomial of the Jacobian (see \cite{Ni}).
\end{rem}

\subsection{The rank $3$ case}
We now work out explicitly Holla's formula for the case of rank
$3$ parabolic bundles under Assumptions \ref{assumption} and
\ref{ass:small-weights}. This implies in particular that all
parabolic flags are full. For full flags we have the following
simplification of the expressions $P_{R^I_k}(t)$.

\begin{prop}
  Assume that all parabolic flags are full.  Then
  $$
    P_{R^I_k}(t) =
    \frac{\prod_{i=1}^{n(R^I_k)}(1-t^{2i})^{n-1}(1+t^{2i-1})^{2g}}
    { (1-t^2)^{n(R^I_k)n}
      \prod_{i=1}^{n(R^I_k)-1}(1-t^{2i})}\, .
  $$
\end{prop}

\begin{proof}
  Since all flags are full, we have $\#\{i \,|\, (R^I_k)^p_i \neq 0\}
  = n(R^I_k)$.
\end{proof}

Note that $P_{R^I_k}(t)$ only depends on $I$ through the rank
$n(R_k^I)$ of $G_k/G_{k-1}$.  For $n_k = n(R_k^I)$ we shall
therefore write
\begin{displaymath}
  P_{n_k}(t) = P_{R^I_k}(t).
\end{displaymath}

We calculate $P_{n_k}$ for $n_k=1,2,3$.  We obtain
\begin{equation}
  \label{eq:678}
  \begin{split}
  P_1(t) &= \frac{(1+t)^{2g}}{(1-t^2)}\, , \\
  P_2(t) &=
  \frac{(1+t^2)^{n-1}(1+t)^{2g}(1+t^3)^{2g}}{(1-t^2)^3}\, , \\
  P_3(t) &=
  \frac{(1+2t^2+2t^4+t^6)^{n-1}
    (1+t)^{2g}(1+t^3)^{2g}(1+t^5)^{2g}}{(1-t^2)^4(1-t^4)}\, .
 \end{split}
\end{equation}

Now rewrite (\ref{eq:1}) as
\begin{equation}\label{eq:42}
  P_t(\Delta,n(R))
  = (1-t^2) \sum_{\mathbf{n}} \sum_{\text{$I$ of type $\mathbf{n}$}}
  t^{2 \sigma'(I)} Q_I(t) \prod_{k=1}^r P_{n_k}(t),
\end{equation}
where
 $$
  Q_I(t) =\frac{t^{2\{M_g(I,\alpha) -\Delta(n(R)-n(R_r^I)) \}}}
  {(t^{2n(R^I_1) + 2n(R^I_2)}-1) \cdots (t^{2n(R^I_{r-1}) +
  2n(R^I_r)}-1)}\, .
 $$

For rank $n(R)=3$, the possible Harder--Narasimhan types
$\mathbf{n}$ are $(3)$, $(1,2)$, $(2,1)$ and $(1,1,1)$. In the
following, we list all the possible intersection matrices $I$
according to the various types for rank $3$. We also give the
corresponding values of $\sigma'_p(I)$.

\subsubsection*{Intersection matrix for type $(3)$}

\begin{center}
\begin{tabular}{|l|c|}
\hline
  $I_{i,k}^p$ & $k=1$ \\
\hline
  $i=1$ & $1$ \\
  $i=2$ & $1$ \\
  $i=3$ & $1$ \\
  \hline
  $\sigma'_p(I)$ & $0$ \\
  \hline
\end{tabular}
\end{center}

\subsubsection*{Intersection matrices for type $(1,2)$}

\begin{center}
\begin{tabular}{|l|cc|}
\hline
  $I_{i,k}^p$ & $k=1$ & $k=2$\\
\hline
  $i=1$ & $1$ & $0$\\
  $i=2$ & $0$ & $1$ \\
  $i=3$ & $0$ & $1$ \\
  \hline
  $\sigma'_p(I)$ & \multicolumn{2}{c}{$0$} \vline\\
  \hline
\end{tabular}
\quad
\begin{tabular}{|l|cc|}
\hline
  $I_{i,k}^p$ & $k=1$ & $k=2$\\
\hline
  $i=1$ & $0$ & $1$\\
  $i=2$ & $1$ & $0$ \\
  $i=3$ & $0$ & $1$ \\
  \hline
  $\sigma'_p(I)$ & \multicolumn{2}{c}{$1$} \vline\\
  \hline
\end{tabular}
\quad
\begin{tabular}{|l|cc|}
\hline
  $I_{i,k}^p$ & $k=1$ & $k=2$\\
\hline
  $i=1$ & $0$ & $1$\\
  $i=2$ & $0$ & $1$ \\
  $i=3$ & $1$ & $0$ \\
  \hline
  $\sigma'_p(I)$ & \multicolumn{2}{c}{$2$} \vline\\
  \hline
\end{tabular}

\end{center}

\subsubsection*{Intersection matrices for type $(2,1)$}

\begin{center}
\begin{tabular}{|l|cc|}
\hline
  $I_{i,k}^p$ & $k=1$ & $k=2$\\
\hline
  $i=1$ & $1$ & $0$\\
  $i=2$ & $1$ & $0$ \\
  $i=3$ & $0$ & $1$ \\
  \hline
  $\sigma'_p(I)$ & \multicolumn{2}{c}{$0$} \vline\\
  \hline
\end{tabular}
\quad
\begin{tabular}{|l|cc|}
\hline
  $I_{i,k}^p$ & $k=1$ & $k=2$\\
\hline
  $i=1$ & $1$ & $0$\\
  $i=2$ & $0$ & $1$ \\
  $i=3$ & $1$ & $0$ \\
  \hline
  $\sigma'_p(I)$ & \multicolumn{2}{c}{$1$} \vline\\
  \hline
\end{tabular}
\quad
\begin{tabular}{|l|cc|}
\hline
  $I_{i,k}^p$ & $k=1$ & $k=2$\\
\hline
  $i=1$ & $0$ & $1$\\
  $i=2$ & $1$ & $0$ \\
  $i=3$ & $1$ & $0$ \\
  \hline
  $\sigma'_p(I)$ & \multicolumn{2}{c}{$2$} \vline\\
  \hline
\end{tabular}

\end{center}

\subsubsection*{Intersection matrices for type $(1,1,1)$}

\begin{center}
\begin{tabular}{|l|ccc|}
\hline
  $I_{i,k}^p$ & $k=1$ & $k=2$ & $k=3$\\
\hline
  $i=1$ & $1$ & $0$ & $0$ \\
  $i=2$ & $0$ & $1$ & $0$\\
  $i=3$ & $0$ & $0$ & $1$\\
  \hline
  $\sigma'_p(I)$ & \multicolumn{3}{c}{$0$} \vline\\
  \hline
\end{tabular}
\quad
\begin{tabular}{|l|ccc|}
\hline
  $I_{i,k}^p$ & $k=1$ & $k=2$ & $k=3$\\
\hline
  $i=1$ & $1$ & $0$ & $0$ \\
  $i=2$ & $0$ & $0$ & $1$\\
  $i=3$ & $0$ & $1$ & $0$\\
  \hline
  $\sigma'_p(I)$ & \multicolumn{3}{c}{$1$} \vline\\
  \hline
\end{tabular}

\medskip
\begin{tabular}{|l|ccc|}
\hline
  $I_{i,k}^p$ & $k=1$ & $k=2$ & $k=3$\\
\hline
  $i=1$ & $0$ & $1$ & $0$ \\
  $i=2$ & $1$ & $0$ & $0$\\
  $i=3$ & $0$ & $0$ & $1$\\
  \hline
  $\sigma'_p(I)$ & \multicolumn{3}{c}{$1$} \vline\\
  \hline
\end{tabular}
\quad
\begin{tabular}{|l|ccc|}
\hline
  $I_{i,k}^p$ & $k=1$ & $k=2$ & $k=3$\\
\hline
  $i=1$ & $0$ & $1$ & $0$ \\
  $i=2$ & $0$ & $0$ & $1$\\
  $i=3$ & $1$ & $0$ & $0$\\
  \hline
  $\sigma'_p(I)$ & \multicolumn{3}{c}{$2$} \vline\\
  \hline
\end{tabular}

\medskip
\begin{tabular}{|l|ccc|}
\hline
  $I_{i,k}^p$ & $k=1$ & $k=2$ & $k=3$\\
\hline
  $i=1$ & $0$ & $0$ & $1$ \\
  $i=2$ & $1$ & $0$ & $0$\\
  $i=3$ & $0$ & $1$ & $0$\\
  \hline
  $\sigma'_p(I)$ & \multicolumn{3}{c}{$2$} \vline\\
  \hline
\end{tabular}
\quad
\begin{tabular}{|l|ccc|}
\hline
  $I_{i,k}^p$ & $k=1$ & $k=2$ & $k=3$\\
\hline
  $i=1$ & $0$ & $0$ & $1$ \\
  $i=2$ & $0$ & $1$ & $0$\\
  $i=3$ & $1$ & $0$ & $0$\\
  \hline
  $\sigma'_p(I)$ & \multicolumn{3}{c}{$3$} \vline\\
  \hline
\end{tabular}
\end{center}

Now we compute the exponent of the power of $t^{2}$ in the
numerator of $Q_{I}(t)$. This is greatly simplified thanks to our
assumption on the degree and weights.

\begin{prop}
  \label{prop:weight-dependence}
  Let $n(R) = 3$. Under Assumption \ref{ass:small-weights}, we have
\begin{align*}
  M_g(I,\alpha) - \Delta(n(R)-n(R_r^I)) =
  \begin{cases}
    3\left[\frac{\Delta}{3}\right]-\Delta +2g+1\, ,
    &\text{for $I$ of type $(1,2)$,} \\
    3\left[\frac{2\Delta}{3}\right]-2\Delta +2g+1\, ,
    &\text{for $I$ of type $(2,1)$,} \\
    3g-1\, ,
    &\text{for $I$ of type $(1,1,1)$,}
  \end{cases}
\end{align*}
\end{prop}

\begin{proof}
  As $\Delta \not\equiv 0 \pmod 3$, we have that
  $n(R^I_{\leq k}) \frac{\Delta}{n(R)}$ is non-integer.
  Since the weights are small, we have that,
  for $k \leq r-1$,
  $$
     \left[n(R^I_{\leq k}) \frac{\Delta+\alpha(R)}{n(R)}
    - \alpha(R^I_{\leq k})\right] + 1
    =  \left[n(R^I_{\leq k}) \frac{\Delta}{n(R)}\right] + 1.
  $$
Substituting into the definition of $M_g(I,\alpha)$, we get the
result. For type $(1,1,1)$, we have used that
$\left[\frac{\Delta}{3}\right] +\left[\frac{2\Delta}{3}\right]-
\Delta = -1$, since $\Delta \not\equiv 0 \pmod 3$.
\end{proof}

Note that, in particular, Proposition~\ref{prop:weight-dependence}
implies that $Q_{I}(t)$ only depends on $I$ through its
Harder--Narasimhan type.  We shall therefore need to calculate
$\sum_{I}t^{2\sigma'(I)}$ for each type.  This is an easy task
using the tables given for the intersection matrices and the fact
that
\begin{math}
  \sum_{I}t^{2\sigma'(I)}
  = \prod_{p \in D} \sum_{I^{p}}t^{2\sigma_p'(I)}
\end{math},
as is easily seen by induction on the number of points in $D$. The
result is:
 \begin{equation}\label{eq:1011}
\begin{aligned}
  \sum_{I}t^{2\sigma'(I)} &= (1+t^{2}+t^{4})^{n} \, ,
    &&\text{for $I$ of type $(1,2)$ and $(2,1)$,} \\
  \sum_{I}t^{2\sigma'(I)} &= (1+2t^{2}+2t^{4}+t^{6})^{n} \, ,
    &&\text{for $I$ of type $(1,1,1)$.}
\end{aligned}
\end{equation}

We can now calculate the contribution to $P_{t}(\Delta,3)$ from
$I$ of type $(1,1,1)$:
\begin{equation}
  \sum_{\text{$I$ of type $(1,1,1)$}}
    t^{2\sigma'(I)} Q_{I}(t) \prod_{k=1}^{r} P_{n_{k}}(t) \label{eq:14}
  = \frac{(1+2t^{2}+2t^{4}+t^{6})^{n}t^{6g-2}}
     {(t^{4}-1)^{2}}P_{1}(t)^{3}.
\end{equation}
Similarly, the contributions to $P_{t}(\Delta,3)$ from $I$ of type
$(1,2)$ and $(2,1)$ are:
\begin{align*}
  \sum_{\text{$I$ of type $(1,2)$}}
     t^{2\sigma'(I)} Q_{I}(t) \prod_{k=1}^{r} P_{n_{k}}(t)
  &= \frac{(1+t^{2}+t^{4})^{n}
     t^{2\left\{3\left[\frac{\Delta}{3}\right]-\Delta +2g+1\right\}}}
     {t^{6}-1}P_{1}(t)P_{2}(t), \\
  \sum_{\text{$I$ of type $(2,1)$}}
    t^{2\sigma'(I)} Q_{I}(t) \prod_{k=1}^{r} P_{n_{k}}(t)
  &= \frac{(1+t^{2}+t^{4})^{n}
     t^{2\left\{3\left[\frac{2\Delta}{3}\right]-2\Delta +2g+1\right\}}}
     {t^{6}-1}P_{1}(t)P_{2}(t).
\end{align*}
Summing the contributions of type $(1,2)$ and type $(2,1)$ some
simplification results because, whenever $\Delta \not\equiv 0
\pmod 3$, one has
\begin{align*}
  t^{2\{3\left[\frac{\Delta}{3}\right]-\Delta +2g+1\}}
  + t^{2\{3\left[\frac{2\Delta}{3}\right]-2\Delta +2g+1\}}
  &= t^{4g+2}(t^{2\{3\left[\frac{\Delta}{3}\right]-\Delta\}}
     + t^{2\{3\left[\frac{2\Delta}{3}\right]-2\Delta\}}) \\
  &= t^{4g+2}(t^{-2}+t^{-4}) \\
  &= t^{4g-2}(1+t^{2})\ .
\end{align*}
Hence we obtain
\begin{equation}\label{eq:17}
  \sum_{\text{$I$ of type $(1,2)$ or $(2,1)$}}
   \hspace{-5mm} t^{2\sigma'(I)} Q_{I}(t) \prod_{k=1}^{r} P_{n_{k}}(t)
      = \frac{(1+t^{2}+t^{4})^{n}
     t^{4g-2}(1+t^{2})}
     {t^{6}-1}P_{1}(t)P_{2}(t).
\end{equation}

Since for $\mathbf{n}=(3)$ we clearly have
\begin{equation}\label{eq:18}
  t^{2\sigma'(I)} Q_{I}(t) = 1,
\end{equation}
we are now in a position to put everything together and calculate
$P_{t}(\Delta,3)$ for $\Delta \not\equiv 0 \pmod 3$.

\begin{prop}
  \label{prop:19}
  Under Assumption \ref{ass:small-weights},
  the Poincar\'e polynomial of the moduli space of stable
  parabolic bundles of rank $3$ is given by
\begin{multline*}
    P_{t}(\Delta,3) =(1+t)^{2g}(1+2t^2+2t^4+t^6)^{n-1} \cdot \\
      \cdot \frac{(1+t^3)^{2g}(1+t^5)^{2g}+
       (1+t)^{4g}(1+t^{2}+t^{4})t^{6g-2}
       -(1+t)^{2g}(1+t^3)^{2g}(1+t^2)^2t^{4g-2}}
     {(1-t^2)^3(1-t^4)}\ .
\end{multline*}
\end{prop}

\begin{proof}
  Substituting (\ref{eq:14}), (\ref{eq:17}) and  (\ref{eq:18})  in
  (\ref{eq:42}) and using (\ref{eq:678}) we obtain
  \begin{displaymath}
  \begin{split}
    P_{t}(\Delta,3) &=
    (1-t^{2})\Bigl(
    P_{3}(t)
    + \frac{(1+t^{2}+t^{4})^{n}t^{4g-2}(1+t^{2})}
     {t^{6}-1}P_{1}(t)P_{2}(t) \\
    &\qquad + \frac{(1+2t^{2}+2t^{4}+t^{6})^{n}t^{6g-2}}
     {(t^{4}-1)^{2}}P_{1}(t)^{3} \Bigr) \\
     &=
      \frac{(1+2t^2+2t^4+t^6)^{n-1}
        (1+t)^{2g}(1+t^3)^{2g}(1+t^5)^{2g}}{(1-t^2)^3(1-t^4)} \\
     &\qquad + \frac{(1+t^{2}+t^{4})^{n}t^{4g-2}(1+t^{2})(1+t)^{2g}
        (1+t^2)^{n-1}(1+t)^{2g}(1+t^3)^{2g}}
     {(t^{6}-1)(1-t^2)^3} \\
     &\qquad + \frac{(1+2t^{2}+2t^{4}+t^{6})^{n}t^{6g-2}(1+t)^{6g}}
     {(1-t^{4})^{2}(1-t^2)^2}\ .
  \end{split}
  \end{displaymath}
Simplifying this expression we obtain the formula stated.
\end{proof}

\section{Betti numbers of the  moduli space
of rank three parabolic Higgs bundles}\label{sec:rank3}

In this section we put everything together to obtain the
Poincar\'e polynomial of the moduli space of rank three parabolic
Higgs bundles.

\subsection{Poincar\'e polynomial}

\begin{thm} \label{thm:end}
  Let $\mathcal{M}$ be the moduli space of rank three parabolic Higgs
  bundles of some fixed degree and weights, over a connected, smooth
  projective complex algebraic curve of genus $g$. If the
  weights are generic (in the sense that there are no properly
  semistable parabolic Higgs bundles), then the Poincar\'e polynomial
  of $\mathcal{M}$ is given by
\begin{multline*}
  P_t(\mathcal{M} ) =
\Coeff_{u^0v^0}  \bigg( (1+2u^2vt^2+2uv^2t^2+u^3v^3t^4)^n  \,
 \cdot \\ \cdot \,
  \frac{t^{2(4g-3+n)}(1+t)^{2g}
  (1+u^2vt)^{2g}(1+uv^2t)^{2g}}
  {u^{3n+6g-8}v^{3n+6g-7}
    (1-u^2v)(1-uv^2)(1-u^2vt^2)(1-uv^2t^2)
    (1-v^3t^2)(1-u^3t^2)} \bigg) \\
 +\Coeff_{x^0} \Bigg( \frac{(1+t)^{4g}(1+xt)^{2g}
 t^{6g-6} x^{2-2g- n}}{(1-t^2)
 (1-x)(1-xt^2)}\cdot \\
 \cdot\bigg(\frac{t^{2g-2+2n}(x+x^2)(1+2t^2+2t^{2}x+t^4x)^n}
   {(1-t^{-2}x)(1-t^2x^2)}
  -\frac{(t^4x+t^8x^2)(1 + 2t^{2} +2t^4x +t^6x )^n}{(1-t^4x)(1-t^8x^2)}
 \bigg)
 \Bigg) \\
 + (1+t)^{2g}(1+2t^2+2t^4+t^6)^{n-1} \cdot \\
      \cdot \frac{(1+t^3)^{2g}(1+t^5)^{2g}+
       (1+t)^{4g}(1+t^{2}+t^{4})t^{6g-2}
       -(1+t)^{2g}(1+t^3)^{2g}(1+t^2)^2t^{4g-2}}
     {(1-t^2)^3(1-t^4)}\ .
\end{multline*}
\end{thm}

\begin{proof}
  It follows from Morse theory, as explained in
  Section~\ref{sec:morse}, that $P_t(\mathcal{M} ) = P_t(\Delta,3) +
  P_t(\Delta,(1,2)) + P_t(\Delta,(2,1)) + P_t(\Delta,(1,1,1))$, where
  the polynomials on the right hand side are given in
  Propositions~\ref{prop:P_t111_nonzero}, \ref{prop:P_t(1,2)_nonzero},
  \ref{prop:P_t(2,1)_nonzero} and \ref{prop:19}.  In order to apply
  these formulas we need to choose $\Delta_0 \in \{1,2\}$.  However,
  the contribution $P_t(\Delta,(1,1,1))$ of Proposition
  \ref{prop:P_t111_nonzero} is independent of this choice, by using
  the duality $(u,v)\mapsto (v,u)$.  Also, the contribution
  $P_t(\Delta,(1,2))+P_t(\Delta,(2,1))$ from Propositions
  \ref{prop:P_t(1,2)_nonzero} and \ref{prop:P_t(2,1)_nonzero} is
  independent of the choice of $\Delta_0$ by using that for
  $\Delta_0=1,2$ we have $x^{\Delta_0}+x^{3-\Delta_0} = x+x^2$ and
  $x^{\Delta_0}t^{4\Delta_0}+x^{3-\Delta_0}t^{12-4\Delta_0} =
  t^4x+t^8x^2$.

  Even though the various contributions to the Poincar\'e polynomial
  were calculated for a specific choice of degree and weights (cf.\
  Assumptions~\ref{assumption} and~\ref{ass:small-weights}), we know
  by Proposition \ref{prop:deg-tensor} that the final result is
  independent of this choice.
\end{proof}

\begin{rem}
  Obviously, it is possible to do the computation of
  $P_t(\mathcal{M})$ under different choices of degree and
  weights. Most of the calculations in Sections
  \ref{sec:(1,1,1)}--\ref{sec:rank3} are carried out in general,
  and we have always introduced our Assumption \ref{ass:small-weights}
  as late as possible in each section. Of course, the final
  answer will be the same as the one given in Theorem \ref{thm:end}, though
  the partial contributions of the critical submanifolds of
  different types may differ.
\end{rem}

\subsection{Special low genus cases}

We can calculate the Poincar\'e polynomial of the moduli space of
parabolic Higgs bundles for specific values of $n$ and $g$ by using a
computer algebra system. For instance, if $n=1$ and $g=2$ then
Theorem~\ref{thm:end} gives
\begin{multline*}
  P_t(\mathcal{M} ) =
  \,36\,{t}^{26}+324\,{t}^{25}+1368\,{t}^{24}+3620\,{t}^{23}+6810\,{t}^{22}
  +9860\,{t}^{21}+11670\,{t}^{20}+11876\,{t}^{19}\\
  +10860\,{t}^{18}+9224\,{t}^{17}+7408\,{t}^{16}+5688\,{t}^{15}
  +4216\,{t}^{14}+3036\,{t}^{13}+2134\,{t}^{12}+1464\,{t}^{11} \\
  +981\,{t}^{10}+640\,{t}^{9}+401\,{t}^{8}+244\,{t}^{7}+144\,{t}^{6}
  +80\,{t}^{5}+42\,{t}^{4}+20\,{t}^{3}+9\,{t}^{2}+4\,t+1\ ,
\end{multline*}
and when $n=2$ and $g=2$ we obtain
\begin{multline*}
  P_t(\mathcal{M} ) =
  \,252\,{t}^{32}+2416\,{t}^{31}+10848\,{t}^{30}+30540\,{t}^{29}
  +61178\,{t}^{28}+94368\,{t}^{27} \\ +119187\,{t}^{26}
  +129952\,{t}^{25}
  +127737\,{t}^{24}+116656\,{t}^{23}+100849\,{t}^{22}\\+83564\,{t}^{21}
  +66925\,{t}^{20} +52100\,{t}^{19}+39605\,{t}^{18}
  +29504\,{t}^{17}+21572\,{t}^{16} \\ +15472\,{t}^{15}+10884\,{t}^{14}
  +7496\,{t}^{13}+5043\,{t}^{12}+3312\,{t}^{11}+2113\,{t}^{10} \\
  +1308\,{t}^{9}
  +782\,{t}^{8}+448\,{t}^{7}+247\,{t}^{6}+128\,{t}^{5}
  +62\,{t}^{4}+28\,{t}^{3}+11\,{t}^{2}+4\,t+1\ .
\end{multline*}
Another example is the Poincar\'e polynomial of $\mathcal{M}$ for
$g=1$ and $n=1$,
\begin{displaymath}
  P_t(\mathcal{M} ) =
  \,6\,{t}^{8}+18\,{t}^{7}+24\,{t}^{6}+20\,{t}^{5}+13\,{t}^{4}
  +8\,{t}^{3}+4\,{t}^{2}+2\,t+1\ .
\end{displaymath}


\begin{rem}\label{rem:hausel}
  In \cite{Ha2} Hausel conjectured a formula for the Poincar\'e
  polynomial of the moduli space of stable Higgs bundles of any rank.
  Hausel has informed us of an analogous conjecture for the the
  Poincar\'e polynomial of the moduli space of {\em parabolic}
  Higgs bundles of any rank: in the case of rank three the mixed
  Hodge polynomial of the corresponding character variety is
  \begin{eqnarray*}
 H^n_3(q,t) &=&
  \frac{((q t^2+1) (q^2 t^4+q t^2+1))^n (q^3t^5+1)^{2g}
  (q^2t^3+1)^{2g}}{(q^3t^6-1)(q^3t^4-1)(q^2t^4-1)(q^2t^2-1)} \\
 & & + \,  \frac{ (q^3t^6(q+1)(q^2+q+1))^n  q^{6g-6} t^{12g-12}
 (q^3 t+1)^{2g}
 (q^2t+1)^{2g}}{(q^3t^2-1)(q^3-1)(q^2t^2-1)(q^2-1)}\\
 & &+ \,\frac{(q^2t^4 (2q^2t^2+qt^2+q+2))^n q^{4g-4} t^{8g-8}
 (q^3t^3+1)^{2g}(q t+1)^{2g}}{(q^3t^4-1)(q^3t^2-1)(q
 t^2-1)(q-1)} \\
 & & + \,\frac{6^n (q t^2)^{3 n} q^{6 g-6} t^{12 g-12}
 (q t+1)^{4g}}{3(q t^2-1)^2 (q-1)^2}\\
 & & -\,\frac{(3 q^2 t^4 (q t^2+1))^n q^{4 g-4} t^{8 g-8}
 (q^2 t^3+1)^{2 g}(q t+1)^{2 g}}{(q^2 t^4-1) (q^2 t^2-1) (q t^2-1) (q-1)} \\
 & & -\,\frac{(3 q^3 t^6 (q+1))^n q^{6 g-6} t^{12 g-12} (q^2 t+1)^{2 g}
 (q t+1)^{2 g}}{(q^2 t^2-1) (q^2-1) (q t^2-1) (q-1)} \, .
 \end{eqnarray*}
  Hausel conjectures that the Poincar\'e polynomial of the moduli space
  of rank three parabolic Higgs bundles with $n$ marked points and full
  flags is obtained from this polynomial by the substitution
  $P_t(\mathcal{M} ) =H^n_3(1,t)$.
  The formulas are difficult to compare in general, but in
  computer calculations Hausel's formula provides the same result as
  ours in all the cases that we have checked. Thus our results
  provide evidence for this conjecture.
\end{rem}


We finish by considering the case when $X$ has genus zero. It is
easy to see that under our Assumption \ref{ass:small-weights} of small
weights, the moduli space of stable parabolic bundles on $X$ is empty,
because any parabolically stable bundle would have to be stable. However,
our results show that there are non-empty critical submanifolds of the
moduli space of parabolic Higgs bundles for $n \geq 3$: for example,
if $g=0$ and $n=3$, our calculations show that
\begin{displaymath}
  P_t(\mathcal{M} ) = \,7\,t^2 + 1\ ,
\end{displaymath}
where the only contribution is from critical submanifolds of type
$(1,1,1)$. This means that stable parabolic Higgs bundles exist and
hence the moduli space is non-empty.  Following our general
description of the critical submanifolds of
Subsection~\ref{sec:(1,1,1)-description} one can explicitly describe the
critical submanifolds of type $(1,1,1)$ and thus get examples of
stable parabolic Higgs bundles.  One such example is given as follows.
Define parabolic line bundles $L_1 = \mathcal{O}(1)$, $L_2 =
\mathcal{O}$ and $L_3 = \mathcal{O}$ with small weights $\alpha_i(p)$
on $L_i$, such that $\alpha_1(p) < \alpha_2(p) < \alpha_3(p)$ at each
marked point.  Let $E = L_1 \oplus L_2 \oplus L_3$.  Then any map from
$L_i$ to $L_{i+1}$ is strongly parabolic and we can define a Higgs
field $\Phi$ with non-zero components in $\SPH(L_1,L_{2}\otimes K(D))
= \Hom(L_1,L_2\otimes K(D)) = \Gamma(\mathcal{O}(n-3))$ and
$\SPH(L_2,L_3 \otimes K(D)) = \Hom(L_2,L_3\otimes K(D)) =
\Gamma(\mathcal{O}(n-2))$.  Clearly the resulting parabolic Higgs
bundle is stable.  When $n=3$ this parabolic Higgs bundle is a minimum
of the Morse function, as follows from Lemma~\ref{lem:index-111}, and
the critical submanifold consisting of such parabolic Higgs bundles is
easily seen to be isomorphic to $\mathbb{P}^1$.  {}From our
description of the critical submanifolds of
Subsection~\ref{sec:(1,1,1)-description} one sees that there are six
other critical submanifolds, all consisting of parabolic Higgs bundles
with the same underlying vector bundle but with different
distributions of the weights.  All these other critical submanifolds
consist of one point and have index 2.  Of course these observations
check with our calculation of the Poincar\'e polynomial.

One can give a very explicit description of the moduli space in the
$n=3$ and $\Delta=0$ case (this is of course different from the
$\Delta=1$ moduli space considered in the previous paragraph but, as
we know, has the same Betti numbers).  This is done by means of the
Hitchin map (\cite{H2}), which exhibits the moduli space as an
elliptic fibration over $\CC$ (in fact an ALG manifold \cite{CK}).  To carry
this out, consider the general case where the bundle over $\PP^1$ is
trivial and the three points are $0$, $1$ and $\infty$. The Higgs
field (twisting by $K(3)=\cO(1)$) can be written as
$$
\Phi=Az+B(z-1)
$$
where $A$, $B$ and $A+B$ are nilpotent since these are the residues at the
parabolic points. That means that $\Tr \Phi=0$, since $\Tr A=\Tr B=0$, and
 $\Tr \Phi^2=0$ since $\Tr A^2=\Tr B^2=0$ and $\Tr(A+B)^2=0$,
which means $\Tr AB=0$; also
$\Tr \Phi^3=cz(z-1)$ using that  $\Tr A^3=\Tr B^3=0$ and $\Tr(A+B)^3=0$.
The spectral curve (\cite{H2}) has the form
$$
w^3=kz(z-1)
$$
which is a cubic curve invariant by $\ZZ/3$ by multiplying $w$ by a
cube root of unity. As $k$ varies in $\CC$ we have the elliptic
fibration with an $E_6$ curve at $k=0$.  Of course the Hitchin map is
just $(E,\Phi) \mapsto k$ and, in particular, the nilpotent cone is
the $E_6$ curve.

For higher values of $n$, there are also contributions from critical
submanifolds of type $(1,2)$ and $(2,1)$. For instance, for $g=0$ and
$n=4$ our formula gives
\begin{displaymath}
  P_t(\mathcal{M} ) =
  \,271\,{t}^{8}+144\,{t}^{6}+43\,{t}^{4}+9\,{t}^{2}+1\ ,
\end{displaymath}
with non-zero contributions from critical submanifolds of type
$(1,2)$.  For $g=0$ and $n=5$ critical submanifolds of both type
$(1,2)$ and $(2,1)$ contribute and one obtains
\begin{displaymath}
  P_t(\mathcal{M} ) =
  \,4645\,{t}^{14}+3791\,{t}^{12}+1926\,{t}^{10}+762\,{t}^{8}
  +249\,{t}^{6}+63\,{t}^{4}+11\,{t}^{2}+1\ .
\end{displaymath}

\section{The fixed determinant case}
 \label{sec:fixed-det}

The goal of this Section is to calculate the Poincar\'e polynonial of
the moduli space of rank $3$ parabolic Higgs bundles with fixed
determinant. We follow our calculation for the non-fixed determinant
case closely and only point out the main differences.  The final
result is given in Theorem~\ref{thm:end-fixed}. As a corollary we
obtain the fact that fixed determinant moduli space has Euler
characteristic zero---note that this is not the case for the usual
fixed determinant Higgs bundle moduli space, cf.\ \cite{H}, \cite{G1}
and \cite{Ha2}.

\subsection{Preliminaries}
\label{sec:preliminaries-fixed}

Let $E$ be a rank $r$ parabolic bundle with degree $\Delta$ and
weights $\a_i(p)$ with multiplicities $m_i(p)$. Then the
determinant $\Lambda^r E$ is a parabolic bundle, with degree
$\bar{\Delta}=\Delta+\sum_{p\in D} \big[\sum_i
m_i(p)\alpha_i(p)\big]$ and weights $\sum_i
m_i(p)\alpha_i(p)-\big[\sum_i m_i(p)\alpha_i(p)\big]$, for $p\in
D$ (in particular, under Assumption \ref{ass:small-weights}, we
have $\bar{\Delta}=\Delta$). Now, for any choice of weights, the
moduli space of rank $1$ parabolic Higgs bundles of degree
$\bar\Delta$ is naturally identified with the total space of
the cotangent bundle to the Jacobian of degree $\bar\Delta$ line
bundles on $X$. Consider the ``determinant map'' from the moduli
space of stable rank $r$ parabolic Higgs bundles $\mathcal{M}$ to
$T^*\Jac^{\bar\Delta}(X)$:
\begin{equation}
  \label{eq:3}
\begin{aligned}
  \det\colon \mathcal{M} &\to T^*\Jac^{\bar\Delta}(X)\ ,\\
  (E,\Phi) &\mapsto (\Lambda^r E, \Tr \Phi)\ .
\end{aligned}
\end{equation}
Let $\Lambda$ be a fixed line bundle of degree $\bar\Delta$. By
definition, the fibre of $\det$ over $(\Lambda,0)$ is the  moduli
space of stable parabolic Higgs bundles with fixed determinant
$\Lambda$:
\begin{displaymath}
  \mathcal{M}^{\Lambda} = {\det}^{-1}(\Lambda,0)\ .
\end{displaymath}

We shall need the following analogue of
Proposition~\ref{prop:deg-tensor}: it is not hard to see that the
proof, including the relevant parts of \cite{T2}, goes over to the
fixed determinant case.

\begin{prop} \label{prop:deg-tensor-fixed}
  Fix the rank $r$. For different choices of the determinant bundle
  $\Lambda$ and generic weights, the moduli spaces
  $\mathcal{M}^{\Lambda}$ have the same Betti numbers.
  \qed
\end{prop}

\begin{rem}\label{rem:cover}
The group of $r$-torsion points in
the Jacobian,
\begin{math}
  \Gamma_r = \{L \;|\; L^r = \mathcal{O} \},
\end{math}
acts on $\mathcal{M}^{\Lambda}$ by tensor product:
\begin{displaymath}
  (E,\Phi) \mapsto (E\otimes L,\Phi)\ .
\end{displaymath}
We also have an action of $\Gamma_r$ on $T^*\Jac^l(X)$ given by
\begin{displaymath}
  (M,\alpha) \mapsto (M\otimes L^{-1},\alpha)\ ,
\end{displaymath}
for $L \in \Gamma_r$ and, via this action, the covering
\begin{displaymath}
  \begin{aligned}
    T^*\Jac^l(X) &\to T^*\Jac^l(X)\ ,\\
    (M,\alpha) &\mapsto (M^r,\alpha)\ ,
  \end{aligned}
\end{displaymath}
can be viewed as a principal $\Gamma_r$-bundle.  As done in
Atiyah--Bott \cite{atiyah-bott:1982} for ordinary bundles, we can use
this to express $\mathcal{M}$ as a fibred product
\begin{displaymath}
  \begin{aligned}
    \mathcal{M}^{\Lambda}\x_{\Gamma_r}T^*\Jac^0
      &\overset{\cong}{\to}\mathcal{M}\ ,\\
    \bigl((E,\Phi),(L,\alpha)\bigr)
      &\mapsto (E\otimes L,\Phi+\alpha\, \Id)\ .
  \end{aligned}
\end{displaymath}
It follows that the rational cohomology of $\mathcal{M}$ is
isomorphic to the $\Gamma_r$-invariant part of the cohomology of
$\mathcal{M}^{\Lambda}\x_{\Gamma_r}T^*\Jac^0$.  But $\Gamma_r$
acts trivially on the cohomology of $T^*\Jac^0$ and, therefore,
\begin{displaymath}
  H^*(\mathcal{M};\QQ) \cong
  H^*(\mathcal{M}^{\Lambda};\QQ)^{\Gamma_r} \otimes H^*(\Jac^0;\QQ)\ ,
\end{displaymath}
where we write
\begin{displaymath}
  H^*(\mathcal{M}^{\Lambda};\QQ)
  = H^*(\mathcal{M};\QQ)^{\Gamma_r}
  \oplus
  H^*(\mathcal{M}^{\Lambda};\QQ)^{\mathrm{var}}
\end{displaymath}
as the direct sum of the $\Gamma_r$-invariant part and the
non-invariant part, or \emph{variant} part in the terminology
of \cite{HT2}.  It follows from this that
\begin{displaymath}
  P_t(\mathcal{M}) = P_t(\mathcal{M}^{\Lambda})(1+t)^{2g}
\end{displaymath}
if and only if $\Gamma_r$ acts trivially on
$H^*(\mathcal{M}^{\Lambda};\QQ)$: in fact,
\begin{equation}\label{eq:2}
  P_t(\mathcal{M}^{\Lambda})(1+t)^{2g}
    -P_t(\mathcal{M})
  =P_t^{\mathrm{var}}(\mathcal{M}^{\Lambda})(1+t)^{2g}\ ,
\end{equation}
where
  $P_t^{\mathrm{var}}(\mathcal{M}^{\Lambda}) = \sum
  t^i\dim(H^i(\mathcal{M}^{\Lambda};\QQ)^{\mathrm{var}})$ is the
  Poincar\'e polynomial corresponding to the variant part of the
  cohomology.
\end{rem}

\subsection{Morse indices}
\label{sec:morse-indices-fixed}

The $S^1$-action on $\mathcal{M}$ restricts to the fixed determinant
moduli space $\mathcal{M}^{\Lambda}$ and the Morse theory explained in
Section~\ref{sec:morse} can be applied to this latter space.  Thus,
the restriction of $f$ to $\mathcal{M}^{\Lambda} \subset \mathcal{M}$
gives a perfect Bott--Morse function.  The characterization of the
critical points of the Morse function (i.e., the fixed points of the
$S^1$-action) and their stability given in
Propositions~\ref{prop:fixed=vhs} and \ref{prop:vhs-stab} remains
valid.  Hence, for each critical submanifold $\mathcal{N} \subseteq
\mathcal{M}$, there is a corresponding critical submanifold
$\mathcal{N}^{\Lambda} \subseteq \mathcal{M}^{\Lambda}$ and the
determinant map \eqref{eq:3} restricts to give a fibration
\begin{equation}\label{eq:5}
  \det\colon \mathcal{N} \to \Jac^{\bar\Delta}(X)
\end{equation}
with fibre over $\Lambda$ equal to $\mathcal{N}^{\Lambda}$.  Note that
there is no need to map to $T^*\Jac^{\bar\Delta}(X)$ because for any
parabolic complex variation of Hodge structure $(\bigoplus E_l, \Phi)$
we have $\Tr \Phi = 0$.

\begin{rem}\label{rem:cover-critical}
  We have a description of $\mathcal{N}$ as a fibred product
  $\mathcal{N}^{\Lambda}\x_{\Gamma_r}\Jac^0$, analogous to the one given
  in Remark~\ref{rem:cover} for $\mathcal{M}$.  Thus we also have an
  analogous description of the relation between the cohomology of
  $\mathcal{N}$ and that of $\mathcal{N}^{\Lambda}$.
\end{rem}

The deformation theory of $\mathbf{E}=(E,\Phi)$ in the fixed
determinant moduli space is governed by the complex
\begin{align*}
  C_0^\bullet (\mathbf{E}) : \quad
  \PE_0(E) &\overset{[-,\Phi]}{\longrightarrow} \SPE_0(E) \otimes K(D) \\
  f &\longmapsto (f \otimes 1)\Phi - \Phi f\ ,
\end{align*}
where the subscript $0$ indicates trace zero (cf.\
Proposition~\ref{prop:deformation}).  Now let $\mathbf{E}=(\bigoplus
E_l,\Phi)$ be a fixed point of the $S^1$-action.  In order to
determine the weight spaces of the infinitesimal circle action on the
tangent space we modify the subcomplexes $C^\bullet(\mathbf{E})_l$
defined in Subsection~\ref{sec:fixed} to be subcomplexes
$C_0^\bullet(\mathbf{E})_l$ of trace zero endomorphisms.  Note that
$C_0^\bullet(\mathbf{E})_l = C^\bullet(\mathbf{E})_l$ unless $l=0$ or
$l=-1$.  The calculation  of the Morse indices now proceeds
analogously to the non-fixed determinant case of
Section~\ref{sec:morse} and, in particular, we obtain the following
result.

\begin{prop}\label{prop-index-fixed}
  Let the parabolic Higgs bundle $\mathbf{E}=(E,\Phi)$ represent a
  critical point of the restriction of $f$ to $\mathcal{M}^{\Lambda}
  \subset \mathcal{M}$. Then the Morse index of $f$ at this point is
  given by the formula of Proposition~\ref{prop-index}.
\end{prop}

\begin{proof}
  The Morse index equals the real dimension of the space
  $\bigoplus_{l>0}\HH^1(C_0^\bullet(E,\Phi)_l)$.  But, as pointed out
  above, $C_0^\bullet(\mathbf{E})_l = C^\bullet(\mathbf{E})_l$ for
  $l>0$.  This proves the proposition.  Alternatively, we could have
  appealed to the invariance of $f$ under the action of the Jacobian
  on $\mathcal{M}$ by tensor product.
\end{proof}

\begin{rem}\label{rem:laumon-fixed}
  The analogue of Theorem~\ref{thm:laumon} also holds in the fixed
  determinant case, with an analogous proof.
\end{rem}

\subsection{Critical submanifolds of type $(1,1,1)$}
\label{sec:type-111-fixed}

In this section we describe the critical submanifolds of type
$(1,1,1)$ and their contribution to the Poincar\'e polynomial in
the fixed determinant case.  We shall use the notations of
Section~\ref{sec:(1,1,1)}. Note that the description given in
Subsection~\ref{sec:(1,1,1)-description} of the parabolic Higgs
bundles which corresponds to critical points of type $(1,1,1)$
remains valid. Likewise, the characterization of stability given
in Proposition~\ref{prop:par-stability111} is the same.  Thus,
fixing $d_1$, $m$ and $\varpi$, the fixed determinant critical
submanifold is the fibre of the map $\det$ defined in
\eqref{eq:5}:
\begin{displaymath}
  \mathcal{N}^{\Lambda}_{(1,1,1)}(d_1,m,\varpi) =
  {\det}^{-1}(\Lambda)\ .
\end{displaymath}
The description of the critical submanifolds now proceeds as in
\cite{G1} (cf.\ also Hausel--Thaddeus \cite{HT2} for the case of
general rank $r$) to give us the following fixed determinant analogue
of Proposition~\ref{prop:boundsmm}.

\begin{prop}
  \label{prop:fixed-critical-111}
  The critical submanifold
  $\mathcal{N}^{\Lambda}_{(1,1,1)}(d_1,m,\varpi)$ is given by the
  pull-back diagram
  \begin{displaymath}
  \begin{CD}
    \mathcal{N}^{\Lambda}_{(1,1,1)}(d_1,m,\varpi)
      @>>> \Jac^{d_3}(X) \\
    @VVV @VVV \\
    S^{m_1}X \times S^{m_2}X @>>> \Jac^{m_1+2m_2}(X)\ ,
  \end{CD}
  \end{displaymath}
  where the vertical map on the left
  is given by $(L_1 \oplus L_2 \oplus L_3, \Phi_1, \Phi_2) \mapsto
  (\mathrm{div}(\Phi_1),\mathrm{div}(\Phi_2))$,
  the map in the bottom line is $(D_1,D_2) \mapsto
  \mathcal{O}(D_1+2D_2)$ and the vertical map on the right is $L_3 \mapsto
  \Lambda^{-1}\otimes L_3^3\otimes K^3(3D-S_1-2S_2)$.
  Moreover, $\mathcal{N}^{\Lambda}_{(1,1,1)}(d_1,m,\varpi)$ is
  non-empty if and only if
  \begin{align*}
    3m &> 2\Delta - F(\alpha,\varpi),\\
    3d_{1} &> \Delta - G(\alpha,\varpi), \\
    2d_{1} -m &\leq n - s_{1} + 2g-2,\\
    2m - d_1 &\leq \Delta  + n - s_2 + 2g-2,
  \end{align*}
  where $F$ and $G$ were defined in \eqref{eq:def-FG}.
\end{prop}

\begin{proof}
  Given a parabolic Higgs bundle $(L_1 \oplus L_2 \oplus L_3, \Phi_1,
  \Phi_2)$ of type $(1,1,1)$, let $M_i$ be the line
  bundle associated to the divisor $D_i = \mathrm{div}(\Phi_i)$:
\begin{align*}
  M_1 &= L_1^{-1}\otimes L_2\otimes K(D-S_1)\ ,\\
  M_2 &= L_2^{-1}\otimes L_3\otimes K(D-S_2)\ .
\end{align*}
Then
\begin{equation}
  \label{eq:8}
  M_1\otimes M_2^2 = \Lambda^{-1}\otimes L_3^3\otimes K^3(3D-S_1-2S_2)\ .
\end{equation}
Conversely, let $(d_1,m,\varpi)$ be such that $m_i \geq 0$ for
$i=1,2$, take effective divisors $D_i$ of degree $m_i$ and define
$M_i= \mathcal{O}(D_i)$.  Then there is a solution $L_3$ to
\eqref{eq:8}, determined up to the choice of a cube root of the
trivial bundle. Once this choice is made, the isomorphism class of
$(L_1 \oplus L_2 \oplus L_3, \Phi_1, \Phi_2)$ can be recovered from
$(D_1,D_2)$. Now, the first two inequalities of the statement of the
Proposition represent the stability condition for $(L_1 \oplus L_2
\oplus L_3, \Phi_1, \Phi_2)$ and the last two inequalities are
equivalent to $m_i \geq 0$ for $i=1,2$.  Thus we see that, for any
$(d_1,m,\varpi)$ satisfying these conditions, there is a non-empty
critical submanifold, as described in the statement
of the Proposition.
\end{proof}

\begin{rem}\label{rem:fixed-critical-111-m1m2}
  We can also parametrize the critical submanifolds by
  $(m_1,m_2,\varpi)$.  We then have critical submanifolds
  $\mathcal{N}^{\Lambda}_{(1,1,1)}(m_1,m_2,\varpi)
  =\mathcal{N}^{\Lambda}_{(1,1,1)}(d_1,m,\varpi)$, which are non-empty
  if and only if
  \begin{equation}
    \label{eq:13}
  \begin{aligned}
    m_1 +2m_2 &< 6g-6 +3n -s_1 -2s_2 + F(\alpha,\varpi)\ ,\\
    2m_1+m_2 &< 6g-6 +3n -2s_1 -s_2 + G(\alpha,\varpi)\ , \\
    m_{1} &\geq 0\ ,\\
    m_{2}&\geq 0\ ,
  \end{aligned}
  \end{equation}
    and
    \begin{equation}
    m_1+2m_2 +\Delta +s_1+2s_2 \equiv 0 \pmod 3\ .\label{eq:12}\
  \end{equation}
  The conditions \eqref{eq:13} are obtained by formulating the
  conditions of the preceding proposition in terms of $m_1$ and $m_2$,
  and the condition \eqref{eq:12} must be added for it to be
  possible to solve \eqref{eq:8} for $L_3$ (as pointed out in
  \cite{HT2}, this condition was overlooked in \cite{G1}).
\end{rem}

As noted in Remark~\ref{rem:cover-critical}, the rational cohomology of
$\mathcal{N}^{\Lambda}_{(1,1,1)}(d_1,m,\varpi)$ splits in an invariant
part and a variant part, under the action of $\Gamma_3$:
\begin{displaymath}
  H^*(\mathcal{N}^{\Lambda}_{(1,1,1)}(d_1,m,\varpi))
  = H^*(\mathcal{N}^{\Lambda}_{(1,1,1)}(d_1,m,\varpi))^{\Gamma_3}
  \oplus
  H^*(\mathcal{N}^{\Lambda}_{(1,1,1)}(d_1,m,\varpi))^{\mathrm{var}}\ .
\end{displaymath}

\begin{prop}\label{prop:cohomology-111}
  The invariant part of the cohomology of
  $\mathcal{N}^{\Lambda}_{(1,1,1)}(d_1,m,\varpi)$ is given by
  \begin{displaymath}
    H^*(\mathcal{N}^{\Lambda}_{(1,1,1)}(d_1,m,\varpi))^{\Gamma_3}
    \cong H^*(S^{m_1}X \times S^{m_2}X)\ .
  \end{displaymath}
  The variant part of the cohomology is concentrated in degree
  $m_1+m_2$ and has dimension
  \begin{displaymath}
    (3^{2g}-1)\binom{2g-2}{m_1}\binom{2g-2}{m_2}\ .
  \end{displaymath}
\end{prop}

\begin{proof}
  This is essentially \cite[Proposition~3.11]{G1}, cf.\ also
  \cite{HT2}.
\end{proof}

Given this result, we can now find the contribution to the  Poincar\'e
polynomial of $\mathcal{M}^{\Lambda}$ coming from the invariant part
of the critical submanifolds of type $(1,1,1)$.

\begin{prop}\label{prop:invariant-111}
  Under Assumption \ref{ass:small-weights}, the contribution of the
  invariant part of the cohomology of critical submanifolds of type
  $(1,1,1)$ to the Poincar\'e polynomial of $\mathcal{M}^{\Lambda}$ is
  \begin{multline*}
  P^{\Gamma_3}_t(\Lambda,(1,1,1)) =
  \Coeff_{u^0v^0}  \bigg( (1+2u^2vt^2+2uv^2t^2+u^3v^3t^4)^n  \,
  \cdot \\ \cdot \,
  \frac{t^{2(4g-3+n)}
  (1+u^2vt)^{2g}(1+uv^2t)^{2g}}
  {u^{3n+6g-9+\Delta_0}v^{3n+6g-6-\Delta_0}
    (1-u^2v)(1-uv^2)(1-u^2vt^2)(1-uv^2t^2)
    (1-v^3t^2)(1-u^3t^2)} \bigg)\ ,
  \end{multline*}
  where $\Delta_0 \in \{1,2\}$ is the remainder modulo $3$ of
  $\Delta=\deg(\Lambda)$.
\end{prop}

\begin{proof}
  The proof proceeds exactly as in Section~\ref{sec:(1,1,1)}, except
  that we omit the factor $(1+t)^{2g}$ coming from the Jacobian (cf.\
  Propositions~\ref{prop:boundsmm}, \ref{prop-index-fixed} and
  \ref{prop:cohomology-111}).
\end{proof}

It remains to find the contribution from the variant part
of the critical submanifolds of type $(1,1,1)$.

\begin{prop}\label{prop:variant-111}
  Under Assumption \ref{ass:small-weights}, the contribution of the
  variant part of the cohomology of critical submanifolds of type
  $(1,1,1)$ to the Poincar\'e polynomial of $\mathcal{M}^{\Lambda}$ is
  \begin{displaymath}
    P^{\mathrm{var}}_t(\Lambda,(1,1,1)) =
    2\cdot 6^{n-1}(3^{2g}-1)t^{12g-12+6n}(t+1)^{4g-4}\ .
  \end{displaymath}
\end{prop}

\begin{proof}
  It is convenient to parametrize the critical submanifolds by
  $(m_1,m_2,\varpi)$ as explained in
  Remark~\ref{rem:fixed-critical-111-m1m2}.  In order to do the
  calculation, we therefore need to express the Morse index in terms
  of these invariants.  Using Proposition~\ref{prop-index-fixed} and
  Lemma~\ref{lem:index-111} we obtain
  \begin{displaymath}
    \lambda_{(m_1,m_2,\varpi)} = 16g-16+6n-2m_1-2m_2\ .
  \end{displaymath}
  Hence, using Proposition~\ref{prop:cohomology-111}, the contribution
  from the variant part of the cohomology of critical submanifolds of
  type $(1,1,1)$ for fixed $\varpi$ is
  \begin{equation}
    \label{eq:6}
    \sum
    t^{16g-16+6n-m_1-m_2}(3^{2g}-1)\binom{2g-2}{m_1}\binom{2g-2}{m_2}\ ,
  \end{equation}
  where the sum is over $(m_1,m_2)$ satisfying the conditions
  \eqref{eq:13} and \eqref{eq:12}. Note that the terms in the sum are
  only non-zero when $0\leq m_i \leq 2g-2$. But Assumption
  \ref{ass:small-weights} implies that the region defined by
  \eqref{eq:13} contains all such $(m_1,m_2)$. Therefore we can sum
  over all $(m_1,m_2)$, subject to the condition \eqref{eq:12}.  For
  this, let $\xi = e^{2\pi\ima/3}$, then $\sum_{j=1}^3\xi^{j\nu}$
  equals $3$ if $\nu \equiv 0 \pmod 3$ and zero otherwise. It follows
  that we can rewrite \eqref{eq:6} as
  \begin{align*}
    \label{eq:7}
    &\quad\frac{1}{3}\sum_{m_1,m_2}\sum_{j=1}^{3}
      \xi^{j(m_1+2m_2 +\Delta +s_1+2s_2)}
      t^{16g-16+6n-m_1-m_2}(3^{2g}-1)
      \binom{2g-2}{m_1}\binom{2g-2}{m_2} \\
    &=\frac{(3^{2g}-1)t^{16g-16+6n}}{3}\sum_{j=1}^{3}
      \xi^{j(\Delta+s_1+2s_2)}
      (1+t^{-1}\xi^{j})^{2g-2}(1+t^{-1}\xi^{2j})^{2g-2} \\
    &=\frac{(3^{2g}-1)t^{12g-12+6n}}{3}
      \bigl((\xi^{\Delta+s_1+2s_2}+\xi^{2(\Delta+s_1+2s_2)})(t^2-t+1)^{2g-2}
      +(t+1)^{4g-4}\bigr)\ .
    \end{align*}
    It remains to do the sum over $\varpi\in(S_3)^n$.  For this we use
    Table~\ref{tab:si-sigma} to obtain
    \begin{displaymath}
      \sum_{\varpi}\xi^{\Delta+s_1+2s_2}
      =\sum_{\varpi}\xi^{2(\Delta+s_1+2s_2)}=0\ .
    \end{displaymath}
    Since the number of elements of $(S_3)^n$ is $6^n$, we therefore
    obtain the result of the statement of the Proposition.
\end{proof}

\subsection{Parabolic triples of fixed determinant}
\label{sec:triples-fixed}

Now we want to describe the moduli spaces of parabolic triples
with fixed determinant, that we shall use in the following section
to deal with the critical submanifolds of types $(1,2)$ and
$(2,1)$. We follow the notations of Sections
\ref{sec:stable-triples} and \ref{sec:critical-values-flips}.
Fixing the topological and parabolic types of the triples, there
is a determinant map on the moduli space $\cN_\sigma$ of
$\sigma$-stable parabolic triples,
  \begin{eqnarray*}
  \det: \cN_\sigma & \to & \Jac(X) \\
   T=(E_1,E_2,\phi) &\mapsto & \det(E_1) \otimes \det(E_2).
  \end{eqnarray*}
We define the moduli space of $\sigma$-stable parabolic triples
with fixed determinant $\Lambda$ as
 $$
 \cN_\sigma^{\Lambda}={\det}^{-1}(\Lambda).
 $$

In order to state the deformation theory of the parabolic triples
with fixed determinant, we need to introduce the following
subcomplex of $C^{\bullet}(T,T)$,
 $$
  C^{\bullet}_0(T,T):
  \big(\PH(E_{1}, E_{1}) \oplus  \PH(E_{2}, E_{2})\big)_0
  \overset{c}{\too}
  \SPH(E_{2}, E_{1}(D)),
 $$
where $\big(\PH(E_{1}, E_{1}) \oplus  \PH(E_{2}, E_{2})\big)_0 $
is defined as the kernel of the map
 $$
 \begin{array}{ccc}
  \PH(E_{1}, E_{1}) \oplus \PH(E_{2}, E_{2}) &\to& \cO, \\
   (a_1,a_2) &\mapsto &\Tr(a_1) + \Tr(a_2).
 \end{array}
 $$
We have the following result.

\begin{thm}\label{thm:smoothdim-fixed}
Let $T=(E_1,E_2,\phi)$ be a $\sigma$-stable parabolic triple with
determinant $\Lambda$.
\begin{itemize}
 \item[(i)] The Zariski tangent space at the point defined by $T$
 in the moduli space of stable triples with fixed determinant is
 isomorphic to $\HH^{1}(C^{\bullet}_0(T,T))$.
 \item[(ii)] If $\HH^{2}(C^{\bullet}_0(T,T))= 0$, then the moduli
 space of $\sigma$-stable parabolic triples with fixed determinant
 is smooth in a neighbourhood of the point defined by $T$.
 \item[(iii)] If $\phi$ is injective or surjective then
 $T=(E_1,E_2,\phi)$ defines a smooth point in the moduli
 space $\cN_\sigma^\Lambda$.
\end{itemize}
\end{thm}

\begin{proof}
Items (i) and (ii) follow from Theorem \ref{thm:smoothdim}. For
(iii), let us define the complex $C^\bullet_{\rm det}: \cO \to 0$.
This complex is embedded in $C^\bullet (T,T)$ as
 $$
 \begin{array}{ccccc}
C^\bullet_{\rm det} &:& \cO &\longrightarrow & 0 \\
 \downarrow &&\downarrow & & \downarrow \\
  C^\bullet(T,T)&:& \PH(E_{1}, E_{1}) \oplus \PH(E_{2}, E_{2}) &
  \longrightarrow & \SPH(E_{2}, E_{1}(D)),
 \end{array}
 $$
where the left map is $\lambda\mapsto (\lambda\, \Id, \lambda\,
\Id)$. Then it is easy to see that we have a direct sum splitting
of complexes as
 \begin{equation}\label{eqn:dec-fixed}
 C^{\bullet}(T,T)= C^{\bullet}_0(T,T)\oplus
 C^{\bullet}_{\rm det} \, .
 \end{equation}
Now, if $\phi$ is injective or surjective, then the decomposition
(\ref{eqn:dec-fixed}) gives that $0=\HH^2(C^\bullet(T,T))=
\HH^2(C^\bullet_0(T,T))\oplus \HH^2(C^\bullet_{\rm det})$, from
where we get the result stated in (iii).
\end{proof}

In order to study the variation of the moduli spaces
$\cN_\sigma^\Lambda$ when moving $\sigma$, we follow the arguments
of Section \ref{sec:critical-values-flips}. We keep the notations
of that section and work under Assumption \ref{assumption}.
Consider a critical value $\sigma_c$. There is a determinant map
 \begin{eqnarray*}
 \det: B_{\s_c}=\cN_{\s_c}'\times \cN_{\s_c}'' &\to & \Jac(X), \\
  (T',T'') &\mapsto &\det T'\otimes \det T''.
 \end{eqnarray*}
We introduce the following subspace of $B_{\s_c}$,
 \begin{equation}\label{eqn:b-fixed}
 B^{\Lambda}_{\s_c}=
 {\det}^{-1}(\Lambda).
 \end{equation}
The flip loci in the moduli space of parabolic triples with fixed
determinant are given by
 $$
 \cS_{\s_c^{\pm}}^{\Lambda}=
 \cS_{\s_c^{\pm}} \cap \cN_{\s_c^{\pm}}^{\Lambda} \, \subset
 \cN_{\s_c^{\pm}}^{\Lambda}\, .
 $$
The description of $\cS_{\s_c^{\pm}}^{\Lambda}$ follows the arguments
of Subsection \ref{sec:crossing-critical-values}. We get the
following.

\begin{prop}
 \begin{itemize}
 \item[(i)] If $\HH^2(C^\bullet(T',T'))=0$ and
 $\HH^2(C^\bullet(T'',T''))=0$ for every
 $(T',T'')\in B^{\Lambda}_{\s_c}$,
 then
 $B^\Lambda_{\s_c}$ is smooth.
 \item[(ii)] If $\HH^2(C^\bullet(T'',T'))=0$ and
 $\HH^2(C^\bullet(T',T''))=0$ for every $(T',T'')\in B^{\Lambda}_{\s_c}$,
 then $\cS_{\s_c^\pm}= \PP \left(W^\pm|_{B^{\Lambda}_{\s_c}}
 \right)$, where $W^\pm|_{B^{\Lambda}_{\s_c}}$ is
 the restriction of $W^\pm \to B_{\s_c}$ to $B^{\Lambda}_{\s_c}$.
 \end{itemize}
\end{prop}

\begin{proof}
The tangent space to $B^{\Lambda}_{\s_c}$ is given by the
following subcomplex of the complex $C^{\bullet}(T',T') \oplus
C^{\bullet}(T'',T'')$,
 $$
    C^{\bullet}_{d}(T',T''):
    \begin{array}{ll}
    \big( \PH(E_{1}', E_{1}') \oplus  \PH(E_{2}', E_{2}') \oplus \\
    \quad \oplus\PH(E_{1}'', E_{1}'') \oplus  \PH(E_{2}'', E_{2}'')\big)_0
    \end{array}
 \too
    \begin{array}{ll}
    \SPH(E_{2}', E_{1}'(D)) \oplus \\
    \quad \oplus\SPH(E_{2}'', E_{1}''(D))
    \end{array}
 $$
where the $C^0_d(T',T'')$ is the kernel of the map
  $$
  \begin{array}{r}
  \PH(E_{1}', E_{1}') \oplus  \PH(E_{2}', E_{2}')
  \oplus  \PH(E_{1}'', E_{1}'') \oplus  \PH(E_{2}'', E_{2}'')
  \quad \to\quad \cO, \\
 (a_1',a_2',a_1'',a_2'') \qquad\qquad\mapsto \quad
 \Tr(a_1')+\Tr(a_2')+\Tr(a_1'')+\Tr(a_2'').
 \end{array}
 $$
Again there is a splitting of complexes
 $$
 C^{\bullet}(T',T') \oplus C^{\bullet}(T'',T'') = C^\bullet_d(T',T'')
 \oplus C^\bullet_{\rm det},
 $$
where $C^\bullet_{\rm det}\inc C^{\bullet}(T',T') \oplus
C^{\bullet}(T'',T'')$ is given by the map
  \begin{eqnarray*}
  \cO &\to & \PH(E_{1}', E_{1}') \oplus  \PH(E_{2}', E_{2}')
  \oplus  \PH(E_{1}'', E_{1}'') \oplus  \PH(E_{2}'', E_{2}''), \\
\lambda &\mapsto & (\lambda\, \Id,\lambda \, \Id,\lambda\,
\Id,\lambda \, \Id).
 \end{eqnarray*}

This proves that $\HH^2(C^\bullet_d)=0$ and hence that
$B^\Lambda_{\s_c}$ is smooth. The second item follows from
Proposition \ref{prop:5.9}.
\end{proof}

\begin{prop} \label{prop:blow-up-fixed}
Assume that $\cN_{\s_c^\pm}^{\Lambda}$ and $B_{\s_c}^{\Lambda}$
are smooth, and that $\HH^2(C^\bullet(T'',T'))=0$ and
$\HH^2(C^\bullet(T',T''))=0$ for every $(T',T'')\in
B_{\s_c}^{\Lambda}$. Let $\widetilde{\cN}_{\s_c^\pm}^\Lambda$ be
the blow-up of $\cN_{\s_c^\pm}^{\Lambda}$ along
$\cS_{\s_c^\pm}^{\Lambda}$. Then
 $$
 \widetilde{\cN}_{\s_c^+}^\Lambda\cong \widetilde{\cN}_{\s_c^-}^\Lambda\, .
 $$
\end{prop}

\begin{proof}
The proof of Proposition \ref{prop:normal-bundle} needs some
slight modifications to the situation of fixed determinant. The
complex $C^\bullet(\cT,\cT)$ used to compute the tangent bundle to
$\cN_{\sigma^+_{c}}$ should be substituted by
$C^\bullet_0(\cT,\cT)$, which computes the tangent bundle to
$\cN_{\sigma^+_c}^\Lambda$. Likewise, the complex
$C^\bullet(\cT',\cT')\oplus C^\bullet(\cT'',\cT'')$ should be
substituted by the complex $C^\bullet_d(\cT',\cT'')$ introduced
above, which deals with the tangent bundle to $B^\Lambda_{\s_c}$.
Also, the piece
 $$
 \PH_U(\cE_1,\cE_1) \oplus \PH_U(\cE_2,\cE_2)
 $$
in the complex computing the tangent bundle to $\PP W^+$ must be
substituted by the kernel $\big(\PH_U(\cE_1,\cE_1) \oplus
\PH_U(\cE_2,\cE_2)\big)_0$ of
 \begin{eqnarray*}
  \PH_U(\cE_1,\cE_1) \oplus \PH_U(\cE_2,\cE_2) &\to &\cO, \\
  (a_1,a_2 ) &\mapsto &\Tr(a_1)+\Tr(a_2).
 \end{eqnarray*}
Taking this into account, we reach the conclusion that the normal
bundle to $\mathcal{S}_{\s_c^\pm}^\Lambda$ in
$\mathcal{N}_{\s_c^\pm}^\Lambda$ is isomorphic to $\left( p^*
W^\mp \ox \mathcal{O}_{\PP
W^\pm}(-1)\right)|_{B^{\Lambda}_{\s_c}}$.

Using this last fact, the arguments of Proposition \ref{prop:5.12}
carry over verbatim to prove the stated isomorphism.
\end{proof}

Finally we apply Proposition \ref{prop:blow-up-fixed} to compute
the Poincar{\'e} polynomial of the moduli spaces of $\s$-stable
triples with fixed determinant for the case of ranks $r_1=2$ and
$r_2=1$. We follow the notations of Section
\ref{sec:thaddeus-program}. The description of the flip loci also
holds in this situation. Let $\s_c>\s_m$ be a critical value. Then
we have the following equality of Poincar{\'e} polynomials:
 \begin{equation}\label{eqn:Pt-a-fixed}
  P_t\big(\cN_{\s_c^-}^\Lambda\big) -P_t\big(\cN_{\s_c^+}^\Lambda\big) =
  P_t\big( \PP \left(W^-_{\s_c}|_{B^\Lambda_{\s_c}}\right)\big) - P_t\big(\PP
  \left( W^+_{\s_c}|_{B^\Lambda_{\s_c}}\right)\big)\, ,
 \end{equation}
where $W^\pm_{\s_c}|_{B^\Lambda_{\s_c}}$ is a projective fibration
over $B^{\Lambda}_{\s_c}$ with fibres projective spaces of
dimension $w^\pm_{\s_c}-1$. But $B^{\Lambda}_{\s_c}
=\det^{-1}(\Lambda)$ where the determinant map is
  $$
  \begin{array}{rcccl}
  \cN'_{\s_c} \times \cN''_{\s_c}
   &= &\Jac^{d_M} X \times (\Jac^{d_2} X \x
  S^{N}X) & \to &\Jac (X), \\
  && (M, L, Z) &\mapsto & M\otimes L\otimes L(Z).
  \end{array}
  $$
Therefore we have an isomorphism
 $$
 B^{\Lambda}_{\s_c} \cong \Jac\, X\x S^N X.
 $$
The arguments of the proof of Theorem \ref{thm:Pt(Ns)} now give
the following result.

\begin{thm} \label{thm:Pt(Ns)-fixed} Let $\s>\s_m$ be a non-critical value.
For any $\e=\{\e(p)\}_{p\in D}$, $\e(p)\in \{1,2\}$, let
$s_1,s_2,s_3$ and $\bdm$ be defined as in Theorem
\ref{thm:Pt(Ns)}. Then
 $$
 \begin{aligned}
 P_t(\cN_\sigma^{\Lambda})= \sum_{\e} \Coeff_{x^0} \bigg( & \frac{(1+t)^{2g}(1+xt)^{2g}
 t^{2d_1-2d_2+2s_2+2s_3-2\bdm}x^{\bdm}}{(1-t^2)(1-x)(1-xt^2)(1-t^{-2}x)
 x^{d_1-d_2+s_1}} \\
 &-\frac{(1+t)^{2g}(1+xt)^{2g}
 t^{-2d_1+2g-2 +2n-2s_3+4\bdm}x^{\bdm} }{(1-t^2)(1-x)(1-xt^2)(1-t^4x)
 x^{d_1-d_2+s_1}}\bigg)\, .
 \end{aligned}
 $$ \qed
\end{thm}

\subsection{Critical submanifolds of type $(1,2)$ and $(2,1)$}
\label{sec:type-12-21-fixed}

First consider the critical submanifolds of type $(1,2)$. We shall
use the notations of Section~\ref{sec:(1,2)}. Note that the
description given in Section~\ref{sec:(1,2)-description} of the
parabolic Higgs bundles which corresponds to critical points of
type $(1,2)$ remains valid. Thus, fixing $d_1$ and $\varpi$, the
fixed determinant critical submanifold is the fibre of the map
$\det$ defined in \eqref{eq:5}:
 \begin{displaymath}
  \mathcal{N}^{\Lambda}_{(1,2)}(d_1,\varpi) = {\det}^{-1}(\Lambda)\ .
 \end{displaymath}
The characterization of stability given in
Proposition~\ref{prop:august} tells us that
$\mathcal{N}^{\Lambda}_{(1,2)}(d_1,\varpi)$ is isomorphic to the
moduli space of $\sigma$-stable triples (of the appropriate
topological and parabolic type) with fixed determinant
$\Lambda\otimes K^2$, as considered in Subsection
\ref{sec:triples-fixed}, for $\sigma=2g-2$. Therefore Theorem
\ref{thm:Pt(Ns)-fixed} and the computations of
Section~\ref{sec:(1,2)} give the following.

\begin{prop} \label{prop:P_t(1,2)_nonzero-fixed}
Under Assumption \ref{ass:small-weights}, the contribution of the
critical submanifolds of type $(1,2)$ to the Poincar\'e polynomial
of $\cM^{\Lambda}$ is
 \begin{multline*}
  P_t(\Lambda,(1,2)) =\Coeff_{x^0}\bigg( \frac{(1+t)^{2g}(1+xt)^{2g}
 t^{8g-8+2n}x^{5-2g-\Delta_0-n}
 (1+2t^2+2t^{2}x+t^4x)^n}{(1-t^2)(1-x)(1-xt^2)(1-t^{-2}x)(1-t^2x^2)}  \\
 - \frac{(1+t)^{2g}(1+xt)^{2g}
 t^{6g +6 -4\Delta_0} x^{5-2g-\Delta_0-n} (1+ 2 t^{2} +2 t^4x +t^6x )^n}{(1-t^2)
 (1-x)(1-xt^2)(1-t^4x)(1-t^8x^2)} \bigg) \, ,
  \end{multline*}
  where $\Delta_0 \in \{1,2\}$ is the remainder modulo $3$ of
  $\Delta=\deg(\Lambda)$.
\qed
\end{prop}

Now consider the critical submanifolds of type $(2,1)$. Use the
notations of Section~\ref{sec:(2,1)}. Fixing $d_1$ and $\varpi$,
the fixed determinant critical submanifold is the fibre of the map
$\det$ defined in \eqref{eq:5}:
 \begin{displaymath}
  \mathcal{N}^{\Lambda}_{(2,1)}(d_1,\varpi) = {\det}^{-1}(\Lambda)\ .
 \end{displaymath}
Lemma \ref{lem:9.1} holds in this situation, telling us that
$\mathcal{N}^{\Lambda}_{(2,1)}(d_1,\varpi)$ is isomorphic to the
moduli space of $\sigma$-stable triples of type $(1,2)$, with
appropriate degrees and weights, with fixed determinant
$\Lambda^{-1}\otimes K^{-1}(-3D)$, and for $\sigma=2g-2$. The
computations of Section~\ref{sec:(2,1)} together with Theorem
\ref{thm:Pt(Ns)-fixed} yield the following.

\begin{prop} \label{prop:P_t(2,1)_nonzero-fixed}
 Under Assumption \ref{ass:small-weights}, the contribution of the
critical submanifolds of type $(2,1)$ to the Poincar\'e polynomial
of $\cM^{\Lambda}$ is
  \begin{multline*}
  P_t(\Lambda,(2,1)) =\Coeff_{x^0} \bigg( \frac{(1+t)^{2g}(1+xt)^{2g}
 t^{8g-8+2n}x^{2-2g+\Delta_0-n}( 1+2t^2+2t^{2}x+t^4x)^n}{(1-t^2)
 (1-x)(1-xt^2)(1-t^{-2}x)(1-t^2x^2) }  \\
 -  \frac{(1+t)^{2g}(1+xt)^{2g}
 t^{6g-6 +4\Delta_0 } x^{2-2g+\Delta_0- n}(1 + 2t^{2} +2t^4x +t^6x )^n}{(1-t^2)
 (1-x)(1-xt^2)(1-t^4x)(1-t^8x^2)}
 \bigg) \, ,
  \end{multline*}
  where $\Delta_0 \in \{1,2\}$ is the remainder modulo $3$ of
  $\Delta=\deg(\Lambda)$.
\qed
\end{prop}

\begin{rem}\label{rem:cover-12}
  The Poincar\'e polynomials of the critical submanifolds of type
  $(1,2)$ and $(2,1)$ for fixed and non-fixed determinant differ by a
  factor of $(1+t)^{2g}$, coming from the Jacobian. Hence $\Gamma_3$
  acts trivially on the rational cohomology of the fixed determinant
  critical submanifolds of type $(1,2)$ and $(2,1)$ (cf.\
  Remarks~\ref{rem:cover} and \ref{rem:cover-critical}).  The
  triviality of the action can also be seen directly from our
  description of the critical submanifolds as moduli spaces of
  triples, by using the flips picture and arguing as in the proof of
  \cite[Lemma~10.5]{HT2}.
\end{rem}

\subsection{Critical submanifolds of type $(3)$}
\label{sec:type-3-fixed}

The critical points of type $(3)$ are just the stable parabolic
bundles and hence the corresponding critical submanifold is the moduli
space of stable parabolic bundles of fixed determinant $\Lambda$. As
pointed out in Remark~\ref{rem:holla-non-fixed}, the fixed determinant
case was the one studied by Holla \cite{Ho}, and we obtain the
Poincar\'e polynonial of the fixed determinant moduli space by dividing
the formula of Proposition~\ref{prop:19} by $(1+t)^{2g}$.  Thus we
have the following.

\begin{prop}
  \label{prop:type-3-fixed}
  Under Assumption \ref{ass:small-weights}, the Poincar\'e polynomial
  of the moduli space of stable parabolic bundles of rank $3$ of fixed
  determinant $\Lambda$ is given by
\begin{multline*}
    P_{t}(\Lambda,3) = (1+2t^2+2t^4+t^6)^{n-1} \cdot \\
      \cdot \frac{(1+t^3)^{2g}(1+t^5)^{2g}+
       (1+t)^{4g}(1+t^{2}+t^{4})t^{6g-2}
       -(1+t)^{2g}(1+t^3)^{2g}(1+t^2)^2t^{4g-2}}
     {(1-t^2)^3(1-t^4)}\ .
\end{multline*}
\qed
\end{prop}

\begin{rem}\label{rem:cover-3}
  Note, in particular, that $\Gamma_3$ acts trivially on the rational
  cohomology of the moduli space of stable parabolic bundles of rank
  $3$ of fixed determinant (cf.\ Remarks~\ref{rem:cover} and
  \ref{rem:cover-critical}).
\end{rem}

\subsection{Betti numbers of the fixed determinant moduli space}
\label{sec:rank3-fixed}

Finally we put everything together to obtain the Poincar\'e polynomial
of the moduli space of rank three parabolic Higgs bundles of fixed
determinant $\Lambda$.

\begin{thm} \label{thm:end-fixed}
  Let $\mathcal{M}^{\Lambda}$ be the moduli space of rank three
  parabolic Higgs bundles of fixed determinant $\Lambda$ and some
  fixed weights, over a connected, smooth projective complex algebraic
  curve of genus $g$.  If the weights are generic (in the sense that
  there are no properly semistable parabolic Higgs bundles), then the
  Poincar\'e polynomial of $\mathcal{M}^\Lambda$ is given by
\begin{multline*}
  P_t(\mathcal{M}^\Lambda ) =
\Coeff_{u^0v^0}  \bigg( (1+2u^2vt^2+2uv^2t^2+u^3v^3t^4)^n  \,
 \cdot \\ \cdot \,
  \frac{t^{2(4g-3+n)}
  (1+u^2vt)^{2g}(1+uv^2t)^{2g}}
  {u^{3n+6g-8}v^{3n+6g-7}
    (1-u^2v)(1-uv^2)(1-u^2vt^2)(1-uv^2t^2)
    (1-v^3t^2)(1-u^3t^2)} \bigg) \\
 +\Coeff_{x^0} \Bigg( \frac{(1+t)^{2g}(1+xt)^{2g}
 t^{6g-6} x^{2-2g- n}}{(1-t^2)
 (1-x)(1-xt^2)}\cdot \hspace{7cm}\\
 \cdot\bigg(\frac{t^{2g-2+2n}(x+x^2)(1+2t^2+2t^{2}x+t^4x)^n}
   {(1-t^{-2}x)(1-t^2x^2)}
  -\frac{(t^4x+t^8x^2)(1 + 2t^{2} +2t^4x +t^6x )^n}{(1-t^4x)(1-t^8x^2)}
 \bigg)
 \Bigg) \\
 + (1+2t^2+2t^4+t^6)^{n-1} \cdot \hspace{10cm}\\
      \cdot \frac{(1+t^3)^{2g}(1+t^5)^{2g}+
       (1+t)^{4g}(1+t^{2}+t^{4})t^{6g-2}
       -(1+t)^{2g}(1+t^3)^{2g}(1+t^2)^2t^{4g-2}}
     {(1-t^2)^3(1-t^4)} \\
 + 2\cdot 6^{n-1}(3^{2g}-1)t^{12g-12+6n}(t+1)^{4g-4}\ .
\end{multline*}
\qed
\end{thm}

\begin{proof}
  The Theorem follows by an argument analogous to the
  proof of Theorem~\ref{thm:end}, but now using the contributions from
  the fixed determinant critical submanifolds given in
  Propositions~\ref{prop:invariant-111}, \ref{prop:variant-111},
  \ref{prop:P_t(1,2)_nonzero-fixed}, \ref{prop:P_t(2,1)_nonzero-fixed}
  and \ref{prop:type-3-fixed}. Also,
  Proposition~\ref{prop:deg-tensor-fixed} takes the place of
  Proposition~\ref{prop:deg-tensor}.
\end{proof}

\begin{cor}\label{cor:euler-fixed}
  The Euler characteristic of the moduli space of parabolic Higgs
  bundles with fixed determinant $\Lambda$ is
  \begin{displaymath}
    \chi(\mathcal{M}^{\Lambda}) = 0\ .
  \end{displaymath}
\end{cor}

\begin{proof}
  This could be shown by substituting $t=-1$ in the formula of
  Theorem~\ref{thm:end-fixed}.  But it is, in fact, easier to note
  that the Euler characteristic of the moduli space equals the sum of
  the Euler characteristics of the critical submanifolds. Our
  description of these shows that they all have zero Euler
  characteristic. Hence the only potentially non-zero contribution
  comes from the invariant part of the cohomology of the critical
  submanifolds of type $(1,1,1)$, given in
  Proposition~\ref{prop:cohomology-111}.  From MacDonald's formula
  \cite{McD} we have $\chi(S^{m_i}X) = (-1)^{m_i}\binom{2g-2}{m_i}$
  and hence
  \begin{displaymath}
    \chi(\mathcal{M}^{\Lambda}) =
    \sum (-1)^{m_1+m_2}\binom{2g-2}{m_1}\binom{2g-2}{m_2}\ ,
  \end{displaymath}
  where the sum is over all $(m_1,m_2,\varpi)$ satisfying the
  conditions \eqref{eq:13} and \eqref{eq:12}.  This is essentially the
  calculation of the proof of Proposition~\ref{prop:variant-111}, with
  $t$ substituted by $-1$, and gives zero.
\end{proof}

Finally, we can identify the variant part of the
cohomology of $\mathcal{M}^{\Lambda}$ under the action of
$\Gamma_3$---this should be relevant for proving the rank $3$
parabolic version, stated in \cite{HT1}, of the mirror symmetry
Theorem of Hausel--Thaddeus \cite{HT2}.

\begin{thm}
  The variant part of the rational cohomology of
  $\mathcal{M}^{\Lambda}$ has Poincar\'e polynomial
  \begin{displaymath}
    P^{\mathrm{var}}_t(\mathcal{M}^{\Lambda}) =
    2\cdot 6^{n-1}(3^{2g}-1)t^{12g-12+6n}(t+1)^{4g-4}\ .
  \end{displaymath}
\end{thm}

\begin{proof}
  As we have seen in Remarks \ref{rem:cover-12} and \ref{rem:cover-3},
  the critical submanifolds of type $(1,2)$, $(2,1)$ and $(3)$ do not
  contribute to the variant cohomology.  Hence, under
  Assumption~\ref{ass:small-weights}, the variant Poincar\'e
  polynomial $P^{\mathrm{var}}_t(\mathcal{M}^{\Lambda})$ equals the
  contribution coming from critical submanifolds of type $(1,1,1)$,
  given in Proposition~\ref{prop:variant-111}.  But, as we saw in
  \eqref{eq:2},
  \begin{displaymath}
    P_t(\mathcal{M}^{\Lambda})(1+t)^{2g}
      -P_t(\mathcal{M})
    =P_t^{\mathrm{var}}(\mathcal{M}^{\Lambda})(1+t)^{2g}\ ,
  \end{displaymath}
  and we know from Propositions~\ref{prop:deg-tensor} and
  \ref{prop:deg-tensor-fixed} that the left hand side is independent
  of the choice of $\Delta$ and parabolic weights made in
  Assumption~\ref{ass:small-weights}.  Hence the right hand side is
  also independent of this choice.  This finishes the proof.
\end{proof}

\end{document}